\documentclass[11pt]{amsart}
\usepackage[a4paper,margin=24mm]{geometry}
\usepackage{amsmath,amssymb,amsthm,mathtools}
\usepackage{xcolor}
\usepackage{booktabs,tikz}
\usepackage{pgfplots}
\usepackage[section]{placeins}
\AddToHook{cmd/subsection/before}{\FloatBarrier}
\pgfplotsset{compat=1.18}
\usepgfplotslibrary{groupplots}
\usepackage[T1]{fontenc}\usepackage{lmodern,microtype,hyperref}
\usetikzlibrary{arrows.meta,calc}
\allowdisplaybreaks
\numberwithin{equation}{section}
\newtheorem{theorem}{Theorem}[section]
\newtheorem{proposition}[theorem]{Proposition}
\newtheorem{lemma}[theorem]{Lemma}
\newtheorem{corollary}[theorem]{Corollary}
\theoremstyle{definition}\newtheorem{remark}[theorem]{Remark}
\newcommand{\C}{\mathbb C}
\newcommand{\ten}{\otimes_A}
\newcommand{\id}{\mathrm{id}}
\newcommand{\Ric}{\operatorname{Ric}}
\newcommand{\Scal}{\operatorname{Scal}}
\newcommand{\Om}{\Omega}
\newcommand{\extd}{{\rm d}}

\newcommand{\R}{{\mathbb{R}}}

\newcommand{\Z}{{\mathbb{Z}}}

\newcommand{\<}{{\langle}}
\renewcommand{\>}{{\rangle}}

\newcommand{\tens}{\otimes}

\newcommand{\eps}{\epsilon}
\renewcommand{\imath}{\mathrm{i}}

\allowdisplaybreaks
\title{Riemannian geometry and quantum gravity on trees}


\author{Shahn Majid}
\address{School of Mathematical Sciences\\ Queen Mary University of London \\ Mile End Rd, London E1 4NS }
\email{s.majid@qmul.ac.uk}

\begin{document}
\maketitle

\begin{abstract} We show that the tree graphs that admit a quantum Riemannian geometry (QRG) with a positive ratio metric are exactly the Dynkin graphs of finite type $A,D,E$ and affine type $\widetilde D,\widetilde E$. A metric in this context is an assignment of a nonzero real value to each arrow and positive ratio requires that the two directions on any edge have the same sign. It turns out that the metric value in one direction can be chosen freely and the other is then determined. For finite $A, D, E$ type, the construction involves $q$-integers where $q=e^{\imath\pi/h}$ with $h$ the Coxeter number and there is a conjugate pair of $*$-preserving quantum Levi-Civita connections (QLCs). In the case of type $A$, these are necessarily contained in a previously known 1-parameter space of connections.  For affine $\widetilde D,\widetilde E$ type, there is only one QLC,  with real coefficients. Using the quantum Riemannian geometry, we construct a natural Einstein--Hilbert action and hence new quantum gravity models on $A_3,D_4,\widetilde D_4,\widetilde D_5,\widetilde D_6,\widetilde E_6$ in a functional integral approach, other than the 4-star $\widetilde D_4$ which we revisit. \end{abstract}

\section{Introduction}

Quantum or noncommutative Riemannian geometry (QRG) is a more explicit form of noncommutative geometry more generally\cite{Con}  based on algebraic differential forms and bimodule connections\cite{BegMa}. Growing out of the geometry of quantum groups\cite{Dri,Jim,Ma:book}, it has the merit that it works over any field and applies in principle to any algebra $A$ in the role of `coordinate algebra', including noncommutative and finite-dimensional ones. It also applies to graphs as discrete geometries, where the algebra is the commutative one of functions on the set of vertices but differential forms do not commute with functions or graded-commute with each other\cite{Ma:gra}.  A metric is $g\in\Omega^1\tens_A\Omega^1$, where $\Omega^1$ is the $A$-bimodule of differential forms. 

Due to the quadratic nature of the conditions of metric compatibility in the presence of noncommutativity, QRG does not, however, have automatic existence of a `quantum Levi-Civita' connection (QLC) for a given metric. In the graph case a metric reduces to a nonzero real number on every arrow, the `square-length' of the arrow (we consider negative values as of Lorentzian signature) and it becomes an interesting question as to which graphs admit a QRG. Corollary~\ref{cor:univalent}  in the present work, motivated by experience with  the  $A_n$ graph $\bullet$--$\bullet$--$\cdots$--$\bullet$ in \cite{ArgMa:dis} and the $n$-star graph in \cite{BegMa:gra}, states that if a connected graph has more than two vertices and a unary vertex (an external leg or `leaf') then it {\em cannot} admit a QRG with edge-symmetric metric. The metric in one direction cannot equal that in the other direction for all edges if there is a unary vertex.  Tree graphs are obviously in this category and hence a restricted version of the problem is: which tree graphs admit a QRG in this possibly  asymmetric sense. We still want to have the same sign for the two directions of every edge, which we call a {\em positive ratio} metric. Note that we are not restricting the common sign and still allow mixed signature metrics where values in both directions could be  negative on some edges and positive in both directions on other edges. 

Extending the known $A_n$ case, we now fully answer this question for positive ratio metrics: the allowed trees are precisely the finite Dynkin graphs $A_n,D_n,E_6,E_7,E_8$ and the affine Dynkin graphs $\widetilde D_n$, $n\ge4$, and $\widetilde E_6,\widetilde E_7,\widetilde E_8$ as in \cite{Kac}. For $n\ge 5$, $\widetilde D_{n}$ is two trivalent vertices joined by a path of $n-4$ edges while for $n=4$ this contracts to the 4-star. The graph $\widetilde D_5$ is encountered in Physics as a Feynman graph for 2-particle scattering, while $\widetilde D_6$ is the binary tree with 7 vertices. Corollary~\ref{cor:positive-trees} provides this complete classification of underlying trees admitting QRGs with positive ratio metrics. The metric values in each `outward' direction can be arbitrary nonzero reals, with the reverse values then determined. The finite $A, D, E$ types then have a conjugate pair of QLCs,  belonging to the known one-parameter circle of QLCs \cite{ArgMa:dis} in the case of type $A$. By contrast, the  affine $\widetilde D,\widetilde E$ types have a unique QLC with real coefficients. The affine $\widetilde A$ diagrams are polygons rather than trees. Their known geometry, including asymmetric metrics \cite{ArgMa:qua,BSZ}, is recalled in Remark~\ref{rem:poly} for comparison.
There are, however, tree graphs that do not admit a QLC with a positive ratio metric but do admit one with negative ratios on some edges. These have not been classified in general but found case by case up to 7 vertices in Table~\ref{tab:small-trees}. 

The second half of the paper focusses on a key application of QRG to physics, namely to construct quantum gravity models. The idea is that while continuum quantum gravity is, and is likely to remain, intractable, highly noncommutative and/or finite-dimensional $A$ can regularise the infinities in a functional integral approach that appear untameable if our fields are in (a completion of) $C^\infty(M)$ for $M$ a continuum manifold. In the finite-dimensional case, for example, functional integrals now become ordinary finite-dimensional ones. By now, several quantum gravity models have been constructed based on fuzzy spheres\cite{LirMa:fuz} and finite graphs notably the polygon $\Z_n$\cite{Ma:par, ArgMa:qua} and $\Z_2\times\Z_2$\cite{Ma:squ,BliMa:cur} and more recently a first attempt at the 4-star graph (which we revisit here)\cite{Ma:alg4}. We can now add to this repertoire any of the above graphs. We particularly look in most depth at quantum gravity on the binary tree $T_3=\widetilde D_6$, see Section~\ref{sec:quan}. The first step is to compute the Riemann and Ricci curvatures of the QLC, its scalar curvature and hence a suitable Einstein-Hilbert action. For the Ricci tensor we need to `lift' a 2-form into $\Omega^1\tens_A\Omega^1$ for which in our case of $\Omega_{\min}$ there is a natural choice. The second ingredient for this action is integration over the graph, which amounts to a choice of positive weights $\mu$ at every vertex. However, in all cases we can find a choice of $\mu$ that is linear in the metric values $a$.   Aside from the binary tree, we cover in Section~\ref{sec:star} $n$-stars for $n=2,3,4$ where $n=2$ is also $A_3$, $n=3$ is $D_4$ and $n=4$ is $\widetilde D_4$ and is now revisited. This is followed in Section~\ref{sec:join} by adjacent 3-stars $\widetilde D_5$ and in Section~\ref{sec:exte} by the extended 3-star $\widetilde E_6$. 

In all cases we compute the expectation values $\<a\>$ and their relative standard-deviation $\Delta a/\<a\>$, the expectation $\<\Scal\>$ of the Ricci scalar curvature and its standard deviation $\Delta \Scal$. Features recurring across these and earlier models include persistent metric quantum uncertainty, measure dependence, cutoff-sensitive curvature fluctuations and interference effects in Lorentzian action models. We are also starting to see some physical relevance such as to particle physics via Kaluza-Klein theory\cite{LiuMa:kal} and to the  `problem of the cosmological constant'\cite{BliMa:cur}. The latter follows ideas going back to J.A. Wheeler\cite{Whe} that $|\Delta \Scal|$ should be seen as a measure of the vacuum energy of Planck scale spacetime fluctuations and ideas of Carlip, Wang and Unruh\cite{Car,WanUnr} for its impact on local Friedmann equations.  The paper starts with a preliminaries Section~\ref{sec:pre} that provides a short account of the formalism and fixes notation. 

\subsection*{Acknowledgements} This work was supported by a Leverhulme Trust project grant RPG-2024-177 and by access provided by OpenAI through its ChatGPT for Academic Researchers programme.

\noindent{\bf AI declaration:}  ChatGPT(6.1 sol) was used as an extremely competent research assistant for the preparation of proofs and for advanced numerical integrations. The author takes final responsibility for all proofs and reported calculations. 

\noindent{\bf Data declaration:} Experimental data was not used or created in the preparation of this work. 

\noindent{\bf Conflicts declaration:} The author is not aware of any conflicts of interest.

\section{Preliminaries and notation} \label{sec:pre}

Let \(G=(X,E)\) be a finite simple graph, regarded as bidirected, and
\(A=\C(X)\). Write \(x\to y\), abbreviated \(xy\), for a directed
edge. Let \(N(x)\) denote the neighbours of \(x\), and
\(\delta_x\) the characteristic function of the vertex \(x\).
Bidirectedness means that \(x\to y\) occurs if and only if
\(y\to x\) occurs. The one-forms have basis \(\omega_{xy}\), one for
each directed edge \(x\to y\), with
\[
 f\omega_{xy}=f(x)\omega_{xy},\quad
 \omega_{xy}f=f(y)\omega_{xy},\quad
 \extd f=\sum_{x\to y}(f(y)-f(x))\omega_{xy}.
\]
Sums over \(x\to y\) run over all directed edges, including both orientations of each undirected edge.
Write \(\theta=\sum_{x\to y}\omega_{xy}\), so \(\extd f=[\theta,f]\).
The involution is \(\omega_{xy}^*=-\omega_{yx}\). Composable two-paths
\(\omega_{xyz}=\omega_{xy}\ten\omega_{yz}\) form a basis of
\(\Om^1\ten\Om^1\). Our higher calculus is
\[
 \Om_{\min}:\qquad
 \sum_{y\in N(x)\cap N(z)}\omega_{xy}\wedge\omega_{yz}=0
 \quad\text{for every }x,z,
\]
including \(x=z\), with the differential ideal generated by these
relations. It is inner: \(\extd\omega=\theta\wedge\omega+\omega\wedge\theta\).
A real nondegenerate directed metric and its inverse are
\begin{equation}\label{eq:metric-general}
 g=\sum_{x\to y}g_{xy}\omega_{xyx},\qquad
 g_{xy}\in\R\setminus\{0\},\qquad
 (\omega_{xy},\omega_{yx})=g_{yx}^{-1}\delta_x.
\end{equation}
We do not impose edge symmetry or \(\wedge g=0\), but we will be interested in the \emph{positive ratio metrics} where
\(g_{xy}/g_{yx}>0\) on every edge. The common sign of $g_{xy}$ and $g_{yx}$ in this case can still vary from edge to edge allowing for `mixed signature' in a geometric sense.

A bimodule connection has the inner presentation
\[
 \nabla\omega=\theta\ten\omega-\sigma(\omega\ten\theta)+\alpha(\omega),
\]
where \(\sigma\) and \(\alpha\) are bimodule maps. A quantum
Levi--Civita connection (QLC) obeys
\[
 \wedge\nabla=\extd,\qquad
 \bigl(\nabla\otimes\id+(\sigma\otimes\id)(\id\otimes\nabla)\bigr)g=0.
\]
We require star preservation,
\(\nabla(\omega^*)=\sigma\dagger\nabla\omega\), where
\(\dagger(\xi\ten\eta)=\eta^*\ten\xi^*\).
A graph with such a metric and a star-preserving QLC will be called
a quantum Riemannian geometry (QRG). For the real metrics above,
star preservation makes hermitian metric compatibility equivalent
to the displayed tensor metric compatibility
\cite[Section 2.1 and Proposition 3.3]{BegMa:gra}; we use the
hermitian parametrisation below. In the inner presentation,
star preservation requires
\begin{equation}\label{eq:star-pres}
 (\dagger\sigma)^2=\id,\qquad
 \sigma\dagger\alpha(\omega^*)=\alpha(\omega).
\end{equation}
On a tree there are no two-paths with the endpoints of an edge;
therefore \(\alpha=0\) for every bimodule connection.

Since \(\sigma\) is a bimodule map, it preserves the initial and
terminal vertices of a two-path. Its matrix therefore splits into
blocks, one on the span of all \(\omega_{xyz}\) with fixed endpoints
\(x,z\). A \emph{return block} has \(z=x\): its paths
\(x\to y\to x\) return to their starting vertex. Its dimension is
\(\deg(x)\). A \emph{non-return block} has distinct endpoints.
On a tree, each nonempty non-return block contains just one two-path,
since two different intermediate vertices would give a four-cycle.
We call a block with just one basis path a \emph{singleton block};
\(\sigma\) acts on it by a scalar.

For the graph \(G\), write
\begin{equation}\label{PhH}
 P_x=\sum_{y\in N(x)}\omega_{xyx},\qquad
 h_{xy}=\frac{g_{xy}}{g_{yx}},\qquad H_x=\sum_{y\in N(x)}h_{xy}.
\end{equation}
For the local parametrisation of \(\sigma\), introduce complex numbers
\(Q_{xy}\) on the directed edges, with
\(Q_{yx}=\overline{Q_{xy}}\). We then define
\[
 r_{xy}=1+\frac{Q_{xy}}{g_{yx}},\qquad
 r_x=(r_{xy})_{y\in N(x)},\qquad S_x=-I+\boldsymbol1r_x^T.
\]
Here \(r_x\) is a column vector for a fixed ordering of \(N(x)\),
and \(\boldsymbol1\) is the column vector of \(\deg(x)\) ones.
Each neighbour is counted once in these sums. At this point, 
\(S_x\) is just a matrix defined by \(Q_{xy}\) but the following lemma identifies it as the local return block of $\sigma$. 

\begin{lemma}[Local equations on a graph]\label{lem:local}
Let \(G\) have a real nondegenerate directed metric
\eqref{eq:metric-general}.

\textup{(i)} Every star-preserving QLC for \(\Om_{\min}\) has
complex parameters \(Q_{yx}=\overline{Q_{xy}}\) such that
\begin{align}
 \sigma(\omega_{xyx})&=-\omega_{xyx}+r_{xy}P_x,\label{eq:return-sigma}\\
 \sigma(\omega_{xyz})&=\frac{Q_{zy}+g_{zy}}{g_{yx}}\omega_{xyz}
 \quad(x\ne z,\ N(x)\cap N(z)=\{y\}),\label{eq:straight-sigma}
\end{align}
with
\begin{equation}\label{eq:local-star}
 \overline{S_x}S_x=I,\qquad
 r_{xy}\overline{r_{xz}}=h_{xy}h_{xz}
 \quad(y,z\in N(x),\ y\ne z).
\end{equation}
The first equation holds at every vertex $x$, while the second holds whenever
\(N(y)\cap N(z)=\{x\}\).

\textup{(ii)} If \(x\) has degree at least three and
\(N(y)\cap N(z)=\{x\}\) for every pair of distinct neighbours
\(y,z\), then there is a phase \(z_x\) such that
\begin{equation}\label{eq:phase}
 r_{xy}=h_{xy}z_x,\quad |z_x|=1,\quad
 H_x=z_x+\overline{z_x},\quad
 Q_{xy}=g_{xy}z_x-g_{yx}.
\end{equation}
If \(x\) meets a leaf \(l\), then \(g_{lx}/g_{xl}=H_x\).
For an edge \(xy\) whose endpoints both satisfy these hypotheses,
either both phases are \(-1\), or
\begin{equation}\label{eq:branch-edge}
 h_{xy}=1,\qquad z_y=\overline{z_x}.
\end{equation}

\textup{(iii)} If \(G\) is a tree, its star-preserving QLCs are inner (i.e. $\alpha=0$) and exactly given by
 \eqref{eq:return-sigma}--\eqref{eq:straight-sigma} for parameters \(Q_{yx}=\overline{Q_{xy}}\)
obeying both equations in \eqref{eq:local-star} at every vertex.
\end{lemma}
\begin{proof}
\textup{(i)} Equations \eqref{eq:return-sigma}--\eqref{eq:straight-sigma}
are the return-block and singleton-block cases of the torsion and
hermitian-metric parametrisation in
\cite[Propositions 3.1--3.3]{BegMa:gra}. They hold on any graph
admitting a star-preserving QLC, irrespective of its other blocks
or of \(\alpha\).
In \eqref{eq:return-sigma}, \(r_{xy}\) is the coefficient of \(P_x\)
in \((\sigma+\id)(\omega_{xyx})\). Thus the return block of
\(\sigma\), ordered by \(N(x)\), is precisely \(S_x\).
Path reversal fixes each return path, so the first condition in
\eqref{eq:star-pres} on this block is \(\overline{S_x}S_x=I\),
with entrywise conjugation. Expanding \(S_x=-I+\boldsymbol1r_x^T\)
gives
\begin{equation}\label{eq:return-star-scalar}
 \overline{S_x}S_x=I
 \quad\Longleftrightarrow\quad
 \overline{r_{xy}}\left(\sum_{u\in N(x)}r_{xu}-1\right)=r_{xy}
 \quad(y\in N(x)).
\end{equation}
On the reversed singleton paths \(y\to x\to z\) and
\(z\to x\to y\), the two scalar coefficients from
\eqref{eq:straight-sigma}, with one complex conjugated, must
multiply to one. Substitution gives the second of
\eqref{eq:local-star}.

\textup{(ii)} The hypothesis makes every non-return block through
\(x\) a singleton, so the second of \eqref{eq:local-star} holds
for every pair of distinct neighbours. For three or more neighbours,
\(u_y=r_{xy}/h_{xy}\) satisfies \(u_y\bar u_z=1\) for \(y\ne z\).
Thus all \(u_y\) equal one phase; the return equation then gives
\(H_x=2\Re z_x\), by \eqref{eq:return-star-scalar}.
At a leaf the return scalar is
\(\bar z_x-1/h_{xl}\); unit modulus gives \(1/h_{xl}=H_x\).
Finally hermiticity of \(Q\) gives
\(h_{xy}(1+z_x)=1+\bar z_y\). If \(z_x\ne-1\), taking moduli
shows \(h_{xy}(h_{xy}-1)|1+z_x|^2=0\), proving
\eqref{eq:branch-edge}. Otherwise \(z_y=-1\) as well.

\textup{(iii)} On a tree all nonempty non-return blocks are singletons,
and \(\alpha=0\) because there are no triangles. Thus (i) gives
necessity. Conversely, trees have no nondegenerate squares, so the
parametrisation imposes no further metric constraints and gives
torsion-freeness and hermitian metric compatibility.
Equations~\eqref{eq:local-star} give the first condition in
\eqref{eq:star-pres} on every block, and the second is automatic
for \(\alpha=0\). Star preservation then gives tensor metric
compatibility, as recalled above.
\end{proof}

For the quantum gravity part of the paper we also need a bimodule
section \(\iota:\Om^2_{\min}\to\Om^1\ten\Om^1\) of the wedge map.
On a tree, all non-return two-forms vanish, so it suffices to define
\(\iota\) on return two-forms. At a vertex \(x\) with
\(d_x=|N(x)|>0\), a canonical choice is
\begin{equation}\label{eq:uniform-lift}
 \iota_0(\omega_{xy}\wedge\omega_{yx})=\omega_{xyx}-\frac1{d_x}P_x.
\end{equation}
This is zero at a leaf, as required. The images sum to zero and
\(\wedge\iota_0=\id\), so this is a well-defined bimodule section.
It is the zero-sum representative of a return block and is natural
under graph relabelling. It also satisfies
\(\iota_0(\xi^*)=-\dagger\iota_0(\xi)\), the sign coming from the
graded star on two-forms. It is not the only natural choice.
Every bimodule section is of the form
\begin{equation}\label{eq:general-lift}
 \iota_w(\omega_{xy}\wedge\omega_{yx})=\omega_{xyx}-w_{xy}P_x,\qquad
 \sum_{y\in N(x)}w_{xy}=1.
\end{equation}
Real \(w_{xy}\) give the same star compatibility.

\section{Quantum Riemannian geometry on trees}\label{sec:tree}

This section contains our general results regarding QRG on graphs, starting with a result about pendant vertices and then proceeding to the classification of tree graphs that admit QRGs with positive ratio metrics. 

\subsection{Pendant vertices and metric asymmetry}\label{sec:pend}
Here \(G\) is any finite bidirected simple graph, with real
nondegenerate directed metric coefficients. Write $h_{xy}$ and $H_x$ as in (\ref{PhH}) and for 
a leaf \(l\) attached to \(x\), define its inward metric multiple as
\(\rho=g_{lx}/g_{xl}\).

\begin{theorem}[Pendant-vertex constraint]\label{thm:pendant-bound}
Suppose that \(G\) admits a star-preserving QLC for
\(\Om_{\min}\), and that \(l\) is a leaf whose neighbour \(x\)
has degree at least two. Then
\[
 \rho=H_x=\frac1\rho+\sum_{y\in N(x)\setminus\{l\}}
                             \frac{g_{xy}}{g_{yx}}.
\]
In particular \(\rho>1\) if all metric ratios at \(x\) are positive.
If \(x\) is trivalent, with other neighbours \(y,z\) satisfying
\(N(y)\cap N(z)=\{x\}\), then \(|\rho|\le2\).
\end{theorem}
\begin{proof}
Use the parameters \(Q_{xy},r_{xy},S_x\) of
Lemma~\ref{lem:local}\textup{(i)}, whose return-block equations
hold at both \(l\) and \(x\).

Put \(a=h_{xl}=1/\rho\). The leaf return scalar is
\(w=Q_{lx}/g_{xl}\), with \(|w|=1\), so hermiticity gives
\(u:=r_{xl}=1+a\bar w\). For every \(y\ne l\) adjacent to
\(x\), the paths \(l\to x\to y\) and \(y\to x\to l\)
are the unique two-paths with their respective endpoints.
The second equation in \eqref{eq:local-star}, applied to this reversed pair, gives
\[
 u\overline{r_{xy}}=a h_{xy}.
\]
Hence \(u\ne0\) and \(r_{xy}=a h_{xy}/\bar u\).
Equation~\eqref{eq:return-star-scalar} for the leaf column gives
\(\bar u(\sum_{y\in N(x)}r_{xy}-1)=u\). Substituting yields
\[
 |u|^2-u-\bar u+a(H_x-a)=0.
\]
Since \(|u-1|^2=a^2\), this reduces to \(aH_x=1\).
Positivity of all the ratios then gives \(\rho>1/\rho\), hence
\(\rho>1\).

Under the trivalent hypothesis, every pair of distinct neighbours
of \(x\) has \(x\) as its only common neighbour: this follows
from \(l\) being a leaf for pairs involving \(l\), and from the
assumption on \(y,z\) for the remaining pair.
Lemma~\ref{lem:local}\textup{(ii)} therefore gives
\(H_x=z_x+\bar z_x\) with \(|z_x|=1\). Since \(\rho=H_x\),
this proves \(|\rho|\le2\). These constraints concern \(\sigma\);
allowing nonzero \(\alpha\) does not change them.
\end{proof}

\begin{corollary}[Univalent-vertex obstruction]\label{cor:univalent}
A connected graph with more than two vertices and a nondegenerate
edge-symmetric real metric admits no star-preserving QLC for
\(\Om_{\min}\) if it has a univalent vertex.
\end{corollary}
\begin{proof}
Edge symmetry gives \(\rho=1\) and \(H_x=\deg(x)\).
Theorem~\ref{thm:pendant-bound} excludes \(\deg(x)\ge2\), while
\(\deg(x)=1\) would make the connected graph \(K_2\).
\end{proof}

The trivalent upper bound applies to every tree, and to graphs whose other
branches continue arbitrarily provided they do not form a four-cycle
through \(x\). For a tree with a positive ratio metric the Perron argument in
Lemma~\ref{lem:quotient-independent} further gives
\(\sqrt3\le\rho\le2\) at a trivalent vertex attached to a leaf:
the tree contains the three-star, whose Perron eigenvalue is
\(\sqrt3\), and \(\rho\) is the common eigenvalue.

\subsection{Non-existence of QRGs by direct calculation}

Let \(T_k\) be the perfect rooted binary tree with levels
\(0,\ldots,k-1\), and let \(\widehat T_k\) be obtained by adjoining
one extra leaf at its root. Here \(T_2=A_3\) has the known one-phase family of QRGs 
\cite{ArgMa:dis,BegMa:gra}, while
 \(\widehat T_2\) is the three-star and has a QRG in \cite[Theorem 3.4]{BegMa:gra}. We first exclude larger binary trees and their augmented variants, then give the metric obstruction and constructions for other tree families. The seven-vertex case \(T_3\) is constructed in Section~\ref{sec:posi} and studied in Section~\ref{sec:seve}. 

\begin{theorem}\label{thm:no-go}
For \(\Om_{\min}\), no nondegenerate real directed metric on
\(T_k\), \(k\ge4\), or on \(\widehat T_k\), \(k\ge3\), admits a
star-preserving QLC. 
\end{theorem}
\begin{proof}
The degree-three vertices form connected components joined by branch
edges. If one phase is \(-1\), all phases in its component are \(-1\).
At a bottom branch its two leaf ratios \(h\) are \(-1/2\), so its
ratio toward its parent is \(-1\). At its parent the two ratios toward
bottom branches are then \(-1\), forcing the remaining ratio to be
zero since \(H=-2\). This excludes such components whenever the
parent has no leaf. In \(\widehat T_3\), the parent is the root:
its two branch ratios \(-1\) and its external-leaf ratio \(-1/2\)
sum to \(-5/2\), again impossible.

Otherwise branch edges have ratio \(1\), and phases alternate by
conjugation with common \(H\). A bottom branch gives
\(H=1+2/H\), so \(H=2\) or \(-1\). For \(H=2\), its parent's two
branch ratios force its third ratio to be zero. For \(H=-1\), the
third ratio must be \(-3\), contradicting \eqref{eq:branch-edge} if
that neighbour is another branch vertex. This disposes of \(T_k\)
for \(k\ge5\), and \(\widehat T_k\) for \(k\ge4\).

For \(T_4\), both root branches have ratio toward the root \(-3\),
so the root has \(H=-2/3\), whereas the branching vertices have
\(H=-1\). This contradicts the constancy of \(H_x\) proved in
Lemma~\ref{lem:quotient-independent} below.
For \(\widehat T_3\), the root instead requires
\(H=2+1/H\), incompatible with \(H=1+2/H\).
\end{proof}

\subsection{An obstruction independent of the higher calculus}
For positive ratio metrics, the obstruction in Theorem~\ref{thm:no-go}
persists without imposing torsion-freeness. Changing the exterior algebra
therefore cannot remove it while keeping the same first-order calculus,
metric compatibility and star preservation.

\begin{lemma}\label{lem:quotient-independent}
Let a finite tree with a real nondegenerate directed metric carry a
star-preserving metric-compatible bimodule connection, without any
torsion assumption. Write \(B_x\) for its return block,
\(h_{xy}=g_{xy}/g_{yx}\), and \(H_x=\sum_{y\in N(x)}h_{xy}\). Then \(H_x\) is
one common real number \(\lambda\), and
\begin{equation}\label{eq:diagonal-bound}
 |B_x[y,y]|^2=1-h_{xy}(\lambda-h_{xy})\ge0.
\end{equation}
If all \(h_{xy}>0\), there is a positive vertex vector \(v\) with
\(h_{xy}=v_y/v_x\), and \(A_Gv=\lambda v\). Thus \(\lambda\) is the
Perron eigenvalue of the adjacency matrix.
\end{lemma}
\begin{proof}
Write \(B_x[z,y]\) for the coefficient of \(\omega_{xzx}\)
in \(\sigma(\omega_{xyx})\). Hermitian metric compatibility and
star preservation give, independently of higher forms,
\[
 B_x[y,y]=h_{xy}\overline{B_y[x,x]},\quad
 \overline{B_x[y,z]}B_x[z,y]=h_{xy}h_{xz}\ (y\ne z),\quad
 \overline B_xB_x=I.
\]
For clarity, let \(c_{yxz}\) be the scalar on the singleton block
\(\omega_{yxz}\), where \(y\ne z\). The hermitian parametrisation
in \cite[Proposition 3.1]{BegMa:gra} gives
\(g_{zx}\overline{B_x[y,z]}=g_{xy}c_{yxz}\), and reversal gives
\(c_{yxz}\overline{c_{zxy}}=1\). Applying these twice gives the
second identity above; the diagonal hermitian identity gives the first.
Taking the diagonal of the last equation gives
\(D_{xy}:=|B_x[y,y]|^2=1-h_{xy}(H_x-h_{xy})\).
Since \(D_{xy}=h_{xy}^2D_{yx}\), subtraction gives
\(h_{xy}(H_y-H_x)=0\), proving constancy.
On a tree define \(v\) recursively along edges by \(v_y=h_{xy}v_x\).
Positivity and \(\sum_{y\in N(x)}v_y/v_x=\lambda\) prove the final assertion.
\end{proof}

\begin{corollary}\label{cor:quotient-no-rescue}
No further quotient of \(\Om_{\min}\) can produce a star-preserving
QLC for a positive ratio metric on \(T_k\), \(k\ge4\), or
\(\widehat T_k\), \(k\ge3\), if the first-order calculus is unchanged.
\end{corollary}
\begin{proof}
At any bottom branch with two leaf children, the eigenvector equation
gives leaf ratios \(1/\lambda\) and parent ratio
\(h=\lambda-2/\lambda\). Hence \eqref{eq:diagonal-bound} on the
parent entry becomes
\[
 |B_x[\text{parent},\text{parent}]|^2
       =\frac{4-\lambda^2}{\lambda^2}.
\]
For \(T_k\), symmetry and positivity reduce the adjacency problem to
a path on its \(k\) levels with off-diagonal entries \(\sqrt2\), giving
\[
 \lambda=2\sqrt2\cos\frac{\pi}{k+1}>2\quad(k\ge4).
\]
For \(\widehat T_3\), eliminating the leaves gives
\(\lambda^4-5\lambda^2+2=0\), with
\(\lambda^2=(5+\sqrt{17})/2>4\). Larger augmented trees contain this
as a subgraph, so their Perron eigenvalue is at least this value.
In each case the supposed squared modulus is negative.
\end{proof}
The eigenvector condition also explains the depth-dependent metric
ratios suggested by the small examples. For an edge \(x\to y\) from
level \(i\) to level \(i+1\), define its reverse-to-forward factor by
\[
 \phi_i=\frac{g_{yx}}{g_{xy}},\qquad i=0,\ldots,k-2.
\]
Here \(i\) denotes depth, with the root at level zero. The positive
adjacency eigenvector is constant on each level; writing its value
there as \(v_i\), Lemma~\ref{lem:quotient-independent} requires
\(\phi_i=v_i/v_{i+1}\). We use symmetric \(q\)-integers
\[
 (m)_q=\frac{q^m-q^{-m}}{q-q^{-1}},\qquad
 (m)_1=m,\qquad (m+1)_q+(m-1)_q=(q+q^{-1})(m)_q,
\]
where the value at \(q=1\) is by continuity. Up to an overall
positive constant, the level values are
\[
 v_i=2^{-i/2}(i+1)_q,\qquad q=e^{\imath{}\pi/(k+1)}.
\]
Consequently any positive ratio metric admitting such a connection
would have to satisfy
\[
 g_{yx}=\phi_i g_{xy},\qquad
 \phi_i=\sqrt2\,\frac{(i+1)_q}{(i+2)_q}.
\]
These necessary ratios give \(\sqrt2\) on the single level of edges
for \(k=2\), and \(1,2\) on the two levels for \(k=3\).
They are not solutions for larger trees: for \(k\ge4\), the remaining
star and metric equations violate \eqref{eq:diagonal-bound}.
The obstruction applies before imposing torsion and is unchanged by
further degree-two quotients,  but does not apply for metrics with negative ratios.

\subsection{Trees with at most one branching vertex}
\label{sec:path}
Let \(G\) be a finite tree with at least three vertices and at most
one vertex of valency greater than two. Choose
\(c\) to be that vertex if present, and otherwise a middle vertex
of the path. For an even path either of the two middle vertices
may be used. Write \(A_G\) for the adjacency matrix. As in Lemma~\ref{lem:local},
\(h_{xy}=g_{xy}/g_{yx}\) for a real nondegenerate directed metric.

\begin{theorem}\label{thm:one-branch}
The metric admits a star-preserving QLC for \(\Om_{\min}\) if and
only if there are \(\lambda\in[-2,2]\) and a real adjacency
eigenvector \(v\), with no zero entries, such that
\[
 A_Gv=\lambda v,\qquad h_{xy}=\frac{v_y}{v_x}.
\]
For each such metric, the conditions
\[
 r_{cy}=h_{cy}z\quad(y\in N(c)),\qquad
 |z|=1,\qquad z+\bar z=\lambda
\]
select exactly two complex-conjugate QLCs when \(|\lambda|<2\),
and one when \(|\lambda|=2\). If \(G\) has a branching vertex,
these exhaust its star-preserving QLCs.
\end{theorem}
\begin{proof}
Lemma~\ref{lem:quotient-independent} gives \(H_x=\lambda\) at every
vertex. Starting from \(v_c=1\), define \(v_y=h_{xy}v_x\) along
each arm. This gives a nonzero real value at every vertex and
\(A_Gv=\lambda v\), also without sign restrictions.
If \(c\) has degree at least three, \eqref{eq:phase} gives
\(\lambda=z+\bar z\), \(|z|=1\), so \(|\lambda|\le2\).
For a path, the same bound follows from the eigenvector equation
at a vertex where \(|v_x|\) is maximal, since its degree is at most two.

Conversely, choose a unit phase \(z\) with \(z+\bar z=\lambda\).
At \(c\), set
\[
 r_{cy}=h_{cy}z,\qquad Q_{cy}=g_{cy}z-g_{yc},\qquad
 Q_{yc}=\overline{Q_{cy}}.
\]
Since \(H_c=\lambda=z+\bar z\), these satisfy the return and
singleton-path star equations at \(c\), for any valency of \(c\).
The ratios also satisfy
\(|r_{cy}-1|^2=1-h_{cy}(\lambda-h_{cy})\).
Proceed outward along each arm. At a degree-two vertex \(x\), let
\(p\) be its neighbour toward \(c\) and \(y\) its other neighbour.
The incoming \(Q_{xp}\), and hence \(r_{xp}\), is already known.
Define
\[
 r_{xy}=\frac{h_{xp}h_{xy}}{\overline{r_{xp}}},\qquad
 Q_{xy}=g_{yx}(r_{xy}-1),\qquad Q_{yx}=\overline{Q_{xy}}.
\]
The identity
\begin{equation}\label{eq:arm-propagation}
 |r_{xy}-1|^2=1-h_{xy}(\lambda-h_{xy})
\end{equation}
holds at \(c\). It passes across an edge because
\(r_{yx}-1=h_{yx}\overline{r_{xy}-1}\).
At a degree-two vertex put \(a=h_{xp}\), \(b=h_{xy}\),
and \(u=r_{xp}\). Since \(a+b=\lambda\), the incoming identity is
\[
 |u|^2-2\Re u+ab=0.
\]
In particular \(u\ne0\). The outgoing value \(ab/\bar u\)
satisfies the same identity and has product \(ab\) with \(\bar u\).
These identities give both star equations
\eqref{eq:local-star} at \(x\), and continue
\eqref{eq:arm-propagation} along the arm.
At a leaf \(x\), its unique ratio is \(h_{xp}=\lambda\), so
\eqref{eq:arm-propagation} becomes \(|r_{xp}-1|=1\), the
one-dimensional return star equation.

Use \eqref{eq:return-sigma}--\eqref{eq:straight-sigma} to define
\(\sigma\), with \(\alpha=0\). Torsion-freeness follows from the
path-sum relations. To check metric compatibility directly, expand
\[
 \nabla g=\theta\ten g-
 (\sigma\ten\id)(\id\ten\sigma)(g\ten\theta).
\]
On a three-path \(x\to y\to x\to w\), its coefficient is
\[
 g_{wx}\left[r_{xw}-\overline{r_{xw}}
                    \left(\sum_{u\in N(x)}r_{xu}-1\right)\right],
\]
which vanishes by \eqref{eq:return-star-scalar}. On
\(x\to y\to z\to y\), \(z\ne x\), the coefficient is
\[
 \frac{g_{zy}}{h_{yx}}
       \left(h_{yx}h_{yz}-r_{yx}\overline{r_{yz}}\right),
\]
which vanishes by \eqref{eq:local-star} at \(y\).
Every other three-path on a tree has four distinct vertices, and
its coefficient is zero. Thus the construction is
metric-compatible and star-preserving.
Finally \(z\) determines every \(Q\) by the propagation above,
so its two choices give precisely the connections satisfying the
stated condition at \(c\). If \(c\) has degree at least three, that condition
is forced by \eqref{eq:phase}.
\end{proof}

\begin{corollary}\label{cor:qlc-count}
Fix a nondegenerate real directed metric on a finite tree with a
vertex of degree at least three, and suppose it admits a
star-preserving QLC for \(\Om_{\min}\). With \(\lambda=H_x\) as in
Lemma~\ref{lem:quotient-independent}, there are exactly two
complex-conjugate star-preserving QLCs if \(|\lambda|<2\), and
exactly one, with real coefficients in the arrow basis, if
\(|\lambda|=2\).
\end{corollary}
\begin{proof}
Choose a branching vertex \(c\). By Lemma~\ref{lem:local}, its
phase obeys \(|z_c|=1\), \(z_c+\bar z_c=\lambda\), and fixes all
incident \(Q\)'s. Proceed outward through the tree. At a vertex
\(x\) with predecessor \(p\), hermiticity determines \(Q_{xp}\)
and hence \(r_{xp}\). If \(x\) is not a leaf,
the second of \eqref{eq:local-star} forces
\[
 r_{xy}=\frac{h_{xp}h_{xy}}{\overline{r_{xp}}}
 \quad(y\in N(x)\setminus\{p\}),
\]
where \(r_{xp}\ne0\) since the metric is nondegenerate.
Thus a choice of \(z_c\) determines the entire connection.
Complex conjugation preserves all the equations and exchanges the
two distinct phases when \(|\lambda|<2\). When \(|\lambda|=2\),
the initial phase is real and the propagation gives real
coefficients throughout.
\end{proof}

In particular, an \(n\)-star with \(n>4\) admits no
star-preserving QLC for \(\Om_{\min}\), for any nondegenerate real
directed metric, allowing negative ratios as well. Indeed, the leaf and centre
eigenvector equations give \(\lambda^2=n\), whereas
Theorem~\ref{thm:one-branch} requires \(|\lambda|\le2\).
This recovers the star-graph bound in
\cite[Theorem 3.4]{BegMa:gra}.

For the Dynkin diagrams, root the graph at the chosen vertex \(c\).
An arm length counts edges from this vertex to a leaf. The following
values of \(q\) give \(\lambda=q+q^{-1}=(2)_q\) and a positive
adjacency eigenvector, including the affine \(\widetilde E\) diagrams:
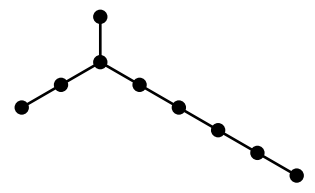
\begin{figure}[tbp]
\centering
\begin{minipage}[c]{.72\textwidth}
\centering
\begin{tabular}{ccc}
\toprule
Graph&Arm lengths&\(q\)\\
\midrule
\(A_n\), \(n\ge3\)&\((\lfloor(n-1)/2\rfloor,\lceil(n-1)/2\rceil)\)&\(e^{\imath{}\pi/(n+1)}\)\\
\(D_n\), \(n\ge4\)&\((1,1,n-3)\)&\(e^{\imath{}\pi/(2n-2)}\)\\
\(E_6\)&\((1,2,2)\)&\(e^{\imath{}\pi/12}\)\\
\(E_7\)&\((1,2,3)\)&\(e^{\imath{}\pi/18}\)\\
\(E_8\)&\((1,2,4)\)&\(e^{\imath{}\pi/30}\)\\
\(\widetilde E_6\)&\((2,2,2)\)&\(1\)\\
\(\widetilde E_7\)&\((1,3,3)\)&\(1\)\\
\(\widetilde E_8\)&\((1,2,5)\)&\(1\)\\
\bottomrule
\end{tabular}
\end{minipage}\hfill
\begin{minipage}[c]{.24\textwidth}

\centering
\begin{tikzpicture}[line width=.5pt]
\path[use as bounding box] (-1.05,-.72) rectangle (1.05,.72);
\coordinate (v0) at (0.00000,-0.15000);
\coordinate (v1) at (0.00000,0.15000);
\coordinate (v2) at (0.00000,0.45000);
\coordinate (v3) at (-0.25981,-0.30000);
\coordinate (v4) at (-0.51962,-0.45000);
\coordinate (v5) at (0.25981,-0.30000);
\coordinate (v6) at (0.51962,-0.45000);
\draw (v0)--(v1);
\draw (v1)--(v2);
\draw (v0)--(v3);
\draw (v3)--(v4);
\draw (v0)--(v5);
\draw (v5)--(v6);
\fill (v0) circle (1.35pt);
\fill (v1) circle (1.35pt);
\fill (v2) circle (1.35pt);
\fill (v3) circle (1.35pt);
\fill (v4) circle (1.35pt);
\fill (v5) circle (1.35pt);
\fill (v6) circle (1.35pt);
\end{tikzpicture}\par\vspace{2pt}
{\scriptsize \(\widetilde E_6\)}

\par\vspace{8pt}
\begin{tikzpicture}[line width=.5pt]
\path[use as bounding box] (-1.05,-.72) rectangle (1.05,.72);
\coordinate (v0) at (0.00000,0.07500);
\coordinate (v1) at (0.00000,0.37500);
\coordinate (v2) at (-0.25981,-0.07500);
\coordinate (v3) at (-0.51962,-0.22500);
\coordinate (v4) at (-0.77942,-0.37500);
\coordinate (v5) at (0.25981,-0.07500);
\coordinate (v6) at (0.51962,-0.22500);
\coordinate (v7) at (0.77942,-0.37500);
\draw (v0)--(v1);
\draw (v0)--(v2);
\draw (v2)--(v3);
\draw (v3)--(v4);
\draw (v0)--(v5);
\draw (v5)--(v6);
\draw (v6)--(v7);
\fill (v0) circle (1.35pt);
\fill (v1) circle (1.35pt);
\fill (v2) circle (1.35pt);
\fill (v3) circle (1.35pt);
\fill (v4) circle (1.35pt);
\fill (v5) circle (1.35pt);
\fill (v6) circle (1.35pt);
\fill (v7) circle (1.35pt);
\end{tikzpicture}\par\vspace{2pt}
{\scriptsize \(\widetilde E_7\)}

\par\vspace{8pt}
\begin{tikzpicture}[line width=.5pt]
\path[use as bounding box] (-1.05,-.72) rectangle (1.05,.72);
\coordinate (v0) at (-0.38971,0.22500);
\coordinate (v1) at (-0.38971,0.52500);
\coordinate (v2) at (-0.64952,0.07500);
\coordinate (v3) at (-0.90933,-0.07500);
\coordinate (v4) at (-0.12990,0.07500);
\coordinate (v5) at (0.12990,-0.07500);
\coordinate (v6) at (0.38971,-0.22500);
\coordinate (v7) at (0.64952,-0.37500);
\coordinate (v8) at (0.90933,-0.52500);
\draw (v0)--(v1);
\draw (v0)--(v2);
\draw (v2)--(v3);
\draw (v0)--(v4);
\draw (v4)--(v5);
\draw (v5)--(v6);
\draw (v6)--(v7);
\draw (v7)--(v8);
\fill (v0) circle (1.35pt);
\fill (v1) circle (1.35pt);
\fill (v2) circle (1.35pt);
\fill (v3) circle (1.35pt);
\fill (v4) circle (1.35pt);
\fill (v5) circle (1.35pt);
\fill (v6) circle (1.35pt);
\fill (v7) circle (1.35pt);
\fill (v8) circle (1.35pt);
\end{tikzpicture}\par\vspace{2pt}
{\scriptsize \(\widetilde E_8\)}
\end{minipage}
\caption{Arm lengths and phases, with the affine \(\widetilde E\) diagrams shown at right.}
\label{fig:dynkin}
\end{figure}
For the finite \(A,D,E\) types, the denominator in
\(q=e^{\imath{}\pi/h}\) is the Coxeter number \(h\) (or dual
Coxeter number, since these types are simply laced). The affine
\(\widetilde E\) diagrams instead have adjacency eigenvalue \(2\) and
\(q=1\) in \eqref{eq:dynkin-centre} below; the finite-type formula
\(\lambda=2\cos(\pi/h)\) does not apply to them.
On an arm of length \(\ell\), write \(x_0=c,x_1,\ldots,x_\ell\)
for its successive vertices. With \(v_c=1\), the eigenvector is
\[
 v_{x_j}=\frac{(\ell+1-j)_q}{(\ell+1)_q},\qquad 0\le j\le\ell.
\]
The \(q\)-integer recurrence gives the eigenvector equation along each arm;
at the centre it is
\begin{equation}\label{eq:dynkin-centre}
 \lambda=\sum_{\text{arms}}\frac{(\ell)_q}{(\ell+1)_q}=q+q^{-1},
\end{equation}
which holds for the listed lengths and values of \(q\).

\begin{corollary}\label{cor:dynkin}
 \(A_n\), \(n\ge3\), \(D_n\), \(n\ge4\), and
\(E_6,E_7,E_8,\widetilde E_6,\widetilde E_7,\widetilde E_8\) all admit
positive ratio metrics with star-preserving QLCs. The outward
coefficients are arbitrary nonzero reals, with reverse coefficients
determined by
\[
 g_{x_{j+1}x_j}=
 \frac{(\ell+1-j)_q}{(\ell-j)_q}
 g_{x_jx_{j+1}},\qquad 0\le j<\ell.
\]
For such a metric, $D_n, E_6,\cdots,E_8$ have exactly two QLCs, with parameter $z=q,q^{-1}$; exactly one for each affine $\widetilde E$ diagram. For $A_n$ there is a modulus 1 space of QLCs\cite{ArgMa:dis} that includes the two QLCs constructed here.  
\end{corollary}
\begin{proof}
The displayed eigenvector is positive, hence is the Perron
eigenvector and is unique up to scale. Theorem~\ref{thm:one-branch}
therefore gives both existence and the complete positive ratio metric 
moduli. Here \(0<\lambda\le2\), with equality precisely for the affine diagrams.
\end{proof}

\begin{remark}
In the affine cases, \(q=1\) and \((m)_1=m\). For example, on
\(\widetilde E_8\), take \(v_c=6\); the eigenvector values are \(3\) on the short arm, \(4,2\) on the next, and
\(5,4,3,2,1\) on the long arm, listed outward from the centre.
The phase \(z=1\) corresponds to \(s=-1\) in the star-graph
convention \(s=-z\); similarly \(D_4\) recovers the three-star
with \(s=-e^{\pm\imath{}\pi/6}\).
For \(n\ge10\), the tree with arms \((1,2,n-4)\) properly contains
\(\widetilde E_8\), so its Perron eigenvalue is strictly greater than \(2\).
Theorem~\ref{thm:one-branch} therefore excludes positive ratio metrics. This also follows from \eqref{eq:arm-bound} in
Corollary~\ref{cor:positive-trees} below.
This does not exclude metrics with negative ratios, which are governed by
the other adjacency eigenvectors in that theorem.
\end{remark}

\begin{remark}
How our solutions for $A_n$ sit inside the modulus 1 family in \cite{ArgMa:dis} proceeds as follows. 
Label the discrete interval \(A_n\) by \(1,\ldots,n\) in path
order. The required positive metric ratios are given by
\[
 q=e^{\imath{}\pi/(n+1)},\qquad v_i=(i)_q,\qquad
 \phi_i:=\frac{g_{i,i+1}}{g_{i+1,i}}=
                  \frac{(i+1)_q}{(i)_q},\qquad 1\le i<n,
\]
with one arbitrary nonzero real metric coefficient per edge and the other determined. 

For the distinguished pair in Theorem~\ref{thm:one-branch}, take
\(c=\lfloor(n+1)/2\rfloor\). For even \(n\), choosing the other
middle vertex exchanges the two solutions: the ratio across the
middle edge is one, and its phase is \(\bar z\). The corresponding values of the free phase \(s\) in \cite{ArgMa:dis} are determined as follows. Let
\(\tau_i=r_{i,i+1}-1\) be the coefficient of
\(\omega_{i,i+1,i}\) in its own image under \(\sigma\).
The free endpoint phase there is \(s=\tau_1\).
Our choice fixes it by the backward recurrence
\[
 \tau_c=\phi_c z-1,\qquad
 \tau_i=\frac{\phi_{i+1}}{1+\tau_{i+1}}-\phi_i
       \quad(i=c-1,\ldots,1),\qquad s=\tau_1.
\]
For \(A_3\) this gives \(s=-z\), recovering the two-star phases
used in Section~\ref{sec:star}. For longer paths it
selects two members of the full phase family, rather than proving
uniqueness. These differ from the real choices \(s=\pm1\) studied in \cite{ArgMa:dis}. 
\end{remark}

We finish by writing the distinguished connections on these diagrams in terms of arrows.
Number the arms by \(i=1,\ldots,d\), where \(d=2\) for \(A_n\)
and \(d=3\) for the finite \(D,E\) and affine \(\widetilde E\) diagrams. On arm \(i\), of length
\(\ell_i\), label the successive vertices by
\(x_{i,0}=c,x_{i,1},\ldots,x_{i,\ell_i}\), and set
\[
 u_{i,j}=\omega_{x_{i,j-1},x_{i,j}},\qquad
 \bar u_{i,j}=\omega_{x_{i,j},x_{i,j-1}},\qquad
 \phi_{i,j}=\frac{(\ell_i+2-j)_q}{(\ell_i+1-j)_q}.
\]
The metric and its free coefficients are
\[
 g=\sum_{i=1}^d\sum_{j=1}^{\ell_i}a_{i,j}
   (u_{i,j}\ten\bar u_{i,j}+\phi_{i,j}\bar u_{i,j}\ten u_{i,j}),
 \qquad a_{i,j}\in\R\setminus\{0\}.
\]
For \(z=q\) or \(q^{-1}\), define the return coefficients by
\begin{align*}
 \tau_{i,1}&=\frac{z}{\phi_{i,1}}-1,\qquad
 \tau'_{i,j}=\phi_{i,j}\overline{\tau_{i,j}},\\
 \tau_{i,j}&=-1+\frac{\phi_{i,j-1}}
 {\phi_{i,j}(1+\phi_{i,j-1}\tau_{i,j-1})}
 \quad(2\le j\le\ell_i).
\end{align*}
Thus all connection coefficients below are specified by \(q,z\)
and the metric values. Terms involving arrows beyond a leaf are
omitted; in particular, index \(2\) is absent on an arm of length one.
On the arrows meeting the centre,
\begin{align*}
 \nabla u_{i,1}
 &=\sum_{k=1}^d\bar u_{k,1}\ten u_{i,1}
   -\tau_{i,1}u_{i,1}\ten\bar u_{i,1}
   -(1+\tau_{i,1})\sum_{k\ne i}u_{k,1}\ten\bar u_{k,1}\\
 &\quad-\frac{a_{i,2}(\phi_{i,2}+\tau'_{i,2})}
                   {a_{i,1}\phi_{i,1}}u_{i,1}\ten u_{i,2},\\
 \nabla\bar u_{i,1}
 &=u_{i,1}\ten\bar u_{i,1}+\bar u_{i,2}\ten\bar u_{i,1}
   -\tau'_{i,1}\bar u_{i,1}\ten u_{i,1}\\
 &\quad-(1+\tau'_{i,1})u_{i,2}\ten\bar u_{i,2}
   -\bar z\sum_{k\ne i}\frac{a_{k,1}}{a_{i,1}}
                  \bar u_{i,1}\ten u_{k,1}.
\end{align*}
For \(2\le j\le\ell_i\), the remaining arrows have
\begin{align*}
 \nabla u_{i,j}
 &=u_{i,j-1}\ten u_{i,j}+\bar u_{i,j}\ten u_{i,j}
   -\tau_{i,j}u_{i,j}\ten\bar u_{i,j}\\
 &\quad-(1+\tau_{i,j})\bar u_{i,j-1}\ten u_{i,j-1}
   -\frac{a_{i,j+1}(\phi_{i,j+1}+\tau'_{i,j+1})}
                   {a_{i,j}\phi_{i,j}}u_{i,j}\ten u_{i,j+1},\\
 \nabla\bar u_{i,j}
 &=u_{i,j}\ten\bar u_{i,j}+\bar u_{i,j+1}\ten\bar u_{i,j}
   -\tau'_{i,j}\bar u_{i,j}\ten u_{i,j}\\
 &\quad-(1+\tau'_{i,j})u_{i,j+1}\ten\bar u_{i,j+1}
   -\frac{a_{i,j-1}(1+\overline{\tau'_{i,j-1}})}{a_{i,j}}
                   \bar u_{i,j}\ten\bar u_{i,j-1}.
\end{align*}
These expressions follow by expanding
\(\nabla\omega=\theta\ten\omega-\sigma(\omega\ten\theta)\)
using \eqref{eq:return-sigma}--\eqref{eq:straight-sigma}; the return
coefficients are \(\tau_{i,j}=Q_{x_{i,j-1},x_{i,j}}/g_{x_{i,j},x_{i,j-1}}\).

\subsection{All trees admitting a QRG with a positive ratio metric}\label{sec:posi}

Recall the four-star with centre \(0\), leaves \(1,\ldots,4\),
and arrows \(u_i=\omega_{0i}\), \(\bar u_i=\omega_{i0}\).
By Theorem~\ref{thm:one-branch}, real metrics admitting a star-preserving QLC are
precisely
\[
 g=\sum_{i=1}^4a_i(u_i\ten\bar u_i+2\eps\bar u_i\ten u_i),
 \qquad a_i\in\R\setminus\{0\},\qquad \eps\in\{1,-1\}.
\]
Each has a unique star-preserving QLC, with centre phase
\(z=\eps\) and \(Q_{0i}=Q_{i0}=-\eps a_i\).
The coefficients \(a_i\) have independent signs; \(\eps\) fixes
the reverse-to-forward ratio. The positive ratio case is the four-star
of \cite[Theorem 3.4]{BegMa:gra}, with phase \(s=-1\) there;
in general \(s=-\eps\). Explicitly,
\begin{align*}
 \nabla u_i&=\sum_{j=1}^4\bar u_j\ten u_i
       +u_i\ten\bar u_i-\frac12\sum_{j=1}^4u_j\ten\bar u_j,\\
 \nabla\bar u_i&=u_i\ten\bar u_i+\eps\bar u_i\ten u_i
       -\eps\sum_{j\ne i}\frac{a_j}{a_i}\bar u_i\ten u_j.
\end{align*}
Thus this connection has real coefficients for both values of
\(\eps\).

Split the centre of the four-star into two vertices, each retaining
two of the original leaves, and join the new vertices by an edge.
Subdividing that joining edge gives a path
\(x_0\to x_1\to\cdots\to x_m\), \(m\ge1\), with two leaves
attached to each endpoint. The resulting tree has \(m+5\) vertices.
Give each joining edge and each outward leaf edge an arbitrary
nonzero real coefficient, and prescribe the reverse coefficients by
\[
 g_{x_{j+1}x_j}=\eps g_{x_jx_{j+1}},\qquad
 g_{lx}=2\eps g_{xl}\quad(l\text{ a leaf attached to }x),
 \qquad\eps\in\{1,-1\}.
\]

\begin{proposition}\label{prop:split-four-star}
Every such metric has a unique star-preserving QLC for
\(\Om_{\min}\), specified by
\[
 Q_{x_jx_{j+1}}=Q_{x_{j+1}x_j}=0,\qquad
 Q_{xl}=Q_{lx}=-\eps g_{xl}.
\]
The metrics here with \(\epsilon=1\) are exactly the positive ratio metrics admitting a star-preserving QLC, and each has a unique such connection. We still have one free nonzero real value per edge. 
\end{proposition}
\begin{proof}
The ratios have \(H_x=2\eps\) at every vertex. The return
parameters are \((1,1/2,1/2)\) at each trivalent endpoint, in
the order joining edge, first leaf, second leaf; they are
\((1,1)\) at an internal joining vertex and \(1-\eps\) at
each leaf. They satisfy both equations in \eqref{eq:local-star}.
The real \(Q\)'s are symmetric, so the same tensor calculation
as in Theorem~\ref{thm:one-branch} gives metric compatibility;
torsion-freeness is built into the parametrisation.

For uniqueness, the phase at a trivalent endpoint is forced to be
\(\eps\), since \(z+\bar z=2\eps\) and \(|z|=1\).
This fixes its incident \(Q\)'s and propagation along the joining
path forces all the other \(Q\)'s. For a positive ratio metric,
Lemma~\ref{lem:quotient-independent} fixes its ratios by the
Perron eigenvector: value \(2\) at every joining vertex and
value \(1\) at each leaf, with eigenvalue \(2\).
\end{proof}

Splitting the four-star preserves its four leaf metric coefficients;
each joining edge has one new free coefficient. Subdivision adds a
bivalent vertex with return parameters \((1,1)\), leaving the
endpoint parameters unchanged. Taking \(m=2\), a single subdivision
gives exactly the seven-vertex binary tree \(T_3\). Its existence
therefore follows from this construction, with its full real-metric
classification and explicit geometry in Section~\ref{sec:seve}.

The affine diagram \(\widetilde D_n\), \(n\ge5\), has two
trivalent vertices, each bearing two leaves, joined by a path of
\(n-4\) edges. It has \(n+1\) vertices. At \(n=4\), identifying
the endpoints of the zero-length joining path gives the 4-star
\(\widetilde D_4\). The first two diagrams with distinct branching
vertices are \(\widetilde D_5\), the adjacent 3-stars, and
\(\widetilde D_6=T_3\).

\begin{corollary}\label{cor:fork}
 \(\widetilde D_n\), \(n\ge4\) all admit positive ratio metrics with a
unique star-preserving QLC for \(\Om_{\min}\). There is one free
nonzero real metric coefficient per edge; the reverse coefficient is
equal on the joining path and twice the outward coefficient on a
leaf edge.
\end{corollary}
\begin{proof}
For \(n\ge5\), apply Proposition~\ref{prop:split-four-star} with
\(m=n-4\) and \(\eps=1\). The case \(n=4\) is the four-star
construction at the start of this subsection.
\end{proof}

\begin{remark}\label{rem:affine}
The trees with arms \((2,2,2)\), \((1,3,3)\) and \((1,2,5)\)
are respectively \(\widetilde E_6,\widetilde E_7,\widetilde E_8\),
with seven, eight and nine vertices. Together with
\(\widetilde D_n\), these account for all affine diagrams in the
tree classification below.
\end{remark}

The spectral classification underlying the following result is the
tree case of Smith's theorem \cite{Smi}. We give a direct argument
to combine it with the QLC constructions above.

\begin{corollary}\label{cor:positive-trees}
A finite tree with at least two vertices admits a positive ratio
metric with a star-preserving QLC for \(\Om_{\min}\) if and only
if its adjacency Perron eigenvalue is at most \(2\).
Up to graph isomorphism, these trees are precisely the finite types
\[
 A_n\ (n\ge2),\qquad D_n\ (n\ge4),\qquad E_6,E_7,E_8,
\]
and the affine types
\[
 \widetilde D_n\ (n\ge4),\qquad
 \widetilde E_6,\widetilde E_7,\widetilde E_8.
\]
\end{corollary}
\begin{proof}
For a QRG with a positive ratio metric, Lemma~\ref{lem:quotient-independent}
identifies the common value \(H_x\) with the Perron eigenvalue.
If the tree has a branching vertex, Lemma~\ref{lem:local} bounds
it by \(2\); paths already have Perron eigenvalue less than \(2\).

To classify trees satisfying this bound, a vertex of degree at
least four gives a four-star subgraph, with Perron eigenvalue \(2\).
A connected proper extension of a connected graph has strictly larger
Perron eigenvalue, so the tree must be exactly that four-star.
If there are at least two trivalent vertices, their joining path,
together with the first edge of each of their four other branches,
is a double-fork subgraph. Its Perron eigenvalue is \(2\), so the
tree must equal this double fork.
The remaining non-path case has one trivalent vertex, with arm
lengths \(a\le b\le c\). Eliminating the arm vertices from
\(2I-A_G\), whose path blocks are positive definite, gives the
Schur complement
\[
 2-\frac{a}{a+1}-\frac{b}{b+1}-\frac{c}{c+1}.
\]
Thus the Perron bound is equivalent to
\begin{equation}\label{eq:arm-bound}
 \frac1{a+1}+\frac1{b+1}+\frac1{c+1}\ge1.
\end{equation}
Its positive integer solutions are \((1,1,c)\), \(c\ge1\),
giving \(D_{c+3}\); \((1,2,c)\), \(2\le c\le5\), giving
\(E_6,E_7,E_8,\widetilde E_8\); and the remaining triples
\((2,2,2)\) and \((1,3,3)\), giving
\(\widetilde E_6\) and \(\widetilde E_7\).

Conversely, Theorem~\ref{thm:one-branch} constructs the QLCs on
the paths with at least three vertices, the tripods and the four-star;
Corollary~\ref{cor:fork} supplies the double forks. On \(A_2\),
equal nonzero real metric coefficients and \(\sigma=\id\),
\(\alpha=0\), give a star-preserving QLC directly.
\end{proof}

For \(\widetilde D_n\), \(n\ge5\), we finish by writing the
resulting connections explicitly. The \(\widetilde D_4\) formulas
are those of the four-star above.
Put \(m=n-4\), and label the joining path
\(x_0,\ldots,x_m\). At each endpoint \(x_s\), \(s\in\{0,m\}\),
label its two leaves \(x_{s1},x_{s2}\). Write
\[
 u_j=\omega_{x_{j-1},x_j}\quad(1\le j\le m),\qquad
 u_{s\alpha}=\omega_{x_s,x_{s\alpha}}
 \quad(s\in\{0,m\},\ \alpha\in\{1,2\}),
\]
with bars denoting reverse arrows. The metrics considered in 
Proposition~\ref{prop:split-four-star} are
\[
 g=\sum_{j=1}^m a_j(u_j\ten\bar u_j+\eps\bar u_j\ten u_j)
  +\sum_{s\in\{0,m\}}\sum_{\alpha=1}^2a_{s\alpha}
       (u_{s\alpha}\ten\bar u_{s\alpha}
                    +2\eps\bar u_{s\alpha}\ten u_{s\alpha}),
\]
where \(a_j,a_{s\alpha}\) are arbitrary nonzero reals and
\(\eps\in\{1,-1\}\). The positive ratio members of this family are precisely those with \(\epsilon=1\), with no sign restriction on the remaining  coefficients.  Expanding \(\nabla=\theta\ten\id-\sigma(\id\ten\theta)\)
gives, on the joining arrows,
\begin{align*}
 \nabla u_j={}&\bar u_j\ten u_j+
 \begin{cases}
  \displaystyle\sum_{\alpha=1}^2
   (\bar u_{0\alpha}\ten u_1-u_{0\alpha}\ten\bar u_{0\alpha}),&j=1,\\
  u_{j-1}\ten u_j-\bar u_{j-1}\ten u_{j-1},&j>1,
 \end{cases}\\
 &-\begin{cases}
  \displaystyle\frac{a_{j+1}}{a_j}u_j\ten u_{j+1},&j<m,\\
  \displaystyle\sum_{\alpha=1}^2\frac{a_{m\alpha}}{a_m}
                         u_m\ten u_{m\alpha},&j=m,
 \end{cases}\\[3pt]
 \nabla\bar u_j={}&u_j\ten\bar u_j+
 \begin{cases}
  \bar u_{j+1}\ten\bar u_j-u_{j+1}\ten\bar u_{j+1},&j<m,\\
  \displaystyle\sum_{\alpha=1}^2
   (\bar u_{m\alpha}\ten\bar u_m-u_{m\alpha}\ten\bar u_{m\alpha}),&j=m,
 \end{cases}\\
 &-\begin{cases}
  \displaystyle\frac{a_{j-1}}{a_j}\bar u_j\ten\bar u_{j-1},&j>1,\\
  \displaystyle\eps\sum_{\alpha=1}^2\frac{a_{0\alpha}}{a_1}
                         \bar u_1\ten u_{0\alpha},&j=1.
 \end{cases}
\end{align*}
For the leaf arrows, let \(\bar\alpha\) denote the other element
of \(\{1,2\}\). Then
\begin{align*}
 \nabla u_{s\alpha}={}&\sum_{\beta=1}^2\bar u_{s\beta}\ten u_{s\alpha}
 +\frac12(u_{s\alpha}\ten\bar u_{s\alpha}
                   -u_{s\bar\alpha}\ten\bar u_{s\bar\alpha})\\
 &+\begin{cases}
  \bar u_1\ten u_{0\alpha}-\frac12u_1\ten\bar u_1,&s=0,\\
  u_m\ten u_{m\alpha}-\frac12\bar u_m\ten u_m,&s=m,
 \end{cases}\\[3pt]
 \nabla\bar u_{s\alpha}={}&u_{s\alpha}\ten\bar u_{s\alpha}
  +\eps\bar u_{s\alpha}\ten u_{s\alpha}
  -\eps\frac{a_{s\bar\alpha}}{a_{s\alpha}}
                         \bar u_{s\alpha}\ten u_{s\bar\alpha}\\
 &-\begin{cases}
  \displaystyle\eps\frac{a_1}{a_{0\alpha}}
                         \bar u_{0\alpha}\ten u_1,&s=0,\\
  \displaystyle\frac{a_m}{a_{m\alpha}}
                         \bar u_{m\alpha}\ten\bar u_m,&s=m.
 \end{cases}
\end{align*}
These formulae also apply when \(m=1\). For each displayed metric,
Proposition~\ref{prop:split-four-star} proves uniqueness of this
star-preserving QLC. All its coefficients are real, as for the
four-star and the \(\widetilde E_8\) connection for positive ratios.

\subsection{Trees admitting any QRG through to seven vertices}\label{sec:smal}

\begin{table}[tbp]
\centering
\begin{tabular}{cp{8cm}cc}
\toprule
Vertices&Tree&\shortstack{Positive ratio\\metric}&\shortstack{Any real\\metric}\\
\midrule
\(\ge 3\)&\(A_n\) (paths; $A_3$ is the 2-star)&Yes&Yes\\
\(\ge 4\)&\(D_n\) (\(D_4\) is the 3-star)&Yes&Yes\\
\(\ge 5\)& \(\widetilde D_{n-1}\) ($\widetilde D_4$ 4-star, $\widetilde D_5$ adjacent 3-stars, $\widetilde D_6=T_3$ )&Yes&Yes\\
\(6,7\)&\(E_6,E_7\) (similarly $E_8$)&Yes&Yes\\
\(6\)&arms \((1,1,1,2)\)&No&Yes\\
\(6\)&5-star&No&No\\
\(7\)&\(\widetilde E_6\) (similarly $\widetilde E_7, \widetilde E_8$) &Yes&Yes\\
\(7\)&arms \((1,1,1,1,2)\)&No&Yes\\
\(7\)&arms \((1,1,1,3)\)&No&Yes\\
\(7\)&arms \((1,1,2,2)\)&No&Yes\\
\(7\)&adjacent 3-star and 4-star&No&No\\
\(7\)&adjacent 3-stars with one arm subdivided &No&No\\
\(7\)&6-star&No&No\\
\bottomrule
\end{tabular}
\caption{Existence of QRGs on all trees with three to seven vertices.}
\label{tab:small-trees}
\end{table}

Table~\ref{tab:small-trees} lists all trees with three to seven
vertices up to graph isomorphism. Here an \(n\)-star has \(n\)
leaves, and arm lengths count edges from the unique branching
vertex.  We have already classified graphs admitting QRGs with positive ratio metrics 
in Corollary~\ref{cor:positive-trees} but for all QRGs (without the positive ratio
restriction), a case by case analysis is needed, which we report here up to seven vertices. 
As before, we require our metrics to be real and QLCs to be star-preserving. 
The last column means we allow negative metric ratios, meaning there could be edges where the metric value in one direction has a different sign from the value in the other direction.

\begin{figure}[tbp]
\centering
\noindent\makebox[.04\linewidth][l]{(a)}\makebox[.92\linewidth][c]{%
\begin{minipage}[c]{.215\linewidth}\centering
\begin{tikzpicture}[line width=.5pt]
\path[use as bounding box] (-1.15,-0.65) rectangle (1.15,0.65);
\coordinate (v0) at (-0.28750,0.00000);
\coordinate (v1) at (-0.28750,0.57500);
\coordinate (v2) at (-0.86250,0.00000);
\coordinate (v3) at (-0.28750,-0.57500);
\coordinate (v4) at (0.28750,-0.00000);
\coordinate (v5) at (0.86250,-0.00000);
\draw (v0)--(v1);
\draw (v0)--(v2);
\draw (v0)--(v3);
\draw (v0)--(v4);
\draw (v4)--(v5);
\fill (v0) circle (1.35pt);
\fill (v1) circle (1.35pt);
\fill (v2) circle (1.35pt);
\fill (v3) circle (1.35pt);
\fill (v4) circle (1.35pt);
\fill (v5) circle (1.35pt);
\end{tikzpicture}\par\vspace{2pt}
\parbox[t][.55cm][t]{\linewidth}{\centering\scriptsize arms \((1,1,1,2)\)}
\end{minipage}%
\hspace{.015\linewidth}%
\begin{minipage}[c]{.215\linewidth}\centering
\begin{tikzpicture}[line width=.5pt]
\path[use as bounding box] (-1.15,-0.65) rectangle (1.15,0.65);
\coordinate (v0) at (-0.30000,-0.06024);
\coordinate (v1) at (-0.30000,0.57063);
\coordinate (v2) at (-0.90000,0.13471);
\coordinate (v3) at (-0.67082,-0.57063);
\coordinate (v4) at (0.07082,-0.57063);
\coordinate (v5) at (0.30000,0.13471);
\coordinate (v6) at (0.90000,0.32966);
\draw (v0)--(v1);
\draw (v0)--(v2);
\draw (v0)--(v3);
\draw (v0)--(v4);
\draw (v0)--(v5);
\draw (v5)--(v6);
\fill (v0) circle (1.35pt);
\fill (v1) circle (1.35pt);
\fill (v2) circle (1.35pt);
\fill (v3) circle (1.35pt);
\fill (v4) circle (1.35pt);
\fill (v5) circle (1.35pt);
\fill (v6) circle (1.35pt);
\end{tikzpicture}\par\vspace{2pt}
\parbox[t][.55cm][t]{\linewidth}{\centering\scriptsize arms \((1,1,1,1,2)\)}
\end{minipage}%
\hspace{.015\linewidth}%
\begin{minipage}[c]{.215\linewidth}\centering
\begin{tikzpicture}[line width=.5pt]
\path[use as bounding box] (-1.15,-0.65) rectangle (1.15,0.65);
\coordinate (v0) at (-0.45000,0.00000);
\coordinate (v1) at (-0.45000,0.45000);
\coordinate (v2) at (-0.90000,0.00000);
\coordinate (v3) at (-0.45000,-0.45000);
\coordinate (v4) at (0.00000,-0.00000);
\coordinate (v5) at (0.45000,-0.00000);
\coordinate (v6) at (0.90000,-0.00000);
\draw (v0)--(v1);
\draw (v0)--(v2);
\draw (v0)--(v3);
\draw (v0)--(v4);
\draw (v4)--(v5);
\draw (v5)--(v6);
\fill (v0) circle (1.35pt);
\fill (v1) circle (1.35pt);
\fill (v2) circle (1.35pt);
\fill (v3) circle (1.35pt);
\fill (v4) circle (1.35pt);
\fill (v5) circle (1.35pt);
\fill (v6) circle (1.35pt);
\end{tikzpicture}\par\vspace{2pt}
\parbox[t][.55cm][t]{\linewidth}{\centering\scriptsize arms \((1,1,1,3)\)}
\end{minipage}%
\hspace{.015\linewidth}%
\begin{minipage}[c]{.215\linewidth}\centering
\begin{tikzpicture}[line width=.5pt]
\path[use as bounding box] (-1.15,-0.65) rectangle (1.15,0.65);
\coordinate (v0) at (-0.19167,0.19167);
\coordinate (v1) at (-0.19167,0.57500);
\coordinate (v2) at (-0.57500,0.19167);
\coordinate (v3) at (-0.19167,-0.19167);
\coordinate (v4) at (-0.19167,-0.57500);
\coordinate (v5) at (0.19167,0.19167);
\coordinate (v6) at (0.57500,0.19167);
\draw (v0)--(v1);
\draw (v0)--(v2);
\draw (v0)--(v3);
\draw (v3)--(v4);
\draw (v0)--(v5);
\draw (v5)--(v6);
\fill (v0) circle (1.35pt);
\fill (v1) circle (1.35pt);
\fill (v2) circle (1.35pt);
\fill (v3) circle (1.35pt);
\fill (v4) circle (1.35pt);
\fill (v5) circle (1.35pt);
\fill (v6) circle (1.35pt);
\end{tikzpicture}\par\vspace{2pt}
\parbox[t][.55cm][t]{\linewidth}{\centering\scriptsize arms \((1,1,2,2)\)}
\end{minipage}%
}\par\vspace{3pt}
\noindent\makebox[.04\linewidth][l]{(b)}\makebox[.92\linewidth][c]{%
\begin{minipage}[c]{.215\linewidth}\centering
\begin{tikzpicture}[line width=.5pt]
\path[use as bounding box] (-1.15,-0.65) rectangle (1.15,0.65);
\coordinate (v0) at (-0.00000,-0.06070);
\coordinate (v1) at (0.00000,0.57500);
\coordinate (v2) at (-0.60459,0.13574);
\coordinate (v3) at (-0.37366,-0.57500);
\coordinate (v4) at (0.37366,-0.57500);
\coordinate (v5) at (0.60459,0.13574);
\draw (v0)--(v1);
\draw (v0)--(v2);
\draw (v0)--(v3);
\draw (v0)--(v4);
\draw (v0)--(v5);
\fill (v0) circle (1.35pt);
\fill (v1) circle (1.35pt);
\fill (v2) circle (1.35pt);
\fill (v3) circle (1.35pt);
\fill (v4) circle (1.35pt);
\fill (v5) circle (1.35pt);
\end{tikzpicture}\par\vspace{2pt}
\parbox[t][.55cm][t]{\linewidth}{\centering\scriptsize 5-star}
\end{minipage}%
\hspace{.015\linewidth}%
\begin{minipage}[c]{.215\linewidth}\centering
\begin{tikzpicture}[line width=.5pt]
\path[use as bounding box] (-1.15,-0.65) rectangle (1.15,0.65);
\coordinate (v0) at (0.00000,0.00000);
\coordinate (v1) at (0.00000,0.57500);
\coordinate (v2) at (-0.49796,0.28750);
\coordinate (v3) at (-0.49796,-0.28750);
\coordinate (v4) at (-0.00000,-0.57500);
\coordinate (v5) at (0.49796,-0.28750);
\coordinate (v6) at (0.49796,0.28750);
\draw (v0)--(v1);
\draw (v0)--(v2);
\draw (v0)--(v3);
\draw (v0)--(v4);
\draw (v0)--(v5);
\draw (v0)--(v6);
\fill (v0) circle (1.35pt);
\fill (v1) circle (1.35pt);
\fill (v2) circle (1.35pt);
\fill (v3) circle (1.35pt);
\fill (v4) circle (1.35pt);
\fill (v5) circle (1.35pt);
\fill (v6) circle (1.35pt);
\end{tikzpicture}\par\vspace{2pt}
\parbox[t][.55cm][t]{\linewidth}{\centering\scriptsize 6-star}
\end{minipage}%
\hspace{.015\linewidth}%
\begin{minipage}[c]{.215\linewidth}\centering
\begin{tikzpicture}[line width=.5pt]
\path[use as bounding box] (-1.15,-0.65) rectangle (1.15,0.65);
\coordinate (v0) at (-0.25875,0.00000);
\coordinate (v1) at (0.54625,0.00000);
\coordinate (v2) at (-0.25875,0.57500);
\coordinate (v3) at (-0.83375,0.00000);
\coordinate (v4) at (-0.25875,-0.57500);
\coordinate (v5) at (0.83375,-0.49796);
\coordinate (v6) at (0.83375,0.49796);
\draw (v0)--(v1);
\draw (v0)--(v2);
\draw (v0)--(v3);
\draw (v0)--(v4);
\draw (v1)--(v5);
\draw (v1)--(v6);
\fill (v0) circle (1.35pt);
\fill (v1) circle (1.35pt);
\fill (v2) circle (1.35pt);
\fill (v3) circle (1.35pt);
\fill (v4) circle (1.35pt);
\fill (v5) circle (1.35pt);
\fill (v6) circle (1.35pt);
\end{tikzpicture}\par\vspace{2pt}
\parbox[t][.55cm][t]{\linewidth}{\centering\scriptsize adjacent 3-, 4-stars}
\end{minipage}%
\hspace{.015\linewidth}%
\begin{minipage}[c]{.215\linewidth}\centering
\begin{tikzpicture}[line width=.5pt]
\path[use as bounding box] (-1.15,-0.65) rectangle (1.15,0.65);
\coordinate (v0) at (-0.19919,-0.19167);
\coordinate (v1) at (0.42050,-0.19167);
\coordinate (v2) at (-0.42050,0.19167);
\coordinate (v3) at (-0.42050,-0.57500);
\coordinate (v4) at (0.64182,-0.57500);
\coordinate (v5) at (0.64182,0.19167);
\coordinate (v6) at (-0.64182,0.57500);
\draw (v0)--(v1);
\draw (v0)--(v2);
\draw (v0)--(v3);
\draw (v1)--(v4);
\draw (v1)--(v5);
\draw (v2)--(v6);
\fill (v0) circle (1.35pt);
\fill (v1) circle (1.35pt);
\fill (v2) circle (1.35pt);
\fill (v3) circle (1.35pt);
\fill (v4) circle (1.35pt);
\fill (v5) circle (1.35pt);
\fill (v6) circle (1.35pt);
\end{tikzpicture}\par\vspace{2pt}
\parbox[t][.55cm][t]{\linewidth}{\centering\scriptsize adjacent 3-stars\\one edge subdivided}
\end{minipage}%
}\par\vspace{3pt}
\noindent\makebox[.04\linewidth][l]{(c)}\makebox[.92\linewidth][c]{%
\begin{minipage}[c]{.30\linewidth}\centering
\begin{tikzpicture}[line width=.5pt]
\path[use as bounding box] (-1.15,-0.65) rectangle (1.15,0.65);
\coordinate (v0) at (0.00000,0.30312);
\coordinate (v1) at (-0.52500,0.00000);
\coordinate (v2) at (0.52500,0.00000);
\coordinate (v3) at (-1.05000,-0.30312);
\coordinate (v4) at (0.00000,-0.30312);
\coordinate (v5) at (1.05000,-0.30312);
\draw (v0)--(v1);
\draw (v0)--(v2);
\draw (v1)--(v3);
\draw (v1)--(v4);
\draw (v2)--(v4);
\draw (v2)--(v5);
\fill (v0) circle (1.35pt);
\fill (v1) circle (1.35pt);
\fill (v2) circle (1.35pt);
\fill (v3) circle (1.35pt);
\fill (v4) circle (1.35pt);
\fill (v5) circle (1.35pt);
\end{tikzpicture}\par\vspace{2pt}
\parbox[t][.55cm][t]{\linewidth}{\centering\scriptsize six-vertex quotient}
\end{minipage}%
\hspace{.07\linewidth}%
\begin{minipage}[c]{.30\linewidth}\centering
\begin{tikzpicture}[line width=.5pt]
\path[use as bounding box] (-1.15,-0.65) rectangle (1.15,0.65);
\coordinate (v0) at (0.00000,0.40981);
\coordinate (v1) at (0.00000,0.10981);
\coordinate (v2) at (-0.60000,-0.23660);
\coordinate (v3) at (0.60000,-0.23660);
\coordinate (v4) at (-0.90000,-0.40981);
\coordinate (v5) at (-0.30000,-0.40981);
\coordinate (v6) at (0.30000,-0.40981);
\coordinate (v7) at (0.90000,-0.40981);
\draw (v0)--(v1);
\draw (v1)--(v2);
\draw (v1)--(v3);
\draw (v2)--(v4);
\draw (v2)--(v5);
\draw (v3)--(v6);
\draw (v3)--(v7);
\fill (v0) circle (1.35pt);
\fill (v1) circle (1.35pt);
\fill (v2) circle (1.35pt);
\fill (v3) circle (1.35pt);
\fill (v4) circle (1.35pt);
\fill (v5) circle (1.35pt);
\fill (v6) circle (1.35pt);
\fill (v7) circle (1.35pt);
\end{tikzpicture}\par\vspace{2pt}
\parbox[t][.55cm][t]{\linewidth}{\centering\scriptsize \(T_3\) with root leg}
\end{minipage}%
}\par\vspace{3pt}
\caption{(a) Trees  up to 7 vertices admitting real QRGs but none with positive ratio metrics. (b) Trees up to 7 vertices admitting no QRGs. (c) The six-vertex quotient and eight-vertex augmentation of the binary tree \(T_3\) also admit no QRGs.}
\label{fig:tree-variations}
\end{figure}
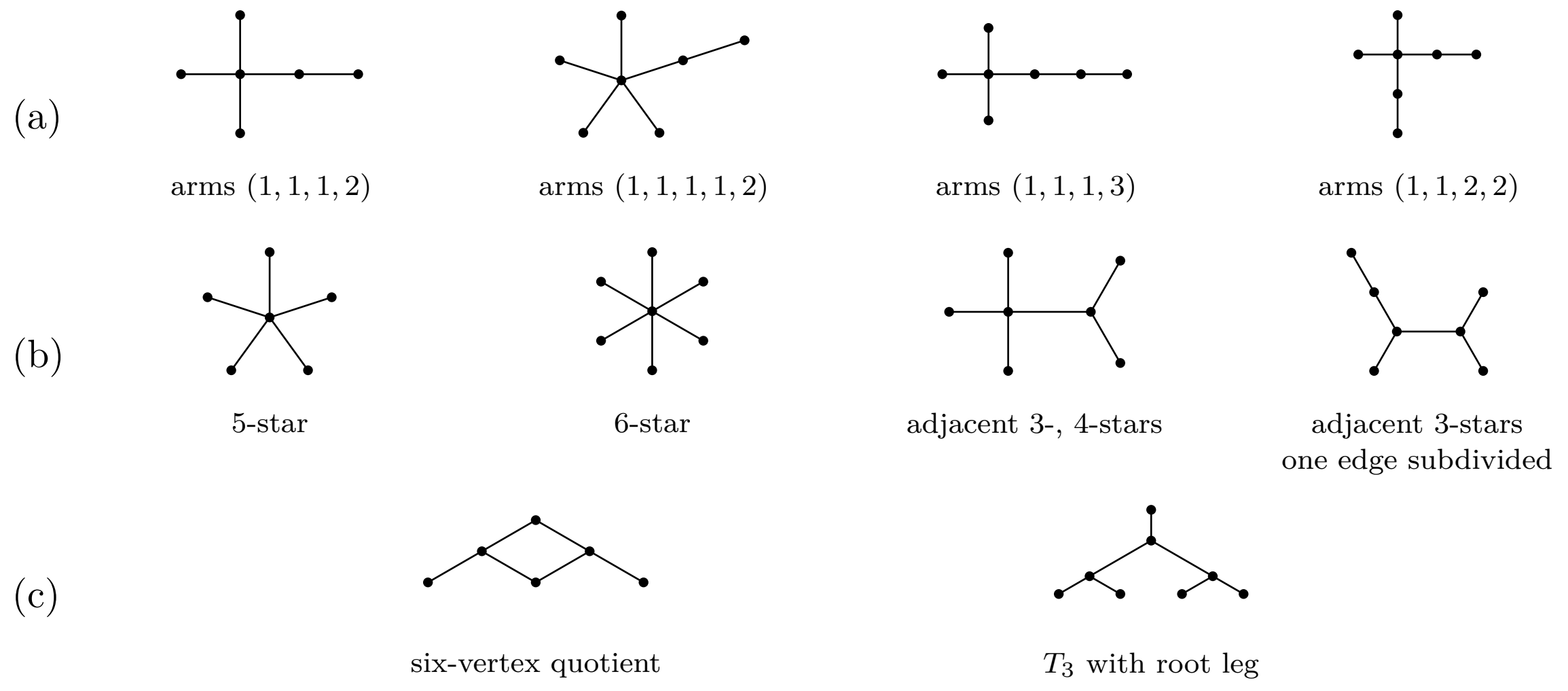

Existence for the four cases requiring negative ratios in the table (shown in Figure~\ref{fig:tree-variations}(a)) follows directly from
Theorem~\ref{thm:one-branch}. In their table order, choose
\[
 \lambda^2=\frac{5-\sqrt{13}}2,\quad
 3-\sqrt5,\quad 3-\sqrt2,\quad\frac{5-\sqrt{17}}2.
\]
The corresponding real adjacency eigenvectors have no zero entries
and change sign. Replacing the minus signs by plus signs gives the
squares of their Perron eigenvalues, greater than four, excluding
positive ratio metrics.

We have already covered non-existence of any QRGs for stars. The remaining non-existence cases in the table (shown in Figure~\ref{fig:tree-variations}(b)) follow from
Lemma~\ref{lem:local}. For adjacent branching vertices, either
\(h_{xy}=1\) or both phases are \(-1\), with \(\lambda=-2\).
For the degree-four/degree-three pair, writing \(h=h_{xy}\) from
the degree-four vertex gives
\[
 h=\lambda-\frac3\lambda,\qquad
 h^{-1}=\lambda-\frac2\lambda.
\]
Neither alternative is possible. For the two trivalent vertices
with one extended leaf, take \(x\) at the end with two unextended
leaves. Then
\[
 h=\lambda-\frac2\lambda,\qquad
 h^{-1}+\frac1\lambda+\frac{\lambda}{\lambda^2-1}=\lambda.
\]
Here \(\lambda=0,\pm1\) would force a zero eigenvector entry.
The choice \(h=1\) leaves only \(\lambda=2\), which fails the
second equation; \(\lambda=-2\) fails it as well.

Two variations of the binary tree \(T_3\) are also shown in 
Figure~\ref{fig:tree-variations}(c) as having no QRG. The augmented tree $\widehat{T_3}$ is relevant to the
Hopf algebra of rooted trees and already excluded by Theorem~\ref{thm:no-go}. Meanwhile, the quotient of $T_3$ on identifying the leaves \(12=21\) is relevant to 
causal sets but the proof of no QRG here requires different methods due to the square equations and is omitted. Beyond the binary trees, their augmented variants and the star graphs, many other trees with eight or more vertices appear to admit no QRG, even allowing negative metric ratios. 

\begin{remark}[Polygons]\label{rem:poly}
The remaining simply-laced affine diagrams are
\(\widetilde A_{n-1}=\Z_n\), \(n\ge3\). For \(n\ge5\), the
directed-metric classification is due to Bochniak, Sitarz and
Zalecki \cite[Theorem 5.4, Remark 5.5 and Proposition 5.7]{BSZ}.
For positive directed coefficients, their result gives
\[
 g_{i,i+1}=a_i,\qquad g_{i+1,i}=\phi a_i,
 \qquad a_i>0,\quad\phi>0,
\]
with one free value per edge and one free reverse factor $\phi$, generalising the previous $\phi=1$ case for edge-symmetric metrics in \cite{ArgMa:qua,Ma:par}. 
Writing \(e^+=\sum_i\omega_{i,i+1}\) and
\(e^-=\sum_i\omega_{i,i-1}\), the unique star-preserving QLC is
\[
 \nabla e^+=\left(1-\frac{a_{i+1}}{a_i}\right)e^+\ten e^+,
 \qquad
 \nabla e^-=\left(1-\frac{a_{i-2}}{a_{i-1}}\right)e^-\ten e^-.
\]
the same as in the edge-symmetric case. Coefficients are functions of the initial vertex \(i\), placed on
the left, and indices are modulo \(n\). The same displayed construction also works for \(\Z_3\), by direct
substitution. For \(n=3\) or \(n\ge5\), \(\Om_{\min}\) coincides with the
Grassmann calculus used in these works so the above still applies.

For \(n=4\), however, the straight-path relation in
\(\Om_{\min}\) is only
\(e^+\wedge e^++e^-\wedge e^-=0\), so the above connection is
not torsion-free. For the square \(\Z_4=\widetilde A_3\), consider
edge-symmetric metrics with equal positive values on opposite edges:
\[
 g_{i,i+1}=g_{i+1,i}=a_i,\qquad
 a_i=\begin{cases}a&i\text{ even},\\ b&i\text{ odd},\end{cases}
 \qquad a,b>0.
\]
For each such metric, the full family of star-preserving QLCs is
parametrised by \(c,d\in\R\) satisfying
\[
 (b-a)(c^2+d^2)=2abc.
\]
This is a circle when \(a\ne b\), and the line \(c=0\) when
\(a=b\). The corresponding (inner) connection is flat and given by
\begin{align*}
 \nabla e^+&=2e^+\ten e^+
             -\left(1-\frac{(-1)^ic+\imath d}{a_i}\right)(e^+-e^-)^{\tens_A 2},\\
 \nabla e^-&=2e^-\ten e^-
             -\left(1+\frac{(-1)^ic+\imath d}{a_{i-1}}\right)(e^+-e^-)^{\tens_A 2}.
\end{align*}
The coefficients are placed on the left and indices are modulo four.
The point \(c=d=0\) is a canonical choice independent of the values of this limited class of metrics. The proof is omitted.  The general directed-metric
problem for the square appears significantly more complicated and may be looked at elsewhere.
\end{remark}

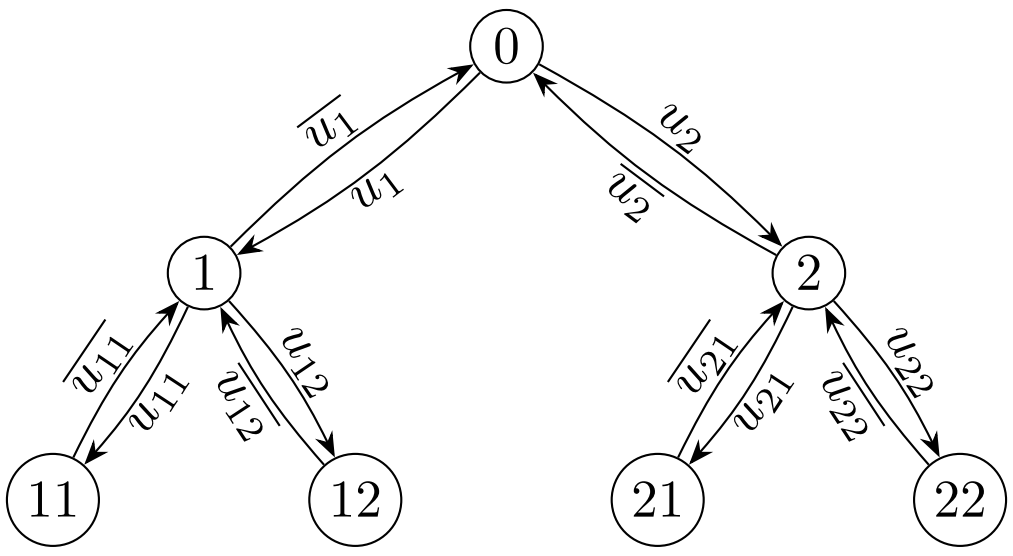
\begin{figure}[tp]
\centering
\begin{tikzpicture}[>=Stealth,every node/.style={font=\small},vertex/.style={circle,draw,inner sep=2pt}]
\node[vertex](r) at (0,3){0};
\node[vertex](b1) at (-2,1.5){1};\node[vertex](b2) at (2,1.5){2};
\node[vertex](l11) at (-3,0){11};\node[vertex](l12) at (-1,0){12};
\node[vertex](l21) at (1,0){21};\node[vertex](l22) at (3,0){22};
\foreach \a/\b/\lab in {r/b1/u_1,b1/l11/u_{11},b2/l21/u_{21}}{
\draw[->,bend left=8] (\a) to node[midway,sloped,below=0.7pt,inner sep=0pt] {$\lab$} (\b);
\draw[->,bend left=8] (\b) to node[midway,sloped,above=0.7pt,inner sep=0pt] {$\overline{\lab}$} (\a);
}
\foreach \a/\b/\lab in {r/b2/u_2,b1/l12/u_{12},b2/l22/u_{22}}{
\draw[->,bend left=8] (\a) to node[midway,sloped,above=0.7pt,inner sep=0pt] {$\lab$} (\b);
\draw[->,bend left=8] (\b) to node[midway,sloped,below=0.7pt,inner sep=0pt] {$\overline{\lab}$} (\a);
}
\end{tikzpicture}
\caption{The twelve directed arrows of the seven-vertex tree.}
\label{fig:tree}
\end{figure}

\section{The seven-vertex binary tree \(\widetilde D_6\)}\label{sec:seve}

In this section we look in more detail at the QRG of the 7-vertex binary tree $T_3=\widetilde D_6$ in the general classification. 

\subsection{Arrows and the minimal exterior algebra}
Label the root \(0\), its children \(1,2\), and the two leaves
of child \(i\) by \(i1,i2\). Thus the leaf labels are
\(11,12,21,22\).
Set
\[
 u_i=\omega_{0i},\quad\bar u_i=\omega_{i0},\qquad
 u_{i\alpha}=\omega_{i,\,i\alpha},\quad
 \bar u_{i\alpha}=\omega_{i\alpha,\,i}.
\]
Bars on these symbols mean the reverse arrow, so \(u_i^*=-\bar u_i\)
and \(u_{i\alpha}^*=-\bar u_{i\alpha}\). Always
\(i,\alpha\in\{1,2\}\). For \(i\in\{1,2\}\), write \(\bar i\) for the other element: \(\bar1=2\), \(\bar2=1\).

Every non-returning two-path is zero in \(\Om^2_{\min}\).
The only remaining relations are
\begin{equation}\label{eq:twoform-relations}
\begin{gathered}
 u_1\wedge\bar u_1+u_2\wedge\bar u_2=0,\qquad
 \bar u_{i\alpha}\wedge u_{i\alpha}=0,\\
 \bar u_i\wedge u_i+u_{i1}\wedge\bar u_{i1}
                         +u_{i2}\wedge\bar u_{i2}=0.
\end{gathered}
\end{equation}
Consequently
\[
 \Om^2_{\min}=\operatorname{span}_{\C}
 \{u_1\wedge\bar u_1,\ u_{i\alpha}\wedge\bar u_{i\alpha}:
                                  i,\alpha\in\{1,2\}\},
 \qquad \dim\Om^2_{\min}=5.
\]
Each basis form is supported at its displayed base vertex on both sides;
in particular \(\Om^2_{\min}\) is central as an \(A\)-bimodule.
The exterior differential is
\(\extd\omega_{xy}=\omega_{xy}\wedge\omega_{yx}
+\omega_{yx}\wedge\omega_{xy}\).
Moreover \(\Om^m_{\min}=0\) for \(m\ge3\): the non-return relations
leave only alternating three-paths, and multiplying the return sum at
\(x\) by \(\omega_{xy}\) kills each remaining path \(xyxy\).

\subsection{Metrics, braiding and connection}

Existence follows from Proposition~\ref{prop:split-four-star} with
\(m=2\). We now classify all nondegenerate real directed metrics,
including the signed cases, in the notation of Figure~\ref{fig:tree}.
\begin{theorem}\label{thm:seven}
The nondegenerate real directed metrics admitting a star-preserving
QLC on this tree are exactly
\begin{equation}\label{eq:seven-metric}
\begin{split}
 g={}&\sum_i a_i(u_i\ten\bar u_i+\eps\bar u_i\ten u_i)+\sum_{i,\alpha}a_{i\alpha}
 (u_{i\alpha}\ten\bar u_{i\alpha}
       +2\eps\bar u_{i\alpha}\ten u_{i\alpha}),
\end{split}
\end{equation}
where $\eps\in\{-1,1\}$ and the six \(a_i,a_{i\alpha}\) are arbitrary nonzero reals.  Each metric has exactly one such QLC.
Note that the reverse-to-forward metric ratios are \(\eps\) on the root edges
and \(2\eps\) on the leaf edges, and that positive directed lengths correspond to \(\eps=1\) and positive real parameters.
\end{theorem}
\begin{proof}
At branch vertex \(i\), write \(H_i=2\Re z_i\). Its leaf ratios are
\(H_i\) by Lemma~\ref{lem:local}. Put
\(\phi_i=g_{i0}/g_{0i}\); then
\(\phi_i=H_i-2/H_i\ne0\). At the root set
\(h_i=1/\phi_i\), \(r_i=1+\bar z_i-h_i\). The singleton star equation
is \(r_1\bar r_2=h_1h_2\). The diagonal return equations therefore give
\[
 |r_i|^2-2\Re r_i+h_1h_2=0
 \quad\Longrightarrow\quad H_i=h_1+h_2.
\]
Thus \(H_1=H_2=H\), \(\phi_1=\phi_2=\phi\), and
\(H=2/\phi\). Together with \(\phi=H-2/H\), this gives
\(H^2=4\), hence \(H=2\eps\), \(\phi=\eps\),
\(z_1=z_2=\eps\). The metric has form \eqref{eq:seven-metric}, and
\(Q\) is forced to be zero on root edges and \(-\eps a_{i\alpha}\)
on leaf edges in either direction. This fixes \(\sigma\) and \(\nabla\).
Conversely Proposition~\ref{prop:split-four-star} with \(m=2\)
constructs the unique star-preserving QLC for every metric in
\eqref{eq:seven-metric}.
\end{proof}

The braiding is explicit on every two-path. On return blocks use
\eqref{eq:return-sigma}, with
\begin{equation}\label{eq:seven-r}
 r_{0i}=1,\qquad r_{i0}=1,\qquad
 r_{i,\,i\alpha}=\tfrac12,\qquad
 r_{i\alpha,\,i}=1-\eps.
\end{equation}
On non-returning paths,
\begin{equation}\label{eq:seven-straight}
 \sigma(\omega_{xyz})=\eps\frac{g_{yz}}{g_{yx}}\omega_{xyz}.
\end{equation}
For reference, the return matrices at the root, a branch vertex, and a
leaf respectively are
\[
 \begin{pmatrix}0&1\\1&0\end{pmatrix},\qquad
 \begin{pmatrix}0&1/2&1/2\\1&-1/2&1/2\\1&1/2&-1/2\end{pmatrix},\qquad
 (-\eps),
\]
where the branch ordering is parent, first child, second child.

Here are the connection values for all twelve arrows. In these and the
subsequent formulas, \(\bar\alpha\) denotes the other element of
\(\{1,2\}\), just as \(\bar i\) does for \(i\).
\begin{align}
\nabla u_i={}&\sum_{m=1}^2\bar u_m\ten u_i-u_{\bar i}\ten\bar u_{\bar i}
 -\sum_{\alpha=1}^2\frac{a_{i\alpha}}{a_i}u_i\ten u_{i\alpha},\label{eq:nabla-u}\\
\nabla\bar u_i={}&u_i\ten\bar u_i+
 \sum_\alpha(\bar u_{i\alpha}\ten\bar u_i-u_{i\alpha}\ten\bar u_{i\alpha})
 -\eps\frac{a_{\bar i}}{a_i}\bar u_i\ten u_{\bar i},\label{eq:nabla-ubar}\\
\nabla u_{i\alpha}={}&u_i\ten u_{i\alpha}
 +\sum_\gamma\bar u_{i\gamma}\ten u_{i\alpha}
 +\frac12\left(u_{i\alpha}\ten\bar u_{i\alpha}
              -u_{i\bar\alpha}\ten\bar u_{i\bar\alpha}
              -\bar u_i\ten u_i\right),\label{eq:nabla-v}\\
\nabla\bar u_{i\alpha}={}&u_{i\alpha}\ten\bar u_{i\alpha}
 +\eps\bar u_{i\alpha}\ten u_{i\alpha}
 -\frac{a_i}{a_{i\alpha}}\bar u_{i\alpha}\ten\bar u_i
 -\eps\frac{a_{i\bar\alpha}}{a_{i\alpha}}\bar u_{i\alpha}\ten u_{i\bar\alpha}.
\label{eq:nabla-vbar}
\end{align}
All indices in sums range over \(\{1,2\}\). The four displayed families
specify respectively two, two, four and four arrows.

\subsection{The full curvature}
We use \(R=(\extd\otimes\id-\id\wedge\nabla)\nabla\).
\begin{proposition}\label{prop:curvature}
The curvature on all arrows is
\begin{align}
 R(u_i)={}&\left(\eps\frac{a_i}{a_{\bar i}}-\frac{a_{i1}+a_{i2}}{2a_i}\right)(u_i\wedge\bar u_i)\ten u_i,\label{eq:R-u}\\
 R(\bar u_i)={}&\sum_{\alpha=1}^2(u_{i\alpha}\wedge\bar u_{i\alpha})\ten
 \left[\left(\eps\frac{a_{\bar i}}{a_i}-\frac{a_i}{a_{i\alpha}}\right)\bar u_i
       +\eps u_{i\alpha}
       -\eps\frac{a_{i\bar\alpha}}{a_{i\alpha}}u_{i\bar\alpha}\right],\label{eq:R-ubar}\\
 R(u_{i\alpha})={}&\frac12(u_{i\alpha}\wedge\bar u_{i\alpha})\ten
 \left[\frac{a_i}{a_{i\alpha}}\bar u_i
       +\left(\frac{a_{i\alpha}}{a_i}-\eps\right)u_{i\alpha}
       +\left(\frac{a_{i\bar\alpha}}{a_i}
              +\eps\frac{a_{i\bar\alpha}}{a_{i\alpha}}\right)u_{i\bar\alpha}\right]
 \nonumber\\
 &+\frac12(u_{i\bar\alpha}\wedge\bar u_{i\bar\alpha})\ten
 \left[-\frac{a_i}{a_{i\bar\alpha}}\bar u_i
       +\left(\frac{a_{i\alpha}}{a_i}
              -\eps\frac{a_{i\alpha}}{a_{i\bar\alpha}}\right)u_{i\alpha}
       +\left(\frac{a_{i\bar\alpha}}{a_i}+\eps\right)u_{i\bar\alpha}\right],\label{eq:R-v}\\
 R(\bar u_{i\alpha})={}&0.\label{eq:R-vbar}
\end{align}
In particular, no member of the metric family is flat.
\end{proposition}
\begin{proof}
Insert \eqref{eq:nabla-u}--\eqref{eq:nabla-vbar}, using
\(\extd\omega_{xy}=(\omega_{xy}\wedge\omega_{yx})+(\omega_{yx}\wedge\omega_{xy})\) and
\eqref{eq:twoform-relations}. In a term
\((\extd\omega_{xy})\ten\omega_{yz}\), only \((\omega_{yx}\wedge\omega_{xy})\) survives the
balanced tensor product. In \(\omega_{xy}\wedge\nabla\omega_{yz}\),
only terms whose first two legs return to \(x\) survive. These rules
give the displayed expressions. The coefficient of
\((u_{i1}\wedge\bar u_{i1})\ten u_{i1}\) in \(R(\bar u_i)\) is the nonzero constant
\(\eps\); \((u_{i1}\wedge\bar u_{i1}),(u_{i2}\wedge\bar u_{i2})\) are independent.
\end{proof}

\subsection{Ricci tensor and scalar}
We fix \(\iota_0\) as in \eqref{eq:uniform-lift} for the following formulas and use the contraction
convention of \cite{Ma:par}:
\begin{equation}\label{eq:ricci-definition}
 \Ric=((\ ,\ )\otimes\id\otimes\id)
       (\id\otimes\iota_0\otimes\id)(\id\otimes R)g,
 \qquad\Scal=(\ ,\ )\Ric.
\end{equation}
Then the complete Ricci tensor can be written compactly as
\begin{equation}\label{eq:ricci}
\begin{split}
\Ric=\sum_i\bigg[&\rho_i u_i\ten\bar u_i+\tau_i\bar u_i\ten u_i\\
 &+\sum_\alpha\bigl(\gamma_{i\alpha}u_i\ten u_{i\alpha}
   +\zeta_{i\alpha}\bar u_{i\alpha}\ten\bar u_i
   +\mu_{i\alpha}\bar u_{i\alpha}\ten u_{i\alpha}
   +\nu_{i\alpha}\bar u_{i\alpha}\ten u_{i\bar\alpha}\bigr)\bigg],
\end{split}
\end{equation}
where
\begin{align*}
 \rho_i&=\frac13\left[-\frac{2a_{\bar i}}{a_i}
       +\eps a_i\left(\frac1{a_{i1}}+\frac1{a_{i2}}\right)\right],
 &\tau_i&=\frac{a_i}{2a_{\bar i}}-\eps\frac{a_{i1}+a_{i2}}{4a_i},\\
 \gamma_{i\alpha}&=\frac13\left(\frac{a_{i\alpha}}{a_{i\bar\alpha}}-1\right),
 &\zeta_{i\alpha}&=\frac{\eps a_i}{3}
              \left(\frac2{a_{i\alpha}}+\frac1{a_{i\bar\alpha}}\right),\\
 \mu_{i\alpha}&=\frac13\left(\eps\frac{a_{i\alpha}}{a_i}
                  -2+\frac{a_{i\alpha}}{a_{i\bar\alpha}}\right),
 &\nu_{i\alpha}&=\frac13\left(\eps\frac{a_{i\bar\alpha}}{a_i}
                  +2\frac{a_{i\bar\alpha}}{a_{i\alpha}}-1\right).
\end{align*}
In particular the scalar is the following function on the seven vertices:
\begin{align}
 \Scal(0)&=\frac13\sum_{i=1}^2
       \left(\frac1{a_{i1}}+\frac1{a_{i2}}
                        -\frac{2\eps a_{\bar i}}{a_i^2}\right),\label{eq:scalar-root}\\
 \Scal(i)&=\frac1{2a_{\bar i}}-\eps\frac{a_{i1}+a_{i2}}{4a_i^2},\label{eq:scalar-branch}\\
 \Scal(i\alpha)&=\frac{\eps}{3a_i}
                    -\frac2{3a_{i\alpha}}+\frac1{3a_{i\bar\alpha}}.
\label{eq:scalar-leaf}
\end{align}
These expressions follow by applying \eqref{eq:uniform-lift} to
Proposition~\ref{prop:curvature}. For example the lift of \((u_{i\alpha}\wedge\bar u_{i\alpha})\)
has coefficient \(-1/3\) on \(\bar u_i\ten u_i\) and coefficient
\(\delta_{\alpha\gamma}-1/3\) on \(u_{i\gamma}\ten\bar u_{i\gamma}\), for \(\gamma=1,2\).
The inverse metric then selects the returning first two tensor legs.
Only the return terms of \eqref{eq:ricci} contribute to the scalar.

\begin{corollary}[Scalar-flat but curved metrics]\label{cor:scalar-flat}
For \(\iota_0\), the scalar vanishes identically exactly when
\[
 a_1=a_2=a\ne0,\qquad a_{i\alpha}=\eps a\quad\text{for all }i,\alpha.
\]
For positive directed lengths this is the one-parameter family with all
six outward weights equal to \(a>0\). The full curvature never vanishes.
\end{corollary}
\begin{proof}
The two leaf equations at each branch force
\(1/a_{i1}=1/a_{i2}=\eps/a_i\). The branch equations then force
\(a_1=a_2\); the root equation follows. Non-flatness is already proved.
\end{proof}

\begin{remark}[Dependence on the lift]

For this binary tree, another choice of lift,  supplied by the connection itself, is
\(\iota_\sigma[\tau]=(\tau-\sigma\tau)/2\) on the nonzero return
quotients, extended by zero on the zero blocks. Here \(\sigma P_x=P_x\)
and \(\sigma^2=I\) on the relevant blocks, so this descends and is a
section. Its weights are \((1/2,1/2)\) at the root and
\((1/2,1/4,1/4)\) at a branch, in parent/children order.
These weights are real, so for either sign of \(\eps\) it satisfies
\(\iota_\sigma(\xi^*)=-\dagger\iota_\sigma(\xi)\).

For \(\iota_\sigma\), the branch scalar in \eqref{eq:scalar-branch}
is unchanged, the root scalar is \(3/4\) of \eqref{eq:scalar-root}, and
\[
 \Scal_\sigma(i\alpha)=
 \frac{\eps}{2a_i}-\frac3{4a_{i\alpha}}+\frac1{4a_{i\bar\alpha}}.
\]
The scalar-flat solutions are the same, but the tensors and scalar functions
are generally different. Thus the lift is part of the curvature data.
\end{remark}

\begin{proposition}[No central lifted Ricci tensor]\label{prop:noncentral} 
For this connection, no bimodule section \(\iota\) of the wedge gives a
central Ricci tensor under \eqref{eq:ricci-definition}.
\end{proposition}
\begin{proof}
Use \eqref{eq:general-lift}, and at branch \(i\) abbreviate its child
weights by \(w_1,w_2\). Let \(X_\alpha\) be the coefficient of
\(u_i\ten u_{i\alpha}\) and \(Y_\alpha\) that of
\(\bar u_{i\alpha}\ten\bar u_i\). The curvature formulas give
\[
 X_\alpha=\frac{a_{i\alpha}}{a_{i\bar\alpha}}w_{\bar\alpha}-w_\alpha,
 \qquad
 Y_\alpha=\frac{\eps a_i}{a_{i\alpha}}(1+X_\alpha).
\]
Centrality would require both coefficients to vanish, since these
paths have distinct endpoints. This is impossible. The assertion
allows arbitrary complex lift weights, not only real ones.
\end{proof}
In particular, no bimodule section gives an Einstein relation
\(\Ric=\Lambda g\), even with \(\Lambda\in A\), since its right-hand
side would be central. At the positive scalar-flat metric, the explicit
tensor is
\[
 \Ric=\sum_{i,\alpha}\left(
 \bar u_{i\alpha}\ten\bar u_i
 +\frac23\bar u_{i\alpha}\ten u_{i\bar\alpha}\right)\ne0.
\]
Scalar-flatness therefore does not mean Ricci-flatness.

\subsection{Einstein--Hilbert action}
The Einstein--Hilbert action on a finite graph is a weighted sum of
the scalar curvature values over its vertices. The vertex measure
\(\mu\) is additional data, chosen by hand to obtain a physically
reasonable action. We then define
\begin{equation}\label{eq:EH-definition}
 S[a]=\sum_{x\in X}\mu(x)\Scal(x)
 \end{equation}
 where we sum over the vertex set $X$.

For our tree, we take \(\eps=1\), retain the lift
\(\iota_0\), and first consider the case of $\mu$ independent of the metric. A natural choice here is
\[
 \mu(0)=3,\qquad \mu(i)=4,\qquad \mu(i\alpha)=6.
\]
as a positive measure that gives a squared-difference expression for the action as follows.

\begin{proposition}\label{prop:EH-action}
The Einstein--Hilbert action for constants $\mu$ as above is
\begin{equation}\label{eq:EH-squares}
 S[a]=-
 \frac{2(a_1+a_2)(a_1-a_2)^2}{a_1^2a_2^2}
 -\sum_{i,\alpha=1}^2
       \frac{(a_i-a_{i\alpha})^2}{a_i^2a_{i\alpha}}.
\end{equation}
For positive metric parameters, \(S[a]\le0\), with equality exactly
when all six parameters are equal.
\end{proposition}
\begin{proof}
Substitution of \eqref{eq:scalar-root}--\eqref{eq:scalar-leaf} and the chosen weights into \eqref{eq:EH-definition} gives
\[
\begin{split}
 S[a]={}&-2\left(\frac{a_2}{a_1^2}+\frac{a_1}{a_2^2}
                         -\frac1{a_1}-\frac1{a_2}\right)-\sum_{i,\alpha=1}^2
       \left(\frac{a_{i\alpha}}{a_i^2}
                         +\frac1{a_{i\alpha}}-\frac2{a_i}\right).
\end{split}
\]
Factoring yields \eqref{eq:EH-squares}. Positivity of the parameters
makes each squared difference contribute a nonpositive term, proving
the assertion.
\end{proof}

The action here has dimensions of inverse length squared (it goes as $1/a$) which is not necessarily wrong but not the only option. In continuum geometry we use $\sqrt{|g|}$ where $g$ is the metric, while experience for the $A_3$ graph \cite[Appendix A]{ArgMa:dis} or the 4-star graph\cite{Ma:alg4} suggests a metric dependence of order $a$. This then results in a dimensionless action. In our case, we consider all positive linear combinations of the square-lengths and search for similar squares of differences and positivity features as above. This leads naturally to
\begin{equation}\label{eq:linear-metric-measure}
 \mu(0)=\frac32(a_1+a_2),\qquad
 \mu(i)=4a_i,\qquad
 \mu(i\alpha)=\frac32(3a_i+a_{\bar i}).
\end{equation}
These weights are positive for positive metric parameters and respect
the rooted tree symmetries. At an equal metric \(a_i=a_{i\alpha}=a\), they become
\(a\) times the original weights \((3,4,6)\).

\begin{proposition}[A dimensionless Einstein--Hilbert action]
\label{prop:linear-metric-action}
The measure \eqref{eq:linear-metric-measure} gives
\begin{equation}\label{eq:linear-metric-action}
 S[a] =-\frac{(a_1-a_2)^2(a_1^2+a_2^2)}{a_1^2a_2^2} -\sum_{i,\alpha}
             \frac{(a_i-a_{i\alpha})^2}{a_i a_{i\alpha}}.
\end{equation}
It is invariant under \(a_i\mapsto ca_i,\ a_{i\alpha}\mapsto ca_{i\alpha}\). For positive metric
parameters it is nonpositive, with equality exactly when all six
parameters are equal.
\end{proposition}
\begin{proof}
Since \(\mu(i\alpha)=\mu(0)+3a_i\),
substitution of \eqref{eq:scalar-root}--\eqref{eq:scalar-leaf}
and collection of the reciprocal leaf terms give
\[
 S[a]=-\left(\frac{a_1}{a_2}-1\right)^2
  -\left(\frac{a_2}{a_1}-1\right)^2
  -\sum_{i,\alpha}
       \left(\frac{a_i}{a_{i\alpha}}
              +\frac{a_{i\alpha}}{a_i}-2\right).
\]
Factoring the summands gives \eqref{eq:linear-metric-action}.
Its homogeneity, sign and zero set follow directly.
\end{proof}

\begin{remark}\label{rem:EH-minus}
For \(\eps=-1\), there is no analogous natural choice of measure
independent of the metric. A natural linear choice is
\[
\begin{gathered}
 \mu_-(0)=\frac32(a_1+a_2),\qquad \mu_-(i)=4a_i,\quad  \mu_-(i\alpha)=\frac92a_i+\frac32a_{\bar i}
                      +4a_{i\alpha}+2a_{i\bar\alpha}.
\end{gathered}
\]
Discarding an additive constant, its Einstein--Hilbert action is
\[
\begin{split}
 S_-[a]&:=\sum_{x\in X}\mu_-(x)\Scal(x)+16=\frac{(a_1-a_2)^2(a_1^2+4a_1a_2+a_2^2)}{a_1^2a_2^2}
   -\sum_{i,\alpha}\frac{(a_i-a_{i\alpha})^2}{a_i a_{i\alpha}}.
\end{split}
\]
This follows by substitution of
\eqref{eq:scalar-root}--\eqref{eq:scalar-leaf} with \(\eps=-1\).
The action is dimensionless and uses the same differences as
Proposition~\ref{prop:linear-metric-action}, but is indefinite for
positive parameters. It vanishes when all six parameters are equal.
\end{remark}

\begin{remark}[Actions for the connection lift] For the lift \(\iota_\sigma\) we have: 

(i) For \(\eps=1\), the constant measure $\mu=4$ gives
\[
 4\sum_{x\in X}\Scal_\sigma(x)=S[a],
\]
the action in Proposition~\ref{prop:EH-action}.  Moreover, the linear measure
\[
 \mu(0)=2(a_1+a_2),\qquad \mu(i)=4a_i,\qquad
 \mu(i\alpha)=3a_i+a_{\bar i}
\]
instead gives \[ \sum_x\mu(x)\Scal_\sigma(x)=S[a],\]
 the
dimensionless action in Proposition~\ref{prop:linear-metric-action}. 

(ii) For \(\eps=-1\)  the linear measure
\[
\begin{gathered}
 \mu(0)=2(a_1+a_2),\qquad \mu(i)=4a_i,\quad \mu(i\alpha)=3a_i+a_{\bar i}
                              +3a_{i\alpha}+a_{i\bar\alpha}.
\end{gathered}
\]
 gives, on addition of  an additive constant,
\[
 \sum_{x\in X}\mu(x)\Scal_\sigma(x)+16=S_-[a],
\]
the action in Remark~\ref{rem:EH-minus}.
\end{remark}

\section{Quantum gravity on the 7-vertex binary tree}\label{sec:quan}

In this section we compute quantum gravity expectation values on the 7-vertex binary tree using the dimensionless Einstein--Hilbert action
$S[a]$ of Proposition~\ref{prop:linear-metric-action}. This is for $\eps=1$ (so opposite arrows have metric coefficients of the same sign)
and the lift $\iota_0$. Following the finite-dimensional approach of \cite{Ma:yet}, we
integrate over the six metric parameters; the QLC is unique for each
metric. The integration measure on the space of metrics is additional
data, distinct from the vertex measure. We consider direct and
Liouville measures and in both cases we regulate the theory with metric values cut off to $\delta<a_i,a_{i\alpha}<L$ as needed.

Broadly similar metric and curvature fluctuations,
and interference effects in the case of the Lorentzian action, were obtained for the
action in Proposition~\ref{prop:EH-action} given by metric-independent vertex weights. The main difference is that 
this action has metric degree $-1$, so after scaling the metric values by $L$, its direct-measure distribution depends on
$GL$ rather than $G$, with the cutoff ratio $\delta/L$ held fixed. For Liouville measure, rescaling the metric values also
rescales the coupling. This action is not invariant under metric inversion. The two actions thus have different exact symmetries and cutoff dependence.

Numerical integrations are checked by repeating the calculation with
finer integration grids at representative couplings across each
plotted range. Unless stated otherwise, the resulting changes are
less than $10^{-6}$ in the dimensionless quantities plotted; where
integration tails are truncated, we also check stability on extending
the integration domain. For the regularised Lorentzian theory we
additionally halve the small imaginary part of the coupling and check
agreement to at least the same accuracy. This accuracy statement
refers to observed numerical stability, not a rigorous error bound.

\subsection{Quantum gravity with direct measure}

We first use the direct measure
\[
 \mathrm{d}^6a=\prod_i\mathrm{d}a_i
                  \prod_{i,\alpha}\mathrm{d}a_{i\alpha}.
\]

For \(G,L>0\) and \(0\le\delta<L\), define
\[
 Z_{\delta,L}(G)=\int_{(\delta,L)^6}
         e^{S[a]/G}\,\mathrm{d}^6a,\qquad
 \< O\>=\frac1{Z_{\delta,L}(G)}
       \int_{(\delta,L)^6}O[a]e^{S[a]/G}\,\mathrm{d}^6a.
\]
The coupling \(G\) is dimensionless. Since \(S[La]=S[a]\),
\[
 Z_{\delta,L}(G)=L^6 Z_{\delta/L,1}(G).
\]
Thus, at \(\delta=0\), the distribution of \(a_i/L,\ a_{i\alpha}/L\) depends only
on \(G\), rather than \(GL\). In particular,
\(\< a_i\>/L,\ \< a_{i\alpha}\>/L\), \(L\<\Scal(x)\>\) and
\(L\Delta\Scal(x)\) are functions of \(G\) alone.

\begin{proposition}
For \(G,L>0\), \(0<Z_{0,L}(G)\le L^6\). The first and second
scalar-curvature moments at each vertex are absolutely convergent
without a lower cutoff. The partition function on \((0,\infty)^6\)
diverges, with the exact scaling
\(Z_{0,L}(G)=L^6 Z_{0,1}(G)\).
\end{proposition}
\begin{proof}
The partition bound follows from \(S\le0\). Put \(E=-S\) and
write \(a_i=rb_i\), \(a_{i\alpha}=rb_{i\alpha}\), where \(r\)
is the maximum of the six metric parameters. The maximum of the six
ratios \(b_i,b_{i\alpha}\) is therefore one. On each of the six faces selected by a maximal
parameter, the measure is \(r^5\,\mathrm{d}r\,\mathrm{d}^5b\).
Every adjacent ratio among the six parameters is bounded by a
polynomial in \(1+E\), using the squared-difference expression for
\(E\). Following the edges joining these parameters bounds each inverse ratio
\(1/b_i\) and \(1/b_{i\alpha}\) by such a polynomial. Consequently
\( |\Scal(x)|^m\le r^{-m}P_m(E)\) for a polynomial \(P_m\).
The product \(P_m(E)e^{-E/G}\) is bounded, and the radial integral
\(\int_0^L r^{5-m}\,\mathrm{d}r\) is finite for \(m<6\).
This proves both moment assertions. The exact scaling and positivity
prove the divergence as \(L\to\infty\).
\end{proof}

This radial separation gives another exact statement: \(r/L\) has
density \(6t^5\), \(0<t<1\), independently of \(G\) and of the
ratios \(b_i,b_{i\alpha}\). As \(G\downarrow0\), the latter concentrate at
\(b_i=b_{i\alpha}=1\), since the action vanishes only at equal parameters.
For either edge parameter \(a=a_i\) or \(a=a_{i\alpha}\),
\[
 \frac{\< a\>}{L}\longrightarrow\frac67,
 \qquad
 \frac{\sqrt{\< a^2\>-\< a\>^2}}
      {\< a\>}\longrightarrow\frac1{4\sqrt3}.
\]
At \(G\to\infty\), the normalised distribution instead tends to
uniform measure on the metric cube, so \(\< a\>/L\to1/2\).
Unlike the metric observables, curvature is singular on the boundary
of that cube; its moments need not have finite strong-coupling limits.

For local metric observables, take \(a(i)=a_i\) at a branch and
\(a(i\alpha)=a_{i\alpha}\) at a leaf. At the root use
\(a(0)=(a_1+a_2)/2\). Root and branch means coincide by symmetry,
but their fluctuations need not: the root variance includes the
correlation between \(a_1\) and \(a_2\). We plot the normalised means
and relative standard deviations
\[
 \frac{\< a(x)\>}{L},\qquad
 \frac{\Delta a(x)}{\< a(x)\>},\qquad
 \Delta a(x)=\sqrt{\< a(x)^2\>-\< a(x)\>^2}.
\]
They depend only on \(G\). At weak coupling all three relative
standard deviations tend to \(1/(4\sqrt3)\); at strong coupling
independent uniform edge parameters give
\[
 \frac{\Delta a(0)}{\< a(0)\>}\longrightarrow\frac1{\sqrt6},
 \qquad
 \frac{\Delta a_i}{\< a_i\>},\quad
 \frac{\Delta a_{i\alpha}}{\< a_{i\alpha}\>}
                 \longrightarrow\frac1{\sqrt3}.
\]

\begin{figure}[tbp]
\centering
\begin{tikzpicture}
\begin{groupplot}[
 group style={group size=2 by 2,horizontal sep=1.8cm,vertical sep=1.7cm},
 width=.44\linewidth,height=5.7cm,
]
\nextgroupplot[xmode=log,xmin=.001,xmax=100,
 xlabel={$G$},grid=major,grid style={gray!20},
 tick label style={font=\small},label style={font=\small},
 legend style={font=\scriptsize,draw=none},
ylabel={$\< a\>/L$},ymin=.48,ymax=.93,
 legend pos=south west]
\addplot[blue!75!black,line width=1.3pt] table[x=G,y=root] {dimensionless-direct-metrics.dat};
\addlegendentry{Root $0$}
\addplot[orange!85!black,dashed,line width=1pt] table[x=G,y=branch] {dimensionless-direct-metrics.dat};
\addlegendentry{Branch $i$}
\addplot[green!45!black,densely dotted,line width=1.3pt] table[x=G,y=leaf] {dimensionless-direct-metrics.dat};
\addlegendentry{Leaf $i\alpha$}
\addplot[gray,densely dashed,forget plot] coordinates {(.001,.857142857142857) (100,.857142857142857)} node[pos=.98,above,anchor=south east,font=\scriptsize] {$6/7$};
\addplot[gray,densely dashed,forget plot] coordinates {(.001,.5) (100,.5)} node[pos=.98,above,anchor=south east,font=\scriptsize] {$1/2$};
\nextgroupplot[xmode=log,xmin=.001,xmax=100,
 xlabel={$G$},grid=major,grid style={gray!20},
 tick label style={font=\small},label style={font=\small},
 legend style={font=\scriptsize,draw=none},
ylabel={$\Delta a/\< a\>$},ymin=.12,ymax=.65]
\addplot[blue!75!black,line width=1.2pt] table[x=G,y=root_relative_sd] {dimensionless-direct-metrics.dat};
\addplot[orange!85!black,dashed,line width=1.2pt] table[x=G,y=branch_relative_sd] {dimensionless-direct-metrics.dat};
\addplot[green!45!black,densely dotted,line width=1.3pt] table[x=G,y=leaf_relative_sd] {dimensionless-direct-metrics.dat};
\addplot[gray,densely dashed,forget plot] coordinates {(.001,.144337567297406) (100,.144337567297406)} node[pos=.98,above,anchor=south east,font=\scriptsize] {$1/(4\sqrt3)$};
\addplot[gray,densely dashed,forget plot] coordinates {(.001,.408248290463863) (100,.408248290463863)} node[pos=.98,above,anchor=south east,font=\scriptsize] {$1/\sqrt6$};
\addplot[gray,densely dashed,forget plot] coordinates {(.001,.577350269189626) (100,.577350269189626)} node[pos=.98,above,anchor=south east,font=\scriptsize] {$1/\sqrt3$};
\nextgroupplot[grid=major,grid style={gray!20},
 tick label style={font=\small},label style={font=\small},
 legend style={font=\scriptsize,draw=none},
xmode=log,xmin=.1,xmax=10,ymin=-1.1,ymax=.1,
 xlabel={$G$},ylabel={$L\<\Scal\>$},legend pos=south west]
\addplot[blue!75!black,line width=1.1pt] table[x=G,y=root] {dimensionless-direct.dat};
\addplot[orange!85!black,dashed,line width=1.1pt] table[x=G,y=branch] {dimensionless-direct.dat};
\addplot[green!45!black,densely dotted,line width=1.3pt] table[x=G,y=leaf] {dimensionless-direct.dat};
\nextgroupplot[grid=major,grid style={gray!20},
 tick label style={font=\small},label style={font=\small},
 legend style={font=\scriptsize,draw=none},
xmin=0,xmax=10,ymin=0,ymax=6.2,xlabel={$G$},ylabel={$L\Delta\Scal$}]
\addplot[blue!75!black,line width=1.1pt] table[x=G,y=root_sd] {dimensionless-direct.dat};
\addplot[orange!85!black,dashed,line width=1.1pt] table[x=G,y=branch_sd] {dimensionless-direct.dat};
\addplot[green!45!black,densely dotted,line width=1.3pt] table[x=G,y=leaf_sd] {dimensionless-direct.dat};
\end{groupplot}
\end{tikzpicture}
\caption{Quantum gravity metric and curvature means and standard deviations with direct measure and cutoff \(0<a<L\). The lower coupling is \(G=0.001\) for metrics and \(G=0.1\) for curvatures. Dashed horizontal lines show the limiting values.}
\label{fig:direct-measure}
\end{figure}
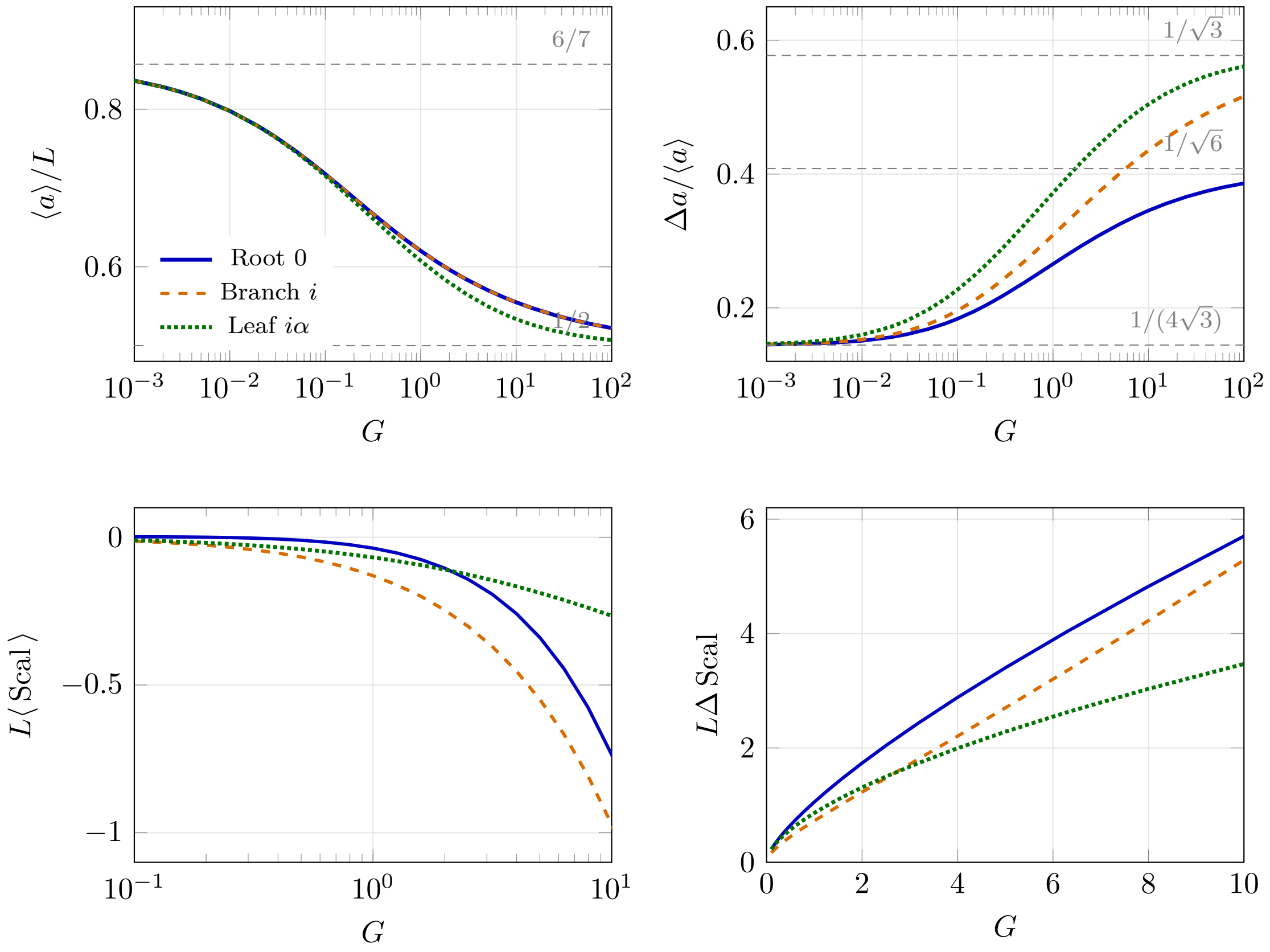

Rooted tree symmetry gives three local curvature means:
\begin{equation}\label{eq:scalar-expectations}
\begin{aligned}
 \<\Scal(0)\>
   &=\frac43\left(\left\langle\frac1{a_{11}}\right\rangle
                  -\left\langle\frac{a_2}{a_1^2}\right\rangle\right),\\
 \<\Scal(i)\>
   &=\frac12\left(\left\langle\frac1{a_1}\right\rangle
                  -\left\langle\frac{a_{11}}{a_1^2}\right\rangle\right),\\
 \<\Scal(i\alpha)\>
   &=\frac13\left(\left\langle\frac1{a_1}\right\rangle
                  -\left\langle\frac1{a_{11}}\right\rangle\right).
\end{aligned}
\end{equation}
Their standard deviations are
\(\Delta\Scal(x)=\sqrt{\<\Scal(x)^2\>-\<\Scal(x)\>^2}\).
The use of local fluctuations as a probe of quantum geometry is
motivated by \cite{BliMa:cur}. At fixed \(a_1,a_2\), the four leaf integrations factorise.
For example, their common one-dimensional kernel is
\[
 I_G(a)=\int_0^L
       \exp\!\left[-\frac1G\left(\frac at+\frac ta-2\right)\right]
       \,\mathrm{d}t.
\]
The remaining two-dimensional integral has weight
\[
 I_G(a_1)^2I_G(a_2)^2
 \exp\!\left[-\frac1G\left\{
       \left(\frac{a_1}{a_2}-1\right)^2
      +\left(\frac{a_2}{a_1}-1\right)^2\right\}\right].
\]
Inserting the corresponding leaf moments gives the local curvature
means and standard deviations shown in Figure~\ref{fig:direct-measure}.
These are moments at individual vertices, not of a sum over vertices. All direct-measure calculations use \(L=1\) numerically; the exact
scaling above restores arbitrary \(L\) in the plotted quantities.

\subsection{Quantum gravity with Liouville measure}\label{sec:liou}
Following \cite[Section 2.4]{Ma:alg4}, the changes
\(a_i=e^{\phi_i}\), \(a_{i\alpha}=e^{\phi_{i\alpha}}\) suggest the Liouville measure
\[
 \mathrm{d}^6\phi=\prod_i\frac{\mathrm{d}a_i}{a_i}\prod_{i,\alpha}\frac{\mathrm{d}a_{i\alpha}}{a_{i\alpha}}.
\]
We retain the metric cutoffs \(\delta<a_i,a_{i\alpha}<L\), as for direct measure;
the corresponding field range is
\(\log\delta<\phi_i,\phi_{i\alpha}<\log L\). The partition function is
\[
 Z^{\mathrm{Liouv}}_{\delta,L}(G)
   =\int_{(\delta,L)^6}e^{S[a]/G}\prod_i\frac{\mathrm{d}a_i}{a_i}\prod_{i,\alpha}\frac{\mathrm{d}a_{i\alpha}}{a_{i\alpha}},
\]
and expectations use the same integral with an observable inserted
and division by \(Z^{\mathrm{Liouv}}_{\delta,L}(G)\).
Scaling \(a_i=Lb_i,\ a_{i\alpha}=Lb_{i\alpha}\) preserves both action and integration measure, so
\[
 Z^{\mathrm{Liouv}}_{\delta,L}(G)
       =Z^{\mathrm{Liouv}}_{\delta/L,1}(G).
\]
Thus \(\< a(x)\>/L\), relative metric uncertainties and
\(L\Delta\Scal(x)\) depend on \(G\) and the ratio \(\delta/L\).
We specialise to reciprocal cutoffs \(\delta=1/L\), \(L>1\), in
fixed metric units. Equivalently \(-\log L<\phi_i,\phi_{i\alpha}<\log L\).
The action in these motivating field variables is
\[
\begin{split}
 S[\phi]={}&-4\cosh(\phi_1-\phi_2)
                       \bigl(\cosh(\phi_1-\phi_2)-1\bigr)-2\sum_{i,\alpha}
                  \bigl(\cosh(\phi_i-\phi_{i\alpha})-1\bigr).
\end{split}
\]
Its quadratic expansion is
\[
 S[\phi]=-2(\phi_1-\phi_2)^2
           -\sum_{i,\alpha}(\phi_i-\phi_{i\alpha})^2
           +O(\|\phi\|^4).
\]
The action and Liouville measure are invariant under common metric
rescaling, or translation of all \(\phi_i,\phi_{i\alpha}\). Without a lower metric
cutoff the common scale gives an infinite factor even when the
relative fields are restricted to a bounded region.

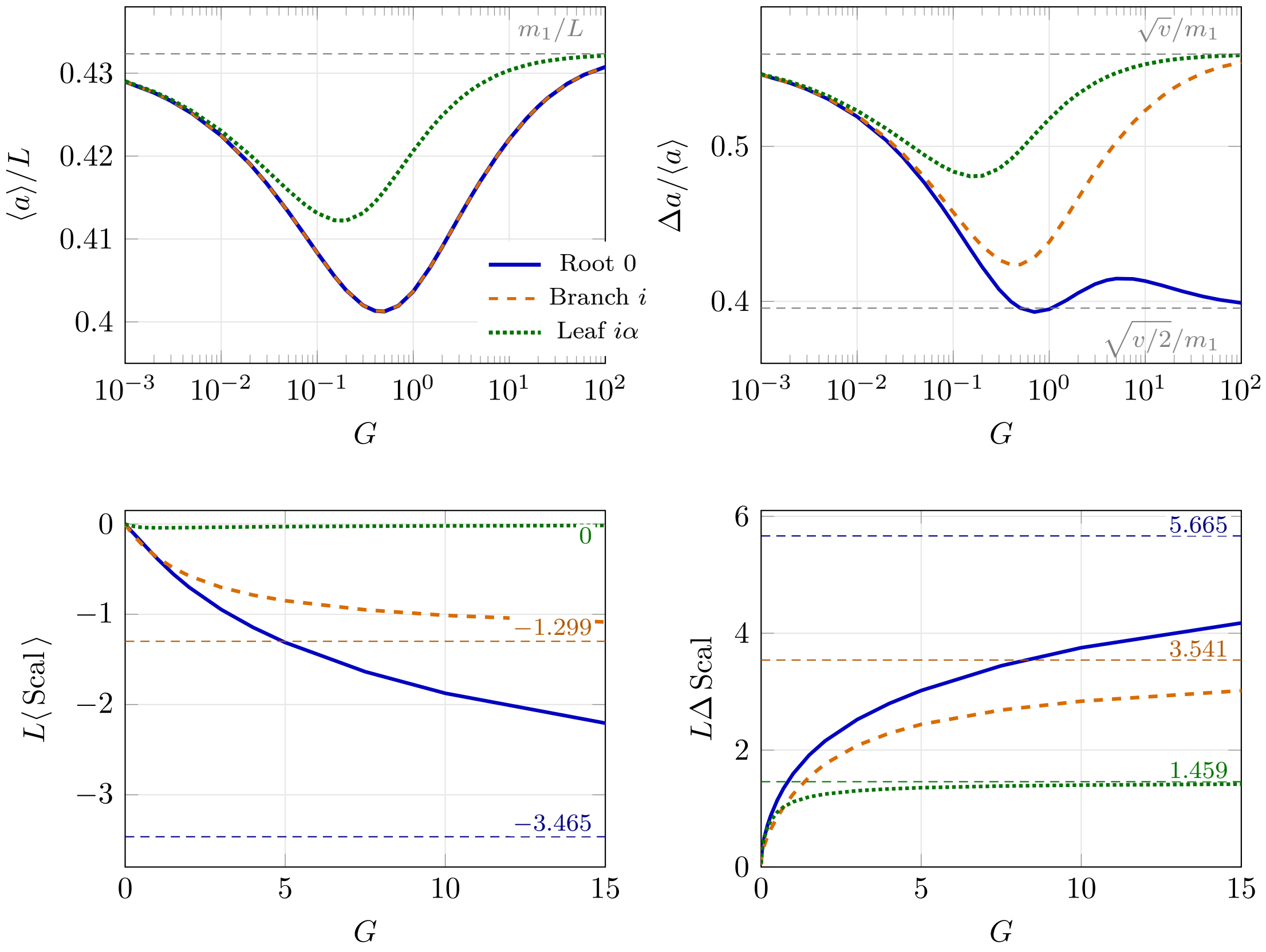
\begin{figure}[tbp]
\centering
\begin{tikzpicture}
\begin{groupplot}[
 group style={group size=2 by 2,horizontal sep=1.8cm,vertical sep=1.7cm},
 width=.44\linewidth,height=5.7cm,
]
\nextgroupplot[xmode=log,xmin=.001,xmax=100,xlabel={$G$},
 grid=major,grid style={gray!20},tick label style={font=\small},
 label style={font=\small},legend style={font=\scriptsize,draw=none},
ylabel={$\< a\>/L$},ymin=.395,ymax=.438,
 legend style={at={(1.12,.02)},anchor=south east,font=\scriptsize,draw=none}]
\addplot[blue!75!black,line width=1.3pt] table[x=G,y=root_over_L] {dimensionless-liouville-metrics.dat};
\addlegendentry{Root $0$}
\addplot[orange!85!black,dashed,line width=1pt] table[x=G,y=branch_over_L] {dimensionless-liouville-metrics.dat};
\addlegendentry{Branch $i$}
\addplot[green!45!black,densely dotted,line width=1.3pt] table[x=G,y=leaf_over_L] {dimensionless-liouville-metrics.dat};
\addlegendentry{Leaf $i\alpha$}
\addplot[gray,densely dashed,forget plot] coordinates {(.001,0.43233235838169) (100,0.43233235838169)}
 node[pos=.98,above,anchor=south east,font=\scriptsize] {$m_1/L$};
\nextgroupplot[xmode=log,xmin=.001,xmax=100,xlabel={$G$},
 grid=major,grid style={gray!20},tick label style={font=\small},
 label style={font=\small},legend style={font=\scriptsize,draw=none},
ylabel={$\Delta a/\< a\>$},ymin=.36,ymax=.59]
\addplot[blue!75!black,line width=1.2pt] table[x=G,y=root_relative_sd] {dimensionless-liouville-metrics.dat};
\addplot[orange!85!black,dashed,line width=1.2pt] table[x=G,y=branch_relative_sd] {dimensionless-liouville-metrics.dat};
\addplot[green!45!black,densely dotted,line width=1.3pt] table[x=G,y=leaf_relative_sd] {dimensionless-liouville-metrics.dat};
\addplot[gray,densely dashed,forget plot] coordinates {(.001,0.55949556343132) (100,0.55949556343132)}
 node[pos=.98,above,anchor=south east,font=\scriptsize] {$\sqrt v/m_1$};
\addplot[gray,densely dashed,forget plot] coordinates {(.001,0.39562310694608) (100,0.39562310694608)}
 node[pos=.98,below,anchor=north east,font=\scriptsize] {$\sqrt{v/2}/m_1$};
\nextgroupplot[xmin=0,xmax=15,xlabel={$G$},
 grid=major,grid style={gray!20},tick label style={font=\small},
 label style={font=\small},legend style={font=\scriptsize,draw=none},
ylabel={$L\<\Scal\>$},ymin=-3.8,ymax=.15,
 legend pos=south west]
\addplot[blue!75!black,line width=1.2pt] table[x=G,y expr={exp(1)*\thisrow{root}}] {dimensionless-liouville.dat};
\addplot[orange!85!black,dashed,line width=1.2pt] table[x=G,y expr={exp(1)*\thisrow{branch}}] {dimensionless-liouville.dat};
\addplot[green!45!black,densely dotted,line width=1.2pt] table[x=G,y expr={exp(1)*\thisrow{leaf}}] {dimensionless-liouville.dat};
\addplot[blue!50!black,densely dashed,forget plot] coordinates {(0,-3.4647008036213) (15,-3.4647008036213)} node[pos=.98,above,anchor=south east,font=\scriptsize,fill=white,inner sep=1pt] {$-3.465$};
\addplot[orange!70!black,densely dashed,forget plot] coordinates {(0,-1.299262801358) (15,-1.299262801358)} node[pos=.98,above,anchor=south east,font=\scriptsize,fill=white,inner sep=1pt] {$-1.299$};
\node[green!45!black,anchor=north east,font=\scriptsize,fill=white,inner sep=1pt] at (axis cs:14.7,0) {$0$};
\nextgroupplot[xmin=0,xmax=15,xlabel={$G$},
 grid=major,grid style={gray!20},tick label style={font=\small},
 label style={font=\small},legend style={font=\scriptsize,draw=none},
ylabel={$L\Delta\Scal$},ymin=0,ymax=6.1]
\addplot[blue!75!black,line width=1.2pt] table[x=G,y=root_Lsd] {dimensionless-liouville.dat};
\addplot[orange!85!black,dashed,line width=1.2pt] table[x=G,y=branch_Lsd] {dimensionless-liouville.dat};
\addplot[green!45!black,densely dotted,line width=1.2pt] table[x=G,y=leaf_Lsd] {dimensionless-liouville.dat};
\addplot[blue!50!black,densely dashed,forget plot] coordinates {(0,5.6648153461862) (15,5.6648153461862)} node[pos=.98,above,anchor=south east,font=\scriptsize,fill=white,inner sep=1pt] {$5.665$};
\addplot[orange!70!black,densely dashed,forget plot] coordinates {(0,3.5410388291581) (15,3.5410388291581)} node[pos=.98,above,anchor=south east,font=\scriptsize,fill=white,inner sep=1pt] {$3.541$};
\addplot[green!45!black,densely dashed,forget plot] coordinates {(0,1.4593441562257) (15,1.4593441562257)} node[pos=.98,above,anchor=south east,font=\scriptsize,fill=white,inner sep=1pt] {$1.459$};
\end{groupplot}
\end{tikzpicture}
\caption{Quantum gravity metric and curvature means and standard deviations for Liouville measure and cutoff  \(e^{-1}<a<e\). The lower coupling is \(G=0.001\). Dashed horizontal lines show limiting values, with the curvature asymptotes  approached at much larger \(G\).}
\label{fig:liouville-measure}
\end{figure}

Inverting all six parameters \(a_i,a_{i\alpha}\) also preserves the action and the
reciprocal-cutoff domain, giving the exact duality
\[
 \< a^p\>=\< a^{-p}\>
 \qquad(p\in\R),
\]
for either \(a=a_i\) or \(a=a_{i\alpha}\).
At weak coupling the metric parameters concentrate at a common value,
distributed with measure \(\mathrm{d}a/a\) on \((1/L,L)\).
Consequently
\[
 \< a\>\longrightarrow
       \frac{L-L^{-1}}{2\log L}\qquad(G\downarrow0).
\]
For comparison with metric observables, attach the local outward
square-length \(a(i)=a_i\) to a branch and
\(a(i\alpha)=a_{i\alpha}\) to a leaf. At the root, which has no
incoming outward arrow, use the symmetric mean
\[
 a(0)=\frac{a_1+a_2}{2}.
\]
There are only two independent edge-parameter means:
\(\< a(0)\>=\< a_i\>\) by symmetry, while
\(\< a_{i\alpha}\>\) can differ. The root standard
deviation is nevertheless distinct because it includes the correlation
between the two root edges:
\[
 \bigl(\Delta a(0)\bigr)^2
  =\frac12\bigl(\< a_1^2\>+\< a_1a_2\>\bigr)
       -\< a_1\>^2.
\]
Put
\[
 m_1=\frac{L-L^{-1}}{2\log L},\qquad
 v=\frac{L^2-L^{-2}}{4\log L}-m_1^2.
\]
All three means tend to \(m_1\) at both weak and strong coupling.
All three standard deviations tend to \(\sqrt v\) as
\(G\downarrow0\), when the edge parameters become equal. At
\(G\to\infty\), independence instead gives
\[
 \Delta a(0)\longrightarrow\sqrt{v/2},\qquad
 \Delta a_i,\ \Delta a_{i\alpha}\longrightarrow\sqrt v.
\]
We report the means as \(\< a(x)\>/L\) and the
uncertainties as \(\Delta a(x)/\< a(x)\>\). The normalised
mean tends at both weak and strong coupling to
\[
 \frac{m_1}{L}=\frac{1-L^{-2}}{2\log L}.
\]
The relative standard deviations tend to \(\sqrt v/m_1\) at weak
coupling; their strong-coupling limits are \(\sqrt{v/2}/m_1\) at the
root and \(\sqrt v/m_1\) at a branch or leaf.
Figure~\ref{fig:liouville-measure} illustrates the intervening behaviour
at cutoff ratio \(\delta/L=e^{-2}\), represented by \(L=e\) and
\(\delta=e^{-1}\) in the reciprocal convention.

The curvature standard deviations start at zero and approach the
strong-coupling limits of independent logarithmically distributed
metric parameters. For cutoff ratio \(\delta/L=e^{-2}\),  Figure~\ref{fig:liouville-measure}
shows  the root, branch and leaf values  \(L \< \Scal\>\) and \(L\Delta\Scal\) and their limiting values. All fixed-cutoff
Euclidean moments are analytic in \(1/G\) for \(G>0\); the curves
show smooth crossovers, rather than a finite-cutoff phase transition.

\subsubsection{Large reciprocal cutoffs}

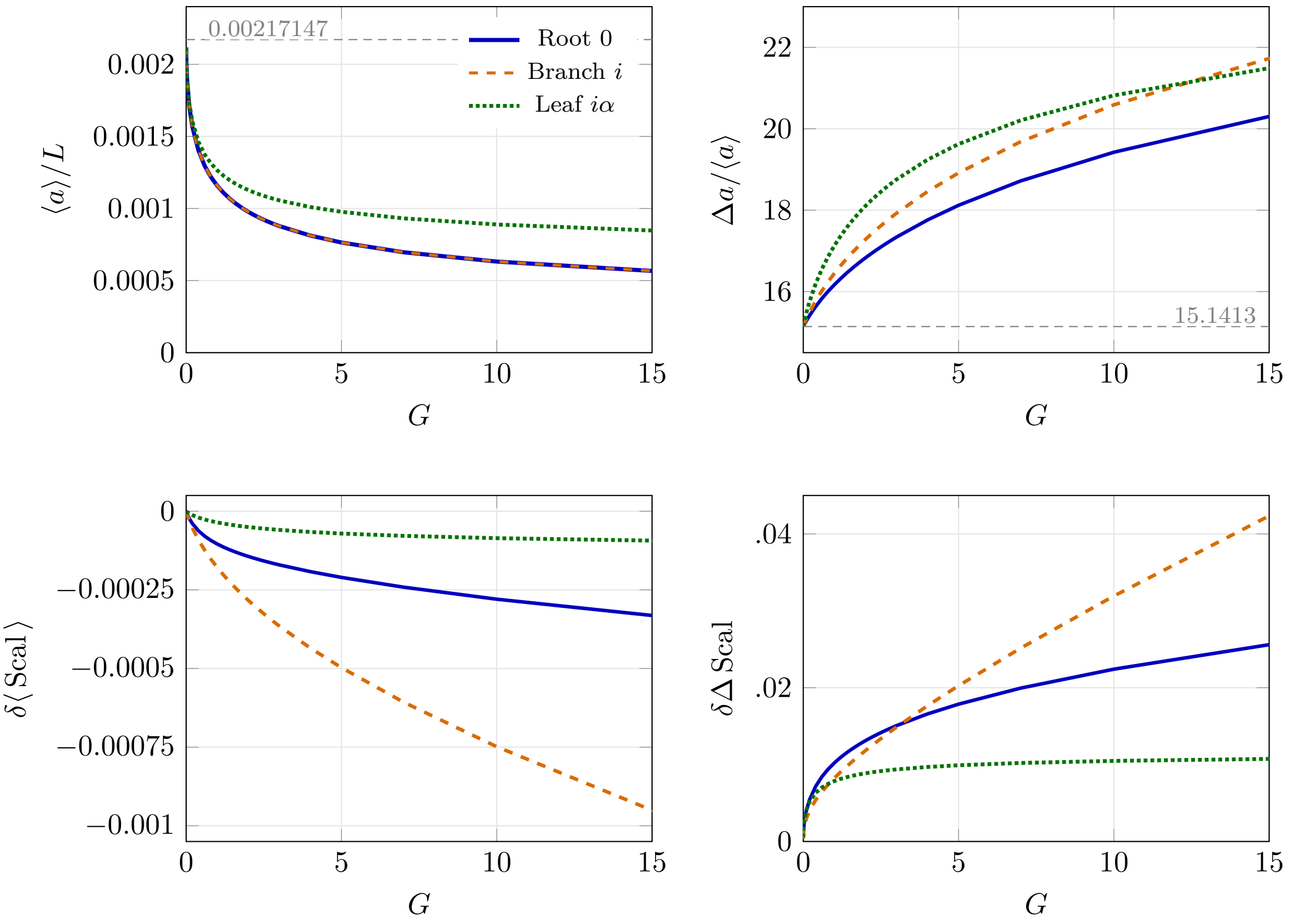
\begin{figure}[tbp]
\centering
\begin{tikzpicture}
\begin{groupplot}[
 group style={group size=2 by 2,horizontal sep=1.8cm,vertical sep=1.7cm},
 width=.44\linewidth,height=5.7cm,
]
\nextgroupplot[xmin=0,xmax=15,xlabel={$G$},
 grid=major,grid style={gray!20},tick label style={font=\small},
 label style={font=\small},legend style={font=\scriptsize,draw=none},
ylabel={$\< a\>/L$},ymin=0,ymax=.0024,
 scaled y ticks=false,ytick={0,.0005,.001,.0015,.002},
 yticklabel style={/pgf/number format/fixed,/pgf/number format/precision=4},
 legend pos=north east]
\addplot[blue!75!black,line width=1.4pt] table[x=G,y=mean_branch] {wide-liouville-metrics.dat};
\addlegendentry{Root $0$}
\addplot[orange!85!black,dashed,line width=1pt] table[x=G,y=mean_branch] {wide-liouville-metrics.dat};
\addlegendentry{Branch $i$}
\addplot[green!45!black,densely dotted,line width=1.3pt] table[x=G,y=mean_leaf] {wide-liouville-metrics.dat};
\addlegendentry{Leaf $i\alpha$}
\addplot[gray,densely dashed,forget plot] coordinates {(0,.00217147240951626) (15,.00217147240951626)} node[pos=.04,above,anchor=south west,font=\scriptsize,fill=white,inner sep=1pt] {$0.00217147$};
\nextgroupplot[xmin=0,xmax=15,xlabel={$G$},
 grid=major,grid style={gray!20},tick label style={font=\small},
 label style={font=\small},legend style={font=\scriptsize,draw=none},
ylabel={$\Delta a/\< a\>$},ymin=14.5,ymax=23]
\addplot[blue!75!black,line width=1.2pt] table[x=G,y=relative_sd_root] {wide-liouville-metrics.dat};
\addplot[orange!85!black,dashed,line width=1.2pt] table[x=G,y=relative_sd_branch] {wide-liouville-metrics.dat};
\addplot[green!45!black,densely dotted,line width=1.3pt] table[x=G,y=relative_sd_leaf] {wide-liouville-metrics.dat};
\addplot[gray,densely dashed,forget plot] coordinates {(0,15.1412849289419) (15,15.1412849289419)} node[pos=.98,above,anchor=south east,font=\scriptsize,fill=white,inner sep=1pt] {$15.1413$};
\nextgroupplot[xmin=0,xmax=15,xlabel={$G$},
 grid=major,grid style={gray!20},tick label style={font=\small},
 label style={font=\small},legend style={font=\scriptsize,draw=none},
ylabel={$\delta\<\Scal\>$},ymin=-.00105,ymax=.00005,
 scaled y ticks=false,ytick={-.001,-.00075,-.0005,-.00025,0},
 yticklabel style={/pgf/number format/fixed,/pgf/number format/precision=5},
 legend pos=south west]
\addplot[blue!75!black,line width=1.2pt] table[x=G,y=delta_scal_root] {wide-liouville-curvature.dat};
\addplot[orange!85!black,dashed,line width=1.2pt] table[x=G,y=delta_scal_branch] {wide-liouville-curvature.dat};
\addplot[green!45!black,densely dotted,line width=1.3pt] table[x=G,y=delta_scal_leaf] {wide-liouville-curvature.dat};
\nextgroupplot[xmin=0,xmax=15,xlabel={$G$},
 grid=major,grid style={gray!20},tick label style={font=\small},
 label style={font=\small},legend style={font=\scriptsize,draw=none},
ylabel={$\delta\Delta\Scal$},ymin=0,ymax=.045,
 scaled y ticks=false,ytick={0,.02,.04},yticklabels={$0$,$.02$,$.04$}]
\addplot[blue!75!black,line width=1.2pt] table[x=G,y=delta_sd_root] {wide-liouville-curvature.dat};
\addplot[orange!85!black,dashed,line width=1.2pt] table[x=G,y=delta_sd_branch] {wide-liouville-curvature.dat};
\addplot[green!45!black,densely dotted,line width=1.3pt] table[x=G,y=delta_sd_leaf] {wide-liouville-curvature.dat};
\end{groupplot}
\end{tikzpicture}
\caption{Results as in Figure~\ref{fig:liouville-measure} for the Liouville measure but with \(10^{-100}<a<10^{100}\). Curvatures are now scaled by \(\delta=1/L\). Dashed horizontal lines show weak-coupling limits; strong coupling lies far beyond the range shown.}
\label{fig:wide-liouville}
\end{figure}

To probe a much larger range of square-lengths, namely $L=10^{100}$ and
$\delta=L^{-1}$, we need new methods as direct numerical integration can be unstable.
The results in Figure~\ref{fig:wide-liouville} show a substantial change from
the narrow range treated before: the metric means are much smaller relative to $L$
and the relative metric fluctuations much larger. The curves remain smooth
throughout the scanned range $0.001\le G\le100$, with no abrupt change in these observables.

Our approach to obtaining these results is to divide every metric parameter by $L$, so the metric 
integration domain is $(e^{-T},1)^6$, where $T=2\log L=200\log10$.
The following separation of the common scale then avoids numerical sampling
across $200$ decades. Write the rescaled parameters as
\[
 b_1=e^A,\qquad b_2=e^{A+d},\qquad
 b_{1\alpha}=e^{A+v_{1\alpha}},\qquad
 b_{2\alpha}=e^{A+d+v_{2\alpha}}.
\]
The Liouville measure is $\mathrm{d}A\,\mathrm{d}d\prod_{i,\alpha}\mathrm{d}v_{i\alpha}$,
and the positive energy $E=-S$ is independent of $A$:
\[
 E(d,v)=(e^d-1)^2+(e^{-d}-1)^2
           +2\sum_{i,\alpha}(\cosh v_{i\alpha}-1).
\]
Let $h=(0,d,v_{11},v_{12},d+v_{21},d+v_{22})$, and put
$h_-=\min h$, $h_+=\max h$, $r=h_+-h_-$.
The six cutoff conditions give precisely
\[
 -T-h_-<A<-h_+.
\]
Hence the common-scale integral is elementary and gives the exact formula
\[
 Z^{\mathrm{Liouv}}_{e^{-T},1}(G)
  =\int_{r<T}(T-r)e^{-E(d,v)/G}\,
                 \mathrm{d}d\prod_{i,\alpha}\mathrm{d}v_{i\alpha}.
\]
For a metric observable $Q$ homogeneous of degree $p>0$, such as a local
square-length or its square, the corresponding numerator is
\[
 \int_{r<T}Q(e^h)\,
       \frac{e^{-p h_+}-e^{-p(T+h_-)}}{p}\,
       e^{-E(d,v)/G}\,\mathrm{d}d\prod_{i,\alpha}\mathrm{d}v_{i\alpha}.
\]
These formulas give both the means and the second moments needed for
relative standard deviations. At the root $Q(b)=(b_1+b_2)/2$; at a
branch or leaf it is the corresponding parameter.

For numerical evaluation, the four leaf integrals are reduced to
one-dimensional conditional moments at fixed parent values. The
relative fields are confined by $E$, even for very large $T$.
After truncating their exponentially suppressed tails, the conditional
moments are independent of $A$ throughout the interior of the
common-scale interval, where its integral is evaluated analytically.
Numerical quadrature is needed only for the relative fields and the
finite boundary layers at the cutoffs. Increasing $T$ therefore does
not require a proportionally finer grid. 
The numerical checks include variation of the relative-field tail cutoff.

The scale integral also controls the large-cutoff limit at fixed $G$.
Dominated convergence gives
\[
 \lim_{T\to\infty}\frac{Z^{\mathrm{Liouv}}_{e^{-T},1}(G)}{T}
 =\left[\int_{\R}e^{-[(e^d-1)^2+(e^{-d}-1)^2]/G}\,\mathrm{d}d\right]
  \left[\int_{\R}e^{-2(\cosh v-1)/G}\,\mathrm{d}v\right]^4.
\]
Both factors are finite, positive and analytic for every $G>0$:
the exponential growth of the energy confines all relative fields.
Thus the partition function per unit logarithmic scale has no
finite-positive-$G$ singularity in this model. The metric numerators
above have finite positive limits, giving $\< a\>/L$
proportional to $T^{-1}$ and relative fluctuations proportional to
$\sqrt T$ at fixed $G$. Their limiting coefficients are also smooth
for $G>0$, since $Q(e^h)e^{-p h_+}$ is bounded for these observables.
These cutoff-dependent limits differ from
first taking $G\to\infty$ at fixed $T$. In particular, the return of
the means to their strong-coupling values in the narrow-cutoff example
need not be visible over any comparable range of $G$ at large $T$.

Scalar curvature can be treated by the same scale integration, now
using homogeneity degrees $-1$ and $-2$ for its first and second
moments. For the reciprocal cutoffs, the useful finite quantities are
$\delta\<\Scal\>=\<\Scal\>/L$ and
$\delta\Delta\Scal=\Delta\Scal/L$.
Indeed, $\Scal(a)=L^{-1}\Scal(b)$, so their evaluation on the
normalised domain uses $e^{-T}\Scal(b)$ and
$e^{-2T}\Scal(b)^2$. The required common-scale integrals are
\[
 e^{-T}\int_{-T-h_-}^{-h_+}e^{-A}\,\mathrm{d}A
       =e^{h_-}-e^{-T+h_+},\qquad
 e^{-2T}\int_{-T-h_-}^{-h_+}e^{-2A}\,\mathrm{d}A
       =\frac{e^{2h_-}-e^{-2T+2h_+}}{2}.
\]
Multiplying these factors by $\Scal(x;e^h)$ and
$\Scal(x;e^h)^2$, respectively, and integrating over $r<T$ with
weight $e^{-E/G}$ gives the two moment numerators. Here
$\Scal(x;e^h)$ denotes the scalar curvature at vertex $x$ evaluated
on the six parameters $e^h$. Division by $Z$ and subtraction of the
squared mean then give
\[
 \delta\Delta\Scal(x)
 =\sqrt{\delta^2\<\Scal(x)^2\>
                    -\bigl(\delta\<\Scal(x)\>\bigr)^2}.
\]

Figure~\ref{fig:wide-liouville} shows negative curvature
means of increasing magnitude and growing curvature fluctuations over
the displayed range. All three vertex types vary smoothly, and the
scan to $G=100$ reveals no abrupt change.
For fixed $G>0$, the scaled first- and second-moment numerators above
have finite limits as $T\to\infty$. Their coefficients are smooth in
$G$, since the confining energy dominates the exponential factors
from the curvature. Together with $Z$ proportional to $T$, this gives
\[
 \delta\<\Scal(x)\>=O(T^{-1}),\qquad
 \delta\Delta\Scal(x)=O(T^{-1/2}).
\]
In reciprocal units the unscaled curvature means are therefore
$O(L/\log L)$, while the standard deviations grow as
$L/\sqrt{\log L}$ with nonzero, vertex-dependent coefficients
at fixed $G$.
Thus smooth dependence on the coupling coexists with divergent
curvature fluctuations as the cutoffs are removed.

\subsection{Quantum gravity with Lorentzian action}\label{sec:lore}
An oscillatory action due to a $\imath{}$ in front of the Einstein-Hilbert action is needed for the Lorentzian case but can also be used for positive metric parameters, as in the finite Euclidean fibre contribution to a Lorentzian
Kaluza--Klein model \cite{LiuMa:kal}. We look at this case here. 

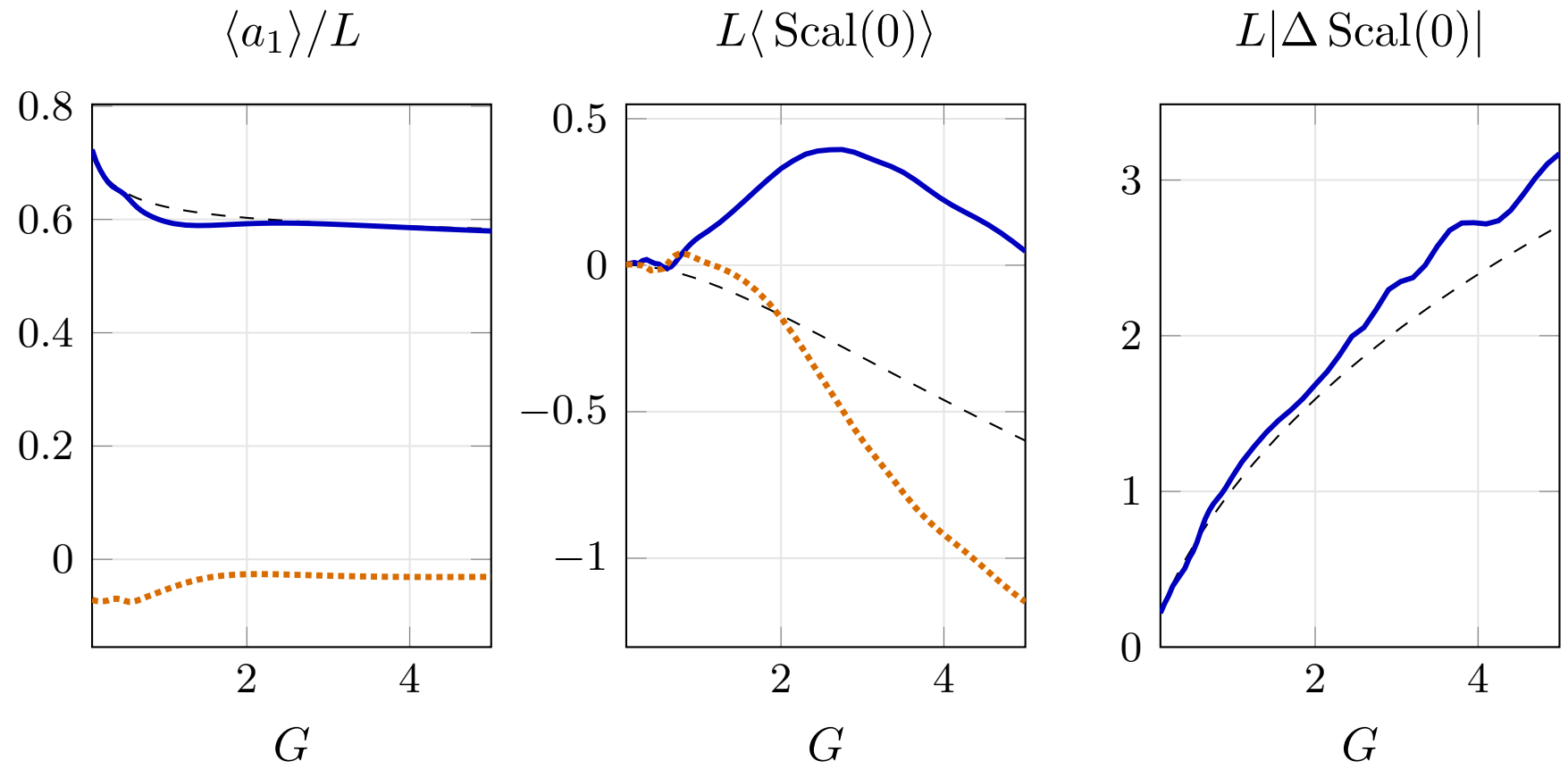
\begin{figure}[tp]
\centering
\begin{tikzpicture}
\begin{groupplot}[
 group style={group size=3 by 1,horizontal sep=1cm},
 width=.28\linewidth,height=5.6cm,xmin=.1,xmax=5,xlabel={$G$},
 grid=major,grid style={gray!20},tick label style={font=\footnotesize},
 title style={font=\small},label style={font=\small},
]
\nextgroupplot[title={$\< a_1\>/L$}]
\addplot[black,dashed] table[x=G,y=a_root] {dimensionless-comparison.dat};
\addplot[blue!75!black,line width=1.1pt] table[x=G,y=aOre] {dimensionless-comparison.dat};
\addplot[orange!85!black,densely dotted,line width=1.3pt] table[x=G,y=aOim] {dimensionless-comparison.dat};
\nextgroupplot[title={$L\<\Scal(0)\>$}]
\addplot[black,dashed] table[x=G,y=root] {dimensionless-comparison.dat};
\addplot[blue!75!black,line width=1.1pt] table[x=G,y=rootOre] {dimensionless-comparison.dat};
\addplot[orange!85!black,densely dotted,line width=1.3pt] table[x=G,y=rootOim] {dimensionless-comparison.dat};
\nextgroupplot[title={$L|\Delta\Scal(0)|$},ymin=0]
\addplot[black,dashed] table[x=G,y=root_sd] {dimensionless-comparison.dat};
\addplot[blue!75!black,line width=1.1pt] table[x=G,y=rootOsd] {dimensionless-comparison.dat};
\end{groupplot}
\end{tikzpicture}
\caption{Lorentzian theory with direct measure and
\(L/10<a_i,a_{i\alpha}<L\). Black dashed curves are for the Euclidean theory; blue solid curves
show the Lorentzian theory real parts or fluctuation magnitudes. Orange
dotted curves show imaginary parts. The plots extend down to $G=0.1$.}
\label{fig:dimensionless-oscillatory}
\end{figure}

With direct measure and
\(0<\delta<a_i,a_{i\alpha}<L\), we use the regularised weight
\[
 e^{-\imath{}S[a]/G_\eta},\qquad
 G_\eta=G-\imath{}\eta,\quad\eta>0,
\]
but with the dimensionless action. The damping factor is again
\(\exp(\eta S/(G^2+\eta^2))\le1\). Both \(G\) and \(\eta\) are
now dimensionless; after scaling by \(L\), the independent parameters
are \(G,\eta,\delta/L\). For a complex observable, use
\begin{equation}\label{eq:complex-standard-deviation}
 \Delta O=\sqrt{\< O^*O\>
                         -\< O^*\>\< O\>},
 \qquad
 |\Delta O|=\sqrt{\left|\< O^*O\>
                         -\< O^*\>\< O\>\right|}.
\end{equation}
The conjugation is applied before integration; with a complex weight,
\(\< O^*\>\) need not be \(\overline{\< O\>}\).
For real scalar curvature, this gives
\( |\Delta\Scal(x)|=\sqrt{|\<\Scal(x)^2\>-\<\Scal(x)\>^2|}\).

Figure~\ref{fig:dimensionless-oscillatory} shows the results for the range \(\delta/L=1/10\), \(L=1\) and \(\eta=10^{-8}\), including comparison with the previous Euclidean version. Smaller lower cutoffs produce more rapid phase oscillations, making curvature fluctuations particularly difficult to resolve numerically to the stated accuracy, requiring substantially finer quadrature or integration methods tailored to these oscillations.

\subsection{Quantum gravity with mixed signature}

So far we have looked at the $a_i,a_{i\alpha}$ all positive, both with or without $\imath{}$ in the action (what we term the Euclidean and Lorentzian forms of the action). Now we allow some of the metric values to become negative, which corresponds to a mixed signature or Lorentzian form of the geometry itself.  Specifically, we choose one positive and one negative outward edge at each branching:
\begin{equation}\label{eq:mixed-signature}
 a_1=x_1,\qquad a_2=-x_2,\qquad
 a_{i1}=x_{i1},\qquad a_{i2}=-x_{i2},\qquad x_i,x_{i\alpha}>0.
\end{equation}
The reverse arrows have the same signs because
\(\eps=1\). Let \(S[x]\) be the previous action as a function of the $x_i,x_{i\alpha}$.  We integrate over $0<\delta<x_i,x_{i\alpha}<L$ much as before.

Note that exchanging the two child choices at every branching
defines an operation
\[
 (Tx)_i=x_{\bar i},\qquad
 (Tx)_{i\alpha}=x_{\bar i\bar\alpha}.
\]
On the metric itself, this is a tree relabelling followed by overall
sign reversal. Since \(S[-a]=S[a]\), we therefore have the symmetry
\[
 S[Tx]=S[x].
\]
The scalar curvature still changes sign under the overall metric
reversal. For equal magnitude integration cutoffs, \(T\) therefore gives
\[
 \<\Scal(0)\>=0,\qquad
 \<\Scal(i)\>=-\<\Scal(\bar i)\>,\qquad
 \<\Scal(i\alpha)\>
       =-\<\Scal(\bar i\bar\alpha)\>.
\]
These identities hold for both Euclidean and Lorentzian forms of action, including the
regularised complex coupling, wherever expectations exist, as they
involve no complex conjugation. In particular, the root mean in the Lorentzian action case vanishes. 

Set $\beta=1/G$ for the Euclidean action and
$\beta=-\imath{}/(G-\imath{}\eta)$ for the regularised Lorentzian action. Put
\[
 I_\beta(x)=\int_\delta^L e^{-\beta(t/x+x/t)}\,\mathrm{d}t,
 \qquad H_\beta(x)=I_\beta(x)I_{-\beta}(x).
\]
The leaf integrations give the exact reduced partition function
\[
 Z_{\delta,L}(\beta)
 =e^{8\beta}\int_\delta^L\!\int_\delta^L
 H_\beta(x)H_\beta(y)
 e^{-\beta\{(1+x/y)^2+(1+y/x)^2\}}\,\mathrm{d}x\,\mathrm{d}y.
\]
The root ratio couples the positive and negative sectors, so there
is still no factorisation into independent three-variable theories. At real Lorentzian coupling $\beta=-\imath{}/G$, the leaf integrals
pair into absolute squares, so $H_\beta$ is nonnegative. However, the
coupling between the two root edges retains a nonconstant complex
phase, so the full partition function is generally complex. For comparison, the metric-independent-measure action of
Proposition~\ref{prop:EH-action} is odd under $T$. At real Lorentzian
coupling its partition function is therefore real and its root
curvature mean generally purely imaginary, instead of zero.

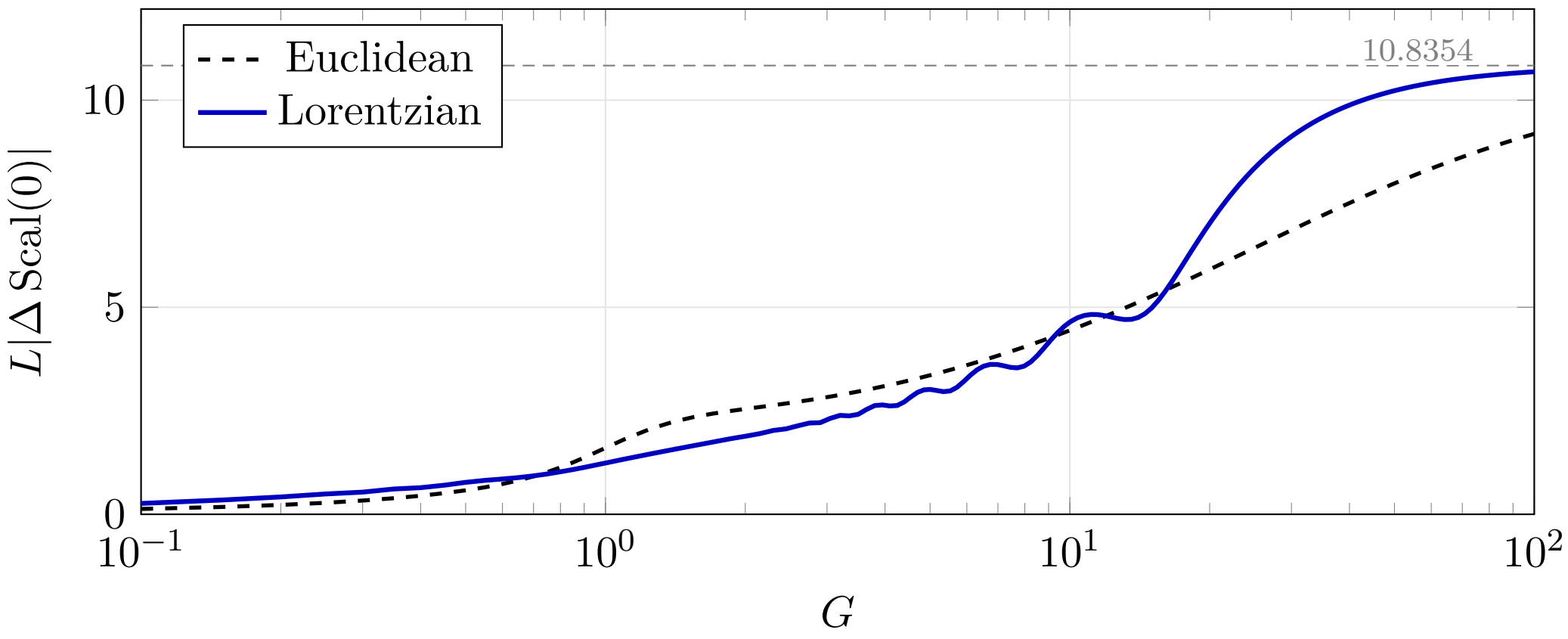
\begin{figure}[tbp]
\centering
\begin{tikzpicture}
\begin{axis}[
 width=.85\linewidth,height=6cm,xmode=log,xmin=.1,xmax=100,
 xlabel={$G$},ylabel={$L|\Delta\Scal(0)|$},ymin=0,ymax=12.2,
 grid=major,grid style={gray!20},legend pos=north west,
]
\addplot[black,dashed,line width=1pt] table[x=G,y=root_sd] {mixed-signature-decade.dat};
\addlegendentry{Euclidean}
\addplot[blue!75!black,line width=1.2pt] table[x=G,y=rootOsd] {mixed-signature-decade.dat};
\addlegendentry{Lorentzian}
\addplot[gray,densely dashed,forget plot] coordinates {(.1,10.83540165551) (100,10.83540165551)}
 node[pos=.96,above,anchor=south east,font=\scriptsize,fill=white,inner sep=1pt] {$10.8354$};
\end{axis}
\end{tikzpicture}
\caption{Mixed-signature root scalar-curvature fluctuations for
Euclidean and Lorentzian actions, with $L/10<x_i,x_{i\alpha}<L$
and $\eta=10^{-8}$ in the Lorentzian case. Scans starts at $G=0.1$. The dashed horizontal
line marks the common strong-coupling limit.}
\label{fig:dimensionless-mixed}
\end{figure}

Both the Euclidean theory and the adopted regularisation of the
Lorentzian action require a lower cutoff $\delta>0$. A leaf of sign opposite
to its parent can approach zero with the other parameters fixed,
making \(S[x]\to\infty\). The exponential in the Euclidean action case,
and the regularising factor when \(\eta>0\) in the Lorentzian one, then grow without bound.
Finite positive cutoffs make all the integrals well-defined.

\begin{figure}[tbp]
\centering
\begin{tikzpicture}
\begin{groupplot}[
 group style={group size=3 by 1,horizontal sep=1cm},
 width=.28\linewidth,height=5.6cm,xmode=log,xmin=.1,xmax=100,xlabel={$G$},
 grid=major,grid style={gray!20},tick label style={font=\footnotesize},
 title style={font=\small},label style={font=\small},
 legend style={font=\scriptsize,draw=none},
]
\nextgroupplot[title={$\< x\>/L$},ymin=0,ymax=1.05,
 legend style={at={(1.24,1.02)},anchor=north east,font=\scriptsize,draw=none}]
\addplot[blue!75!black,line width=1.2pt] table[x=G,y=parent_mean] {mixed-signature-decade.dat};
\addlegendentry{$x_i$}
\addplot[orange!85!black,dashed,line width=1.2pt] table[x=G,y=aligned_mean] {mixed-signature-decade.dat};
\addlegendentry{$x_{11},x_{22}$}
\addplot[green!45!black,densely dotted,line width=1.3pt] table[x=G,y=opposed_mean] {mixed-signature-decade.dat};
\addlegendentry{$x_{12},x_{21}$}
\nextgroupplot[title={$\Delta x/\< x\>$},ymin=0,ymax=1.5]
\addplot[blue!75!black,line width=1.2pt] table[x=G,y=parent_relative_sd] {mixed-signature-decade.dat};
\addplot[orange!85!black,dashed,line width=1.2pt] table[x=G,y=aligned_relative_sd] {mixed-signature-decade.dat};
\addplot[green!45!black,densely dotted,line width=1.3pt] table[x=G,y=opposed_relative_sd] {mixed-signature-decade.dat};
\nextgroupplot[title={$L\<\Scal\>$},
 legend style={at={(1.24,1.02)},anchor=north east,font=\scriptsize,draw=none}]
\addplot[blue!75!black,line width=1.2pt] table[x=G,y=branch_mean] {mixed-signature-decade.dat};
\addlegendentry{Branch $1$}
\addplot[orange!85!black,dashed,line width=1.2pt] table[x=G,y=leaf11_mean] {mixed-signature-decade.dat};
\addlegendentry{Leaf $11$}
\addplot[green!45!black,densely dotted,line width=1.3pt] table[x=G,y=leaf12_mean] {mixed-signature-decade.dat};
\addlegendentry{Leaf $12$}
\end{groupplot}
\end{tikzpicture}
\caption{Euclidean quantum gravity expectations with mixed signature and
$L/10<x_i,x_{i\alpha}<L$. The curvature mean $\< \Scal\>$ is zero at the root and at the vertices $2,22,21$ is the negative of the values shown at $1,11,12$, respectively. The scans  start at $G=0.1$.}
\label{fig:mixed-euclidean-observables}
\end{figure}
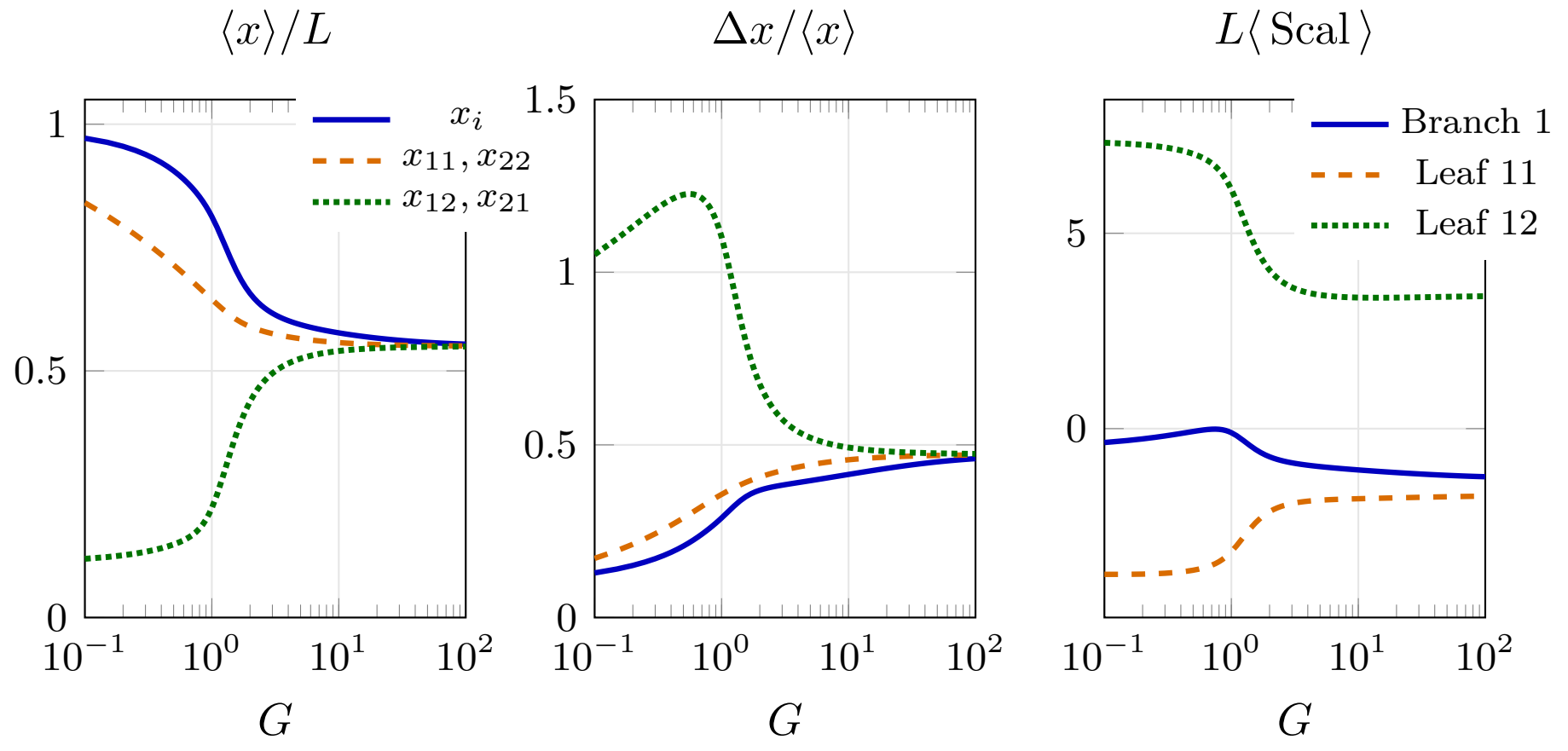

Although the root mean curvature $\<\Scal(0)\>$ vanishes,
its fluctuation magnitude need not. Figure~\ref{fig:dimensionless-mixed}
compares the Euclidean and Lorentzian actions on the decade-wide
range $L/10<x_i,x_{i\alpha}<L$. Both tend at strong coupling to
the common uniform-measure limit marked in the figure.
The Lorentzian curve approaches this limit more rapidly at large $G$.

There are three classes of metric magnitudes: $x_i$ on the root
edges, $x_{11},x_{22}$ on leaf edges whose sign matches that of
their parent edge, and $x_{12},x_{21}$ on those of opposite sign.
The members of each class have the same magnitude statistics by $T$.
Figure~\ref{fig:mixed-euclidean-observables} shows their Euclidean
means and relative uncertainties, together with the independent
curvature means. Means of the signed metric parameters follow by
inserting the signs in \eqref{eq:mixed-signature}.
At strong coupling the direct measure becomes uniform, giving
\[
 \frac{\< x_i\>}{L},\quad
 \frac{\< x_{i\alpha}\>}{L}\longrightarrow\frac{11}{20}=0.55,
 \qquad G\to\infty.
\]
The relative uncertainties tend to $9/(11\sqrt3)\approx0.4724$.

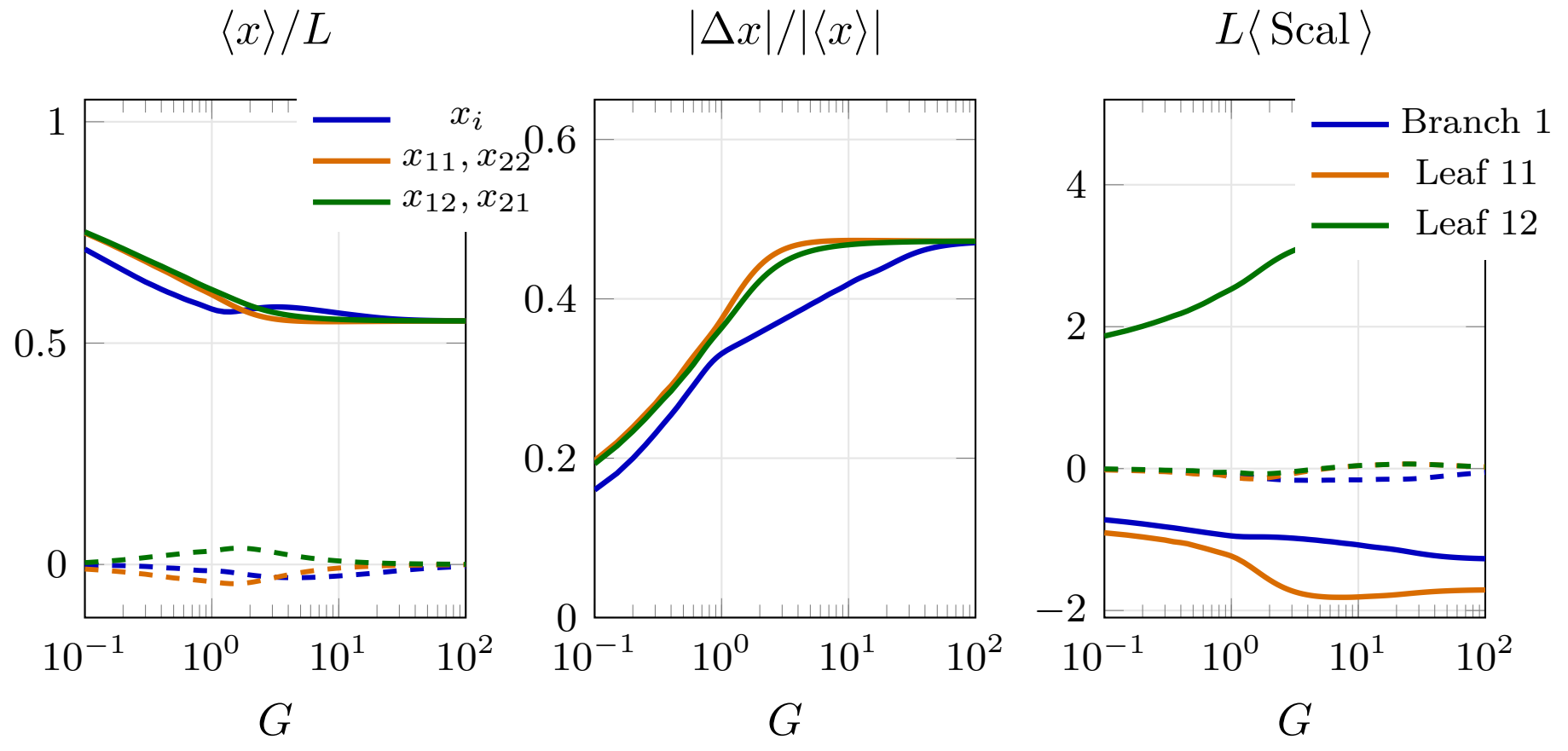
\begin{figure}[tbp]
\centering
\begin{tikzpicture}
\begin{groupplot}[
 group style={group size=3 by 1,horizontal sep=1cm},
 width=.28\linewidth,height=5.6cm,xmode=log,xmin=.1,xmax=100,xlabel={$G$},
 grid=major,grid style={gray!20},tick label style={font=\footnotesize},
 title style={font=\small},label style={font=\small},
 legend style={font=\scriptsize,draw=none},
]
\nextgroupplot[title={$\< x\>/L$},ymin=-.12,ymax=1.05,
 legend style={at={(1.24,1.02)},anchor=north east,font=\scriptsize,draw=none}]
\addplot[blue!75!black,line width=1.2pt] table[x=G,y=parent_re] {mixed-signature-lorentzian.dat};
\addlegendentry{$x_i$}
\addplot[orange!85!black,line width=1.2pt] table[x=G,y=aligned_re] {mixed-signature-lorentzian.dat};
\addlegendentry{$x_{11},x_{22}$}
\addplot[green!45!black,line width=1.2pt] table[x=G,y=opposed_re] {mixed-signature-lorentzian.dat};
\addlegendentry{$x_{12},x_{21}$}
\addplot[blue!75!black,dashed,forget plot,line width=1.1pt] table[x=G,y=parent_im] {mixed-signature-lorentzian.dat};
\addplot[orange!85!black,dashed,forget plot,line width=1.1pt] table[x=G,y=aligned_im] {mixed-signature-lorentzian.dat};
\addplot[green!45!black,dashed,forget plot,line width=1.1pt] table[x=G,y=opposed_im] {mixed-signature-lorentzian.dat};
\nextgroupplot[title={$|\Delta x|/|\< x\>|$},ymin=0,ymax=.65]
\addplot[blue!75!black,line width=1.2pt] table[x=G,y=parent_relative_sd] {mixed-signature-lorentzian.dat};
\addplot[orange!85!black,line width=1.2pt] table[x=G,y=aligned_relative_sd] {mixed-signature-lorentzian.dat};
\addplot[green!45!black,line width=1.2pt] table[x=G,y=opposed_relative_sd] {mixed-signature-lorentzian.dat};
\nextgroupplot[title={$L\<\Scal\>$},ymin=-2.1,ymax=5.2,
 legend style={at={(1.24,1.02)},anchor=north east,font=\scriptsize,draw=none}]
\addplot[blue!75!black,line width=1.2pt] table[x=G,y=branch_re] {mixed-signature-lorentzian.dat};
\addlegendentry{Branch $1$}
\addplot[orange!85!black,line width=1.2pt] table[x=G,y=leaf11_re] {mixed-signature-lorentzian.dat};
\addlegendentry{Leaf $11$}
\addplot[green!45!black,line width=1.2pt] table[x=G,y=leaf12_re] {mixed-signature-lorentzian.dat};
\addlegendentry{Leaf $12$}
\addplot[blue!75!black,dashed,forget plot,line width=1.1pt] table[x=G,y=branch_im] {mixed-signature-lorentzian.dat};
\addplot[orange!85!black,dashed,forget plot,line width=1.1pt] table[x=G,y=leaf11_im] {mixed-signature-lorentzian.dat};
\addplot[green!45!black,dashed,forget plot,line width=1.1pt] table[x=G,y=leaf12_im] {mixed-signature-lorentzian.dat};
\end{groupplot}
\end{tikzpicture}
\caption{Lorentzian analogue of Figure~\ref{fig:mixed-euclidean-observables},
with $\eta=10^{-8}$. Solid and dashed mean curves show real and
imaginary parts, respectively.}
\label{fig:mixed-lorentzian-observables}
\end{figure}

The opposite-sign leaf edges are shorter on average and have larger
relative uncertainties over the range shown. At each branch, the
two leaf curvature means have opposite signs.
The positive ratio terms in the action for opposite-sign parent and
leaf edges favour unequal magnitudes in the Euclidean theory. The Lorentzian counterparts are shown in
Figure~\ref{fig:mixed-lorentzian-observables}. The means are complex,
and relative uncertainties use $|\Delta x|/|\< x\>|$.
The pronounced difference between the two leaf magnitude classes
in the Euclidean theory is much weaker here. The strong-coupling
limits are the same in both theories.

\section{Quantum gravity on $n$-stars and related trees}\label{sec:othe}

We have focussed on quantum gravity on the 7-vertex binary tree
\(T_3=\widetilde D_6\). The same methodology also allows us to
revisit quantum gravity on the 4-star \(\widetilde D_4\)
of \cite{Ma:alg4} with more careful numerics, and on the 2-star
\(A_3\) and 3-star \(D_4\) for comparison. We then study the adjacent
3-stars \(\widetilde D_5\) and the extended 3-star
\(\widetilde E_6\), with arms \((2,2,2)\).
In all cases, we  use the zero-sum
lifting \eqref{eq:uniform-lift} and choose positive vertex measures $\mu$ for the action that are linear in the metric so as to obtain dimensionless actions of the same type as in 
Proposition~\ref{prop:linear-metric-action}. As for the binary tree, a positive ratio metric can have mixed signature. 

Unless stated otherwise, the numerical plots in this section are done in practice at $L=1$, which suffices by scaling arguments as in Section~\ref{sec:quan}. Relative to this, Euclidean-only plots use \(\delta=10^{-10}\) and start
at \(G=0.001\), while Lorentzian plots use direct measure with 
\(\delta=1/10\), \(\eta=10^{-10}\), and start at \(G=0.1\). The numerical
checks previously reported apply here as well.

\subsection{Quantum gravity on the $n$-star revisited}\label{sec:star}

We use the QRG for the $n$-star, $n=2,3,4$  \cite[Theorem 3.4(a)]{BegMa:gra}, taking, for the 2-star or $A_3$, the 
canonical one (which is leaf-independent) rather than the full one-parameter family in \cite{ArgMa:dis}. We 
label the centre $0$ and the leaves $1,\ldots,n$, and set
$u_i=\omega_{0i}$, $\bar u_i=\omega_{i0}$. The positive metric is
\[
 g=\sum_{i=1}^n a_i(u_i\ten\bar u_i+\sqrt n\,\bar u_i\ten u_i),
 \qquad a_i>0.
\]
The common phase is $s=e^{3\pi\imath{}/4}$ for $n=2$,
$s=e^{5\pi\imath{}/6}$ for $n=3$, and $s=-1$ for $n=4$;
one may also take its complex conjugate. In all cases
$s^2+\sqrt n\,s+1=0$. Writing $P_0=\sum_j u_j\ten\bar u_j$,
the braiding and inner connection are specified by
\begin{align*}
 \sigma(u_i\ten\bar u_i)&=-u_i\ten\bar u_i+(1+s^{-1}/\sqrt n)P_0,\\
 \sigma(\bar u_i\ten u_i)&=s\bar u_i\ten u_i,\qquad
 \sigma(\bar u_i\ten u_j)=-\frac{a_j}{a_i}s^{-1}\bar u_i\ten u_j\quad(i\ne j),
 \qquad\alpha=0.
\end{align*}
For the zero-sum lifting \eqref{eq:uniform-lift}, contraction gives
\[
 \Scal(0)=0,\qquad \Scal(i)=\frac1n\sum_{j=1}^n\frac1{a_j}-\frac1{a_i}.
\]
In particular the two conjugate connections have the same real scalar
curvature. One way to see the cancellation of the phase is that the
coefficients of $u_j\wedge\bar u_j\ten u_i$ in $R(u_i)$ are
$-1/\sqrt n$ for $j=i$ and
$(1+s^{-1}/\sqrt n)s^{-1}a_i/a_j$ for $j\ne i$.
The lifting subtracts their mean, and
$\sqrt n\,s^{-1}+s^{-2}=-1$ gives the stated scalar.

Choose the linear vertex weights
\[
 \mu(i)=\frac{n^2}{2}\sum_{j\ne i}a_j.
\]
The central weight does not contribute. This choice gives the uniform
dimensionless action
\begin{equation}\label{eq:n-star-action}
 S[a]=\sum_{i=1}^n\mu(i)\Scal(i)
 =-\frac n2\sum_{i<j}\frac{(a_i-a_j)^2}{a_i a_j}
 =-\frac n2\left[\left(\sum_i a_i\right)\left(\sum_i a_i^{-1}\right)-n^2\right].
\end{equation}
As for the tree, it is nonpositive and vanishes exactly when all
outward metric coefficients agree. For $n=4$ the inward parameters
used in \cite{Ma:alg4} are $g_i=2a_i$ and the Einstein--Hilbert action used there
is $S[a]-24$ in our terms, where the constant has no effect on normalised expectations.
The 2-star is treated here as a member of this symmetric star family;
we defer the other phase choices relevant to the $A_n$ family in \cite{ArgMa:dis}. 

\begin{figure}[tbp]
\centering
\begin{tikzpicture}
\begin{groupplot}[
 group style={group size=2 by 2,horizontal sep=1.7cm,vertical sep=1.4cm},
 width=.44\linewidth,height=5.2cm,xmin=0,xmax=15,restrict x to domain=0:15,xlabel={$G$},xlabel style={yshift=4pt},
 grid=major,grid style={gray!20},tick label style={font=\small},
 label style={font=\small},title style={font=\small},
 legend style={font=\scriptsize,draw=none},
]
\nextgroupplot[title={Direct measure},ylabel={$\< a\>/L$},ymin=.48,ymax=.84,legend pos=north east]
\addplot[orange!85!black,dashed,line width=1.2pt] table[x=G,y=star2_direct_mean] {revisited-star-gravity.dat};
\addlegendentry{2-star}
\addplot[teal!75!black,dashdotted,line width=1.2pt] table[x=G,y=star3_direct_mean] {revisited-star-gravity.dat};
\addlegendentry{3-star}
\addplot[blue!75!black,line width=1.2pt] table[x=G,y=star4_direct_mean] {revisited-star-gravity.dat};
\addlegendentry{4-star}
\nextgroupplot[title={Direct measure},ylabel={$\Delta a/\< a\>$},ymin=.18,ymax=.6]
\addplot[orange!85!black,dashed,line width=1.2pt] table[x=G,y=star2_direct_relative_sd] {revisited-star-gravity.dat};
\addplot[teal!75!black,dashdotted,line width=1.2pt] table[x=G,y=star3_direct_relative_sd] {revisited-star-gravity.dat};
\addplot[blue!75!black,line width=1.2pt] table[x=G,y=star4_direct_relative_sd] {revisited-star-gravity.dat};
\nextgroupplot[title={Liouville measure},ylabel={$\< a\>/L$},ymin=.02,ymax=.045,scaled y ticks=false]
\addplot[orange!85!black,dashed,line width=1.2pt] table[x=G,y=star2_liouville_mean] {revisited-star-gravity.dat};
\addplot[teal!75!black,dashdotted,line width=1.2pt] table[x=G,y=star3_liouville_mean] {revisited-star-gravity.dat};
\addplot[blue!75!black,line width=1.2pt] table[x=G,y=star4_liouville_mean] {revisited-star-gravity.dat};
\nextgroupplot[title={Liouville measure},ylabel={$\Delta a/\< a\>$},ymin=3.2,ymax=4.2]
\addplot[orange!85!black,dashed,line width=1.2pt] table[x=G,y=star2_liouville_relative_sd] {revisited-star-gravity.dat};
\addplot[teal!75!black,dashdotted,line width=1.2pt] table[x=G,y=star3_liouville_relative_sd] {revisited-star-gravity.dat};
\addplot[blue!75!black,line width=1.2pt] table[x=G,y=star4_liouville_relative_sd] {revisited-star-gravity.dat};
\end{groupplot}
\end{tikzpicture}
\caption{Metric means and relative fluctuations for the $n$-star,
$n=2,3,4$, with direct and Liouville measures, using the cutoff ratio of
\cite[Figure 1(b)]{Ma:alg4}. }
\label{fig:revisited-star-metrics}
\end{figure}
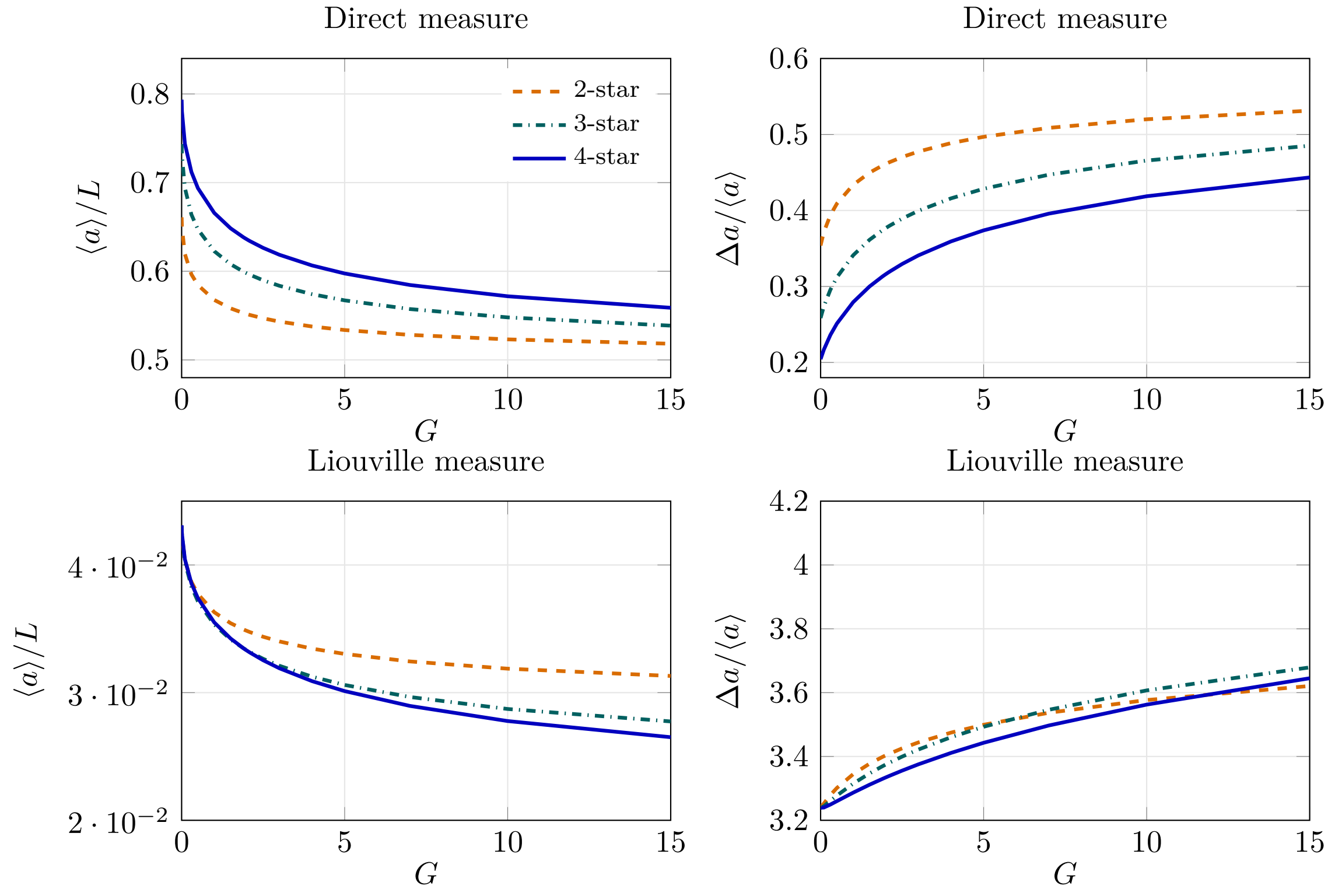

We use $e^{S[a]/G}$ with either direct measure
$\prod_{i=1}^n\mathrm da_i$ or Liouville measure
$\prod_{i=1}^n\mathrm da_i/a_i$, with $\delta<a_i<L$.
All leaves are equivalent in expectation and
\[
 \<\Scal(i)\>=0.
\]
Their curvature fluctuations are nevertheless nonzero. The common
scale can again be integrated exactly.
On a face where the first parameter is largest, write
\[
 (a_1,\ldots,a_n)=L\rho(1,e^{-x_1},\ldots,e^{-x_{n-1}}),\qquad
 m=\max_jx_j,\qquad T=\log(L/\delta).
\]
Then $0\le x_j\le T$, $e^{m-T}<\rho<1$, and $E(x)=-S[a]$
is independent of $\rho$. With
\[
 F_p(m,T)=\int_{e^{m-T}}^1\rho^{p-1}\,\mathrm d\rho
 =\begin{cases}(1-e^{-p(T-m)})/p,&p\ne0,\\T-m,&p=0,\end{cases}
\]
the partition functions are
\begin{align*}
 Z^{\mathrm{dir}}_{\delta,L}
 &=nL^n\int_{[0,T]^{n-1}}e^{-E(x)/G-\sum_jx_j}F_n(m,T)\,\mathrm d^{n-1}x,\\
 Z^{\mathrm{Liouv}}_{\delta,L}
 &=n\int_{[0,T]^{n-1}}e^{-E(x)/G}F_0(m,T)\,\mathrm d^{n-1}x.
\end{align*}
For the $p$th moment of $a_i/L$, insert
$(1+\sum_j e^{-p x_j})/n$ and replace $F_n$ by $F_{n+p}$,
or $F_0$ by $F_p$, respectively. Thus at most three relative
variables remain, with no singular common-scale integration.

\begin{figure}[tbp]
\centering
\begin{tikzpicture}
\begin{groupplot}[
 group style={group size=2 by 1,horizontal sep=1.7cm},
 width=.44\linewidth,height=6cm,xmin=0,xmax=15,restrict x to domain=0:15,xlabel={$G$},xlabel style={yshift=4pt},
 grid=major,grid style={gray!20},tick label style={font=\small},
 label style={font=\small},title style={font=\small},
 legend style={font=\scriptsize,draw=none},scaled y ticks=false,
 yticklabel style={/pgf/number format/fixed,/pgf/number format/precision=2},
]
\nextgroupplot[title={Direct measure},ylabel={$L\Delta\Scal$},ymin=0,legend pos=north west]
\addplot[orange!85!black,dashed,line width=1.2pt] table[x=G,y=star2_direct_curvature_sd] {revisited-star-gravity.dat};
\addlegendentry{2-star}
\addplot[teal!75!black,dashdotted,line width=1.2pt] table[x=G,y=star3_direct_curvature_sd] {revisited-star-gravity.dat};
\addlegendentry{3-star}
\addplot[blue!75!black,line width=1.2pt] table[x=G,y=star4_direct_curvature_sd] {revisited-star-gravity.dat};
\addlegendentry{4-star}
\nextgroupplot[title={Liouville measure},ylabel={$\delta\Delta\Scal$},ymin=0,ytick={0,.02,.04,.06}]
\addplot[orange!85!black,dashed,line width=1.2pt] table[x=G,y=star2_liouville_curvature_sd] {revisited-star-gravity.dat};
\addplot[teal!75!black,dashdotted,line width=1.2pt] table[x=G,y=star3_liouville_curvature_sd] {revisited-star-gravity.dat};
\addplot[blue!75!black,line width=1.2pt] table[x=G,y=star4_liouville_curvature_sd] {revisited-star-gravity.dat};
\end{groupplot}
\end{tikzpicture}
\caption{Leaf scalar-curvature fluctuations for the $n$-star with
direct measure (left) and Liouville measure (right). All scalar-curvature means vanish.} \label{fig:revisited-star-curvature}
\end{figure}
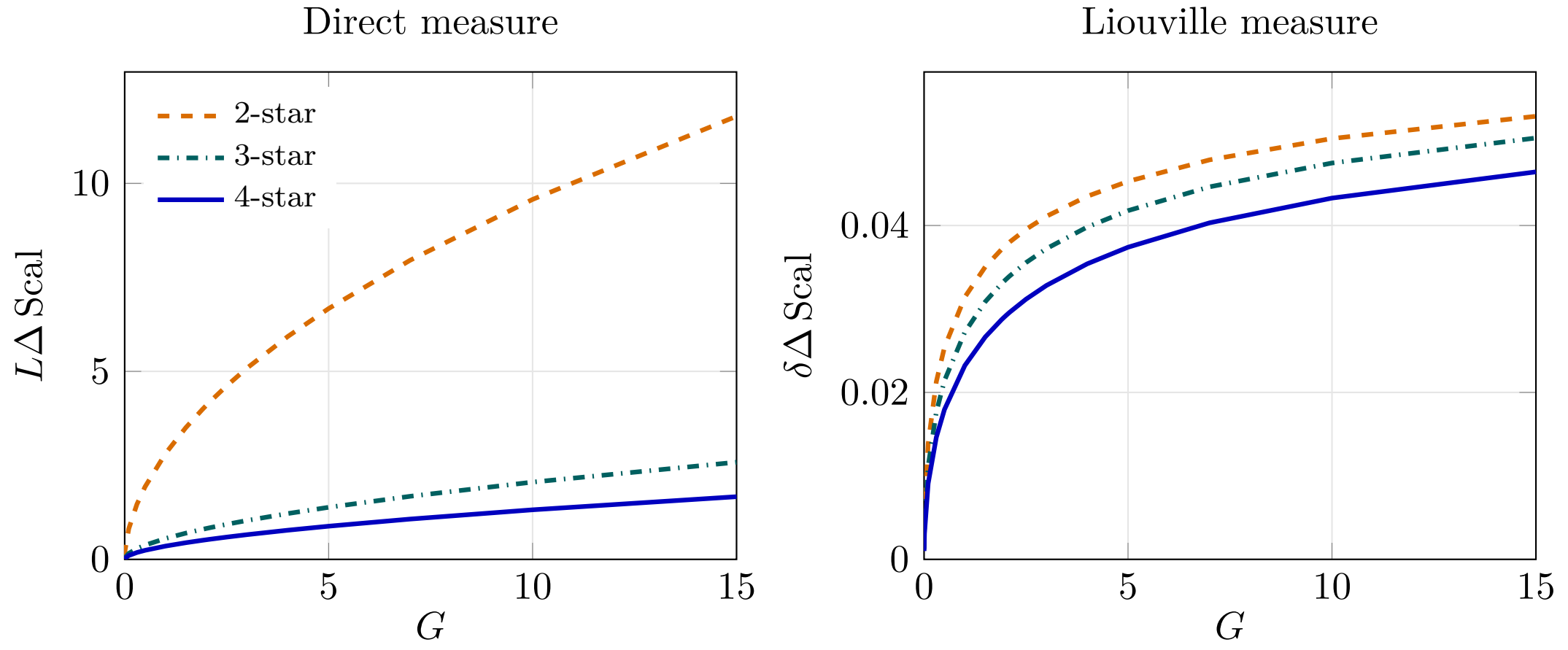

Figure~\ref{fig:revisited-star-metrics} shows smooth behaviour for all
three stars. For direct
measure, $Z$ and the metric moments have finite
limits as $\delta\downarrow0$. In that limit the weak-coupling normalised metric
mean and relative fluctuation tend to $n/(n+1)$ and
$1/\sqrt{n(n+2)}$, respectively. The strong-coupling limits are
$1/2$ and $1/\sqrt3$. For Liouville measure at fixed $T$, both the
weak- and strong-coupling limits are
\[
 \frac{\< a_i\>}{L}\longrightarrow\frac{1-e^{-T}}T,
 \qquad
 \frac{\Delta a_i}{\< a_i\>}
 \longrightarrow\sqrt{\frac{T(1+e^{-T})}{2(1-e^{-T})}-1}.
\]
The intermediate behaviour depends on $n$, as in the tree calculation.
For any fixed $G>0$, the relative energy confines all relative fields
faster than exponentially. Consequently $Z^{\mathrm{Liouv}}/T$ has
a finite positive limit as $T\to\infty$, smoothly depending on $G$.
There is no change of integrability at $G=2$.

The apparent transition in \cite[Figure 1(b)]{Ma:alg4} was carefully revisited and 
traced to a numerical sampling failure in the original
Mathematica integration as follows. At two nearby couplings, $G=1.89690446$
and $1.89993197$, the original calculation of $Z$ gives respectively
$10.89\%$ and $99.999\%$ of its independently checked value.
Below the jump, the smallest sampled $\max_i g_i$ is $0.0627$,
whereas above it the sampling reaches the cutoff $10^{-10}$.
The action is unchanged by common rescaling, while the Liouville
measure $\prod_i\mathrm dg_i/g_i$ gives equal weight to equal
logarithmic intervals of the common scale at fixed ratios.
The omitted small common scales therefore carry substantial weight in $Z$. The sudden
change in the adaptive sampling thus produces a spurious jump in
the normalised mean. Both answers come with slow-convergence warnings but we were unable to attribute
the switch to a particular undocumented internal algorithm. By contrast, integration of the common scale
using our new method as above agrees with an independent numerical integration of all
four logarithmic metric variables and can therefore be considered reliable. 

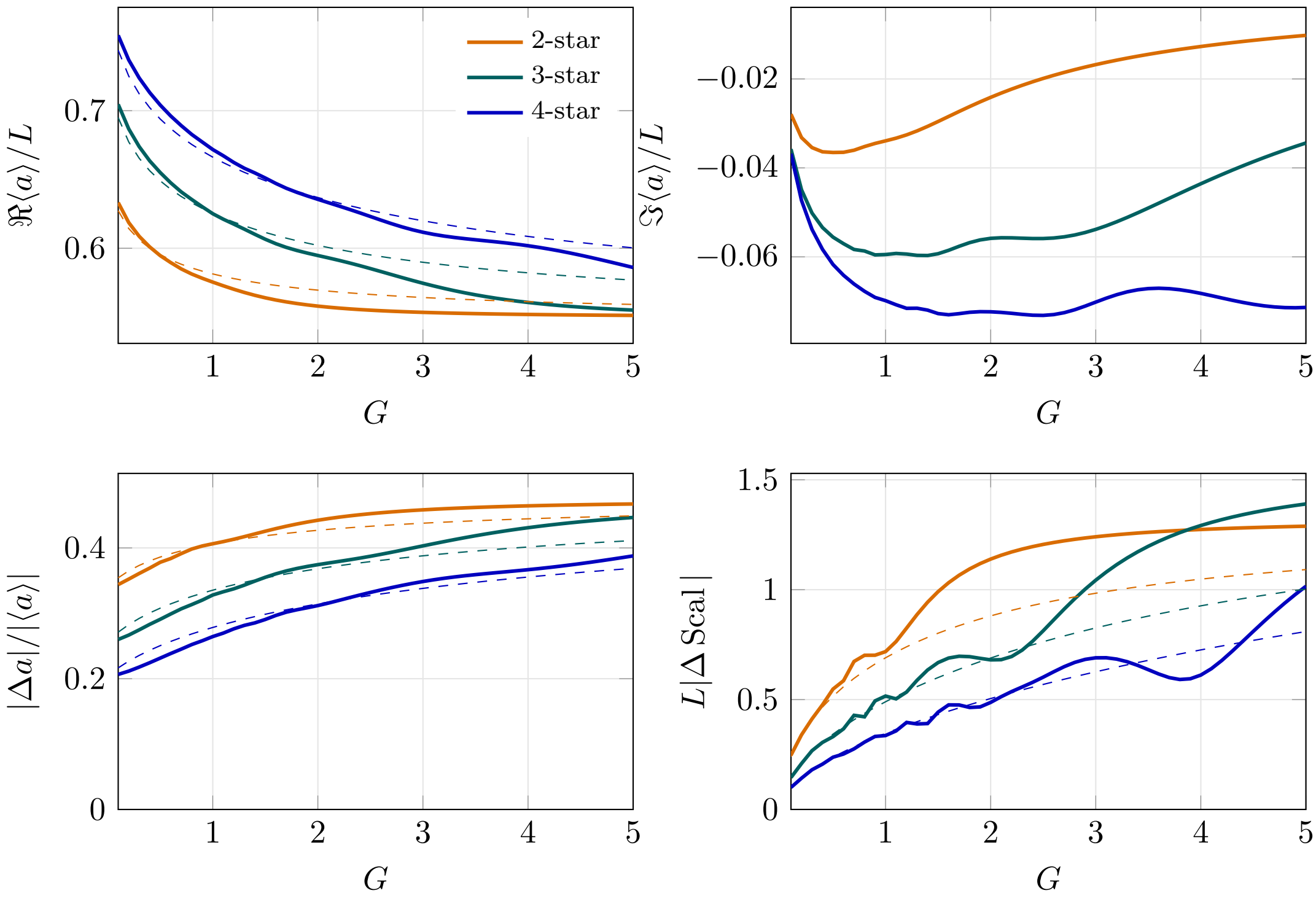
\begin{figure}[tbp]
\centering
\begin{tikzpicture}
\begin{groupplot}[
 group style={group size=2 by 2,horizontal sep=1.7cm,vertical sep=1.4cm},
 width=.44\linewidth,height=5.2cm,xmin=.1,xmax=5,xlabel={$G$},
 grid=major,grid style={gray!20},tick label style={font=\small},
 label style={font=\small},title style={font=\small},
 legend style={font=\scriptsize,draw=none},
]
\nextgroupplot[ylabel={$\Re\< a\>/L$},legend pos=north east]
\addplot[orange!85!black,line width=1.1pt] table[x=G,y=star2_lorentzian_mean_re] {n-star-lorentzian-decade.dat};
\addlegendentry{2-star}
\addplot[teal!75!black,line width=1.1pt] table[x=G,y=star3_lorentzian_mean_re] {n-star-lorentzian-decade.dat};
\addlegendentry{3-star}
\addplot[blue!75!black,line width=1.1pt] table[x=G,y=star4_lorentzian_mean_re] {n-star-lorentzian-decade.dat};
\addlegendentry{4-star}
\addplot[orange!85!black,dashed] table[x=G,y=star2_euclidean_mean] {n-star-lorentzian-decade.dat};
\addplot[teal!75!black,dashed] table[x=G,y=star3_euclidean_mean] {n-star-lorentzian-decade.dat};
\addplot[blue!75!black,dashed] table[x=G,y=star4_euclidean_mean] {n-star-lorentzian-decade.dat};
\nextgroupplot[ylabel={$\Im\< a\>/L$},scaled y ticks=false,
 yticklabel style={/pgf/number format/fixed,/pgf/number format/precision=2}]
\addplot[orange!85!black,line width=1.1pt] table[x=G,y=star2_lorentzian_mean_im] {n-star-lorentzian-decade.dat};
\addplot[teal!75!black,line width=1.1pt] table[x=G,y=star3_lorentzian_mean_im] {n-star-lorentzian-decade.dat};
\addplot[blue!75!black,line width=1.1pt] table[x=G,y=star4_lorentzian_mean_im] {n-star-lorentzian-decade.dat};
\nextgroupplot[ylabel={$|\Delta a|/|\< a\>|$},ymin=0]
\addplot[orange!85!black,line width=1.1pt] table[x=G,y=star2_lorentzian_relative_sd] {n-star-lorentzian-decade.dat};
\addplot[teal!75!black,line width=1.1pt] table[x=G,y=star3_lorentzian_relative_sd] {n-star-lorentzian-decade.dat};
\addplot[blue!75!black,line width=1.1pt] table[x=G,y=star4_lorentzian_relative_sd] {n-star-lorentzian-decade.dat};
\addplot[orange!85!black,dashed] table[x=G,y=star2_euclidean_relative_sd] {n-star-lorentzian-decade.dat};
\addplot[teal!75!black,dashed] table[x=G,y=star3_euclidean_relative_sd] {n-star-lorentzian-decade.dat};
\addplot[blue!75!black,dashed] table[x=G,y=star4_euclidean_relative_sd] {n-star-lorentzian-decade.dat};
\nextgroupplot[ylabel={$L|\Delta\Scal|$},ymin=0]
\addplot[orange!85!black,line width=1.1pt] table[x=G,y=star2_lorentzian_curvature_sd] {n-star-lorentzian-decade.dat};
\addplot[teal!75!black,line width=1.1pt] table[x=G,y=star3_lorentzian_curvature_sd] {n-star-lorentzian-decade.dat};
\addplot[blue!75!black,line width=1.1pt] table[x=G,y=star4_lorentzian_curvature_sd] {n-star-lorentzian-decade.dat};
\addplot[orange!85!black,dashed] table[x=G,y=star2_euclidean_curvature_sd] {n-star-lorentzian-decade.dat};
\addplot[teal!75!black,dashed] table[x=G,y=star3_euclidean_curvature_sd] {n-star-lorentzian-decade.dat};
\addplot[blue!75!black,dashed] table[x=G,y=star4_euclidean_curvature_sd] {n-star-lorentzian-decade.dat};
\end{groupplot}
\end{tikzpicture}
\caption{The $n$-star with Lorentzian action (solid) compared with
Euclidean action (dashed).
The metric means are complex; scalar-curvature means vanish exactly.}
\label{fig:n-star-lorentzian}
\end{figure}

Curvature fluctuations are shown in Figure~\ref{fig:revisited-star-curvature}.
For direct measure the curvature second moment has radial factor
$F_{n-2}$: it diverges logarithmically as $\delta\downarrow0$
for the 2-star, but has a finite limit for the 3- and 4-stars at
fixed $L,G>0$. With Liouville measure, inversion
$a_i\mapsto\delta L/a_i$ preserves the action and measure, giving
\[
 \delta\Delta\Scal(i)=\sqrt{\frac{n-1}{n}
 \left(\frac{\< a_i^2\>}{L^2}
       -\frac{\< a_i a_j\>}{L^2}\right)}\quad(i\ne j).
\]
This also avoids inverse-moment quadrature near the lower cutoff.
At fixed $G>0$, these fluctuations are $O(T^{-1/2})$ as $T\to\infty$.
The root scalar is identically zero in all three models.

Finally, we consider the Lorentzian theory with direct measure,
using the same convention as Section~\ref{sec:lore}, namely
\[
 e^{-\imath{}S[a]/G_\eta},\qquad G_\eta=G-\imath{}\eta,\qquad\eta>0.
\]
The common-scale integration above is unchanged, and the remaining
relative-field weight is $e^{\imath{}E(x)/G_\eta}$.
Permutation symmetry still gives $\<\Scal(i)\>=0$
wherever $Z\ne0$. Metric expectations and curvature second moments
are now complex. We use \eqref{eq:complex-standard-deviation} and
report $|\Delta a|/|\< a\>|$ and $L|\Delta\Scal|$.

For the four-star at small $G$, write $x_j=mt_j$ with fixed
$t_j\in[0,1]$, one of which equals $1$, and deform the $m$ integral
from $[0,T]$ to the contour
\[
 0\longrightarrow\imath{}h\longrightarrow T+\imath{}h
 \longrightarrow T,\qquad 0<h<\pi.
\]
After the common-scale integration, the integrand is entire in $m$,
so Cauchy's theorem preserves the finite-cutoff integral. On the
horizontal segment $\Im E\ge0$, introducing exponential damping.
With this technique, we can set $\delta=1/100$  with results which have weaker interference ripples compared to Figure~\ref{fig:n-star-lorentzian}. We retained $\delta=1/10$ in order to illustrate these effects.

\subsection{Quantum gravity on adjacent 3-stars \(\widetilde D_5\)}\label{sec:join}

Here, we consider the smallest double fork \(\widetilde D_5\) with distinct branching vertices. Label the trivalent vertices \(1,2\), and their leaves
\(11,12,21,22\). Set \(u_i=\omega_{i,\bar i}\),
\(u_{i\alpha}=\omega_{i,i\alpha}\), with bars denoting reverse
arrows. Here \(i,\alpha\in\{1,2\}\), and \(\bar i,\bar\alpha\)
denote the other indices, so \(u_2=\bar u_1\). See Figure~\ref{fig:two-fork}(a).

\begin{figure}[tbp]
\centering
\begin{minipage}[c]{.49\linewidth}\centering
\begin{tikzpicture}[>=Stealth,vertex/.style={circle,draw,inner sep=2pt},every node/.style={font=\small}]
\node[vertex](b1) at (-1.3,0){1};\node[vertex](b2) at (1.3,0){2};
\node[vertex](l11) at (-2.2,1.5588){11};\node[vertex](l12) at (-2.2,-1.5588){12};
\node[vertex](l21) at (2.2,1.5588){21};\node[vertex](l22) at (2.2,-1.5588){22};
\draw[->,bend left=8](b1) to node[midway,auto=left,inner sep=2pt]{$u_1$}(b2);
\draw[->,bend left=8](b2) to node[midway,auto=left,inner sep=2pt]{$u_2$}(b1);
\foreach \a/\b/\lab in {b1/l11/u_{11},b1/l12/u_{12},b2/l21/u_{21},b2/l22/u_{22}}{
\draw[->,bend left=8](\a) to node[midway,auto=left,inner sep=2pt]{$\lab$}(\b);
\draw[->,bend left=8](\b) to node[midway,auto=left,inner sep=2pt]{$\overline{\lab}$}(\a);}
\end{tikzpicture}
\par\small (a) Adjacent 3-stars $\widetilde D_5$.
\end{minipage}\hfill
\begin{minipage}[c]{.49\linewidth}\centering
\begin{tikzpicture}[>=Stealth,vertex/.style={circle,draw,inner sep=2pt},every node/.style={font=\small}]
\node[vertex](r) at (0,0){0};
\node[vertex](b1) at (-1.2124,-.7){1};\node[vertex](b2) at (0,1.4){2};\node[vertex](b3) at (1.2124,-.7){3};
\node[vertex](l1) at (-2.4249,-1.4){$1'$};\node[vertex](l2) at (0,2.8){$2'$};\node[vertex](l3) at (2.4249,-1.4){$3'$};
\foreach \a/\b/\lab in {r/b1/u_1,r/b2/u_2,r/b3/u_3,b1/l1/u'_1,b2/l2/u'_2,b3/l3/u'_3}{
\draw[->,bend left=8](\a) to node[midway,auto=left,inner sep=2pt]{$\lab$}(\b);
\draw[->,bend left=8](\b) to node[midway,auto=left,inner sep=2pt]{$\overline{\lab}$}(\a);}
\end{tikzpicture}
\par\small (b) Extended 3-star $\widetilde E_6$.
\end{minipage}
\caption{The two additional graphs with positive ratio metrics.}
\label{fig:two-fork}\label{fig:three-arm}
\end{figure}

\subsubsection{Calculus and QRG}
All non-returning two-paths vanish in \(\Om^2_{\min}\), and the
return relations are
\[
 u_i\wedge u_{\bar i}+\sum_{\alpha=1}^2
 u_{i\alpha}\wedge\bar u_{i\alpha}=0,\qquad
 \bar u_{i\alpha}\wedge u_{i\alpha}=0.
\]
Thus the four \(u_{i\alpha}\wedge\bar u_{i\alpha}\) form a basis
of \(\Om^2_{\min}\). For every arrow,
\(\extd\omega_{xy}=\omega_{xy}\wedge\omega_{yx}
+\omega_{yx}\wedge\omega_{xy}\), and \(\Om^m_{\min}=0\) for
\(m\ge3\), as for \(T_3\). The positive metric is
\[
 g=a(u_1\ten u_2+u_2\ten u_1)
   +\sum_{i,\alpha}a_{i\alpha}
       (u_{i\alpha}\ten\bar u_{i\alpha}
                        +2\bar u_{i\alpha}\ten u_{i\alpha}),
 \qquad a,a_{i\alpha}>0.
\]
This is Proposition~\ref{prop:split-four-star} with \(m=1\).
Its unique star-preserving QLC has \(\alpha=0\),
\(Q_{12}=0\) and \(Q_{i,i\alpha}=Q_{i\alpha,i}=-a_{i\alpha}\).
The braiding follows from
\eqref{eq:return-sigma}--\eqref{eq:straight-sigma}, with
\(r_{i,\bar i}=1\), \(r_{i,i\alpha}=1/2\) and
\(r_{i\alpha,i}=0\). Explicitly, its connection is
\begin{align*}
 \nabla u_i={}&u_{\bar i}\ten u_i
 +\sum_\alpha\bar u_{i\alpha}\ten u_i
 -\sum_\alpha u_{i\alpha}\ten\bar u_{i\alpha}
 -\sum_\alpha\frac{a_{\bar i\alpha}}a u_i\ten u_{\bar i\alpha},\\
 \nabla u_{i\alpha}={}&u_{\bar i}\ten u_{i\alpha}
 +\sum_\gamma\bar u_{i\gamma}\ten u_{i\alpha}
 +\frac12\left(u_{i\alpha}\ten\bar u_{i\alpha}
 -u_{i\bar\alpha}\ten\bar u_{i\bar\alpha}-u_i\ten u_{\bar i}\right),\\
 \nabla\bar u_{i\alpha}={}&u_{i\alpha}\ten\bar u_{i\alpha}
 +\bar u_{i\alpha}\ten u_{i\alpha}
 -\frac a{a_{i\alpha}}\bar u_{i\alpha}\ten u_i
 -\frac{a_{i\bar\alpha}}{a_{i\alpha}}
          \bar u_{i\alpha}\ten u_{i\bar\alpha}.
\end{align*}
Here every summed index ranges over \(\{1,2\}\).

\subsubsection{Einstein--Hilbert action and quantum gravity}
Contraction of the curvature gives
\[
 \Scal(i)=\frac13\sum_{\alpha=1}^2
                  \left(\frac1{a_{\bar i\alpha}}-\frac{a_{i\alpha}}{a^2}\right),\quad \Scal(i\alpha)=\frac13
       \left(\frac1a-\frac2{a_{i\alpha}}+\frac1{a_{i\bar\alpha}}\right).
\]
The vertex measure
\[
 \mu(i)=3a,\qquad \mu(i\alpha)=6a
\]
therefore gives the Einstein--Hilbert action
\begin{equation}\label{eq:two-fork-action}
 S[a]=-\sum_{i,\alpha=1}^2\frac{(a-a_{i\alpha})^2}{a a_{i\alpha}}.
\end{equation}
It is nonpositive and vanishes precisely when all five parameters
agree. The connection lifting \(\iota_\sigma\) gives the same
action with \(\mu(x)=4a\) at every vertex.

\begin{figure}[tbp]
\centering
\begin{tikzpicture}
\begin{groupplot}[
 group style={group size=2 by 2,horizontal sep=2.8cm,vertical sep=1.5cm},
 width=.39\linewidth,height=5.2cm,xmin=0,xmax=15,xlabel={$G$},xlabel style={yshift=4pt},
 grid=major,grid style={gray!20},tick label style={font=\small},
 label style={font=\small},title style={font=\small},scaled y ticks=false,
 legend style={font=\scriptsize,draw=none,fill=white,fill opacity=.95,text opacity=1,inner sep=2pt},legend columns=2,
]
\nextgroupplot[title={Direct measure},ylabel={$\<a\>/L$},legend style={font=\scriptsize,draw=none,fill=white,fill opacity=.95,text opacity=1,inner sep=2pt,at={(1.18,.98)},anchor=north east}]
\addplot[orange!85!black,line width=1.1pt] table[x=G,y=forks_direct_inner_mean_re] {additional-tree-euclidean.dat};
\addlegendentry{Adjacent: $a$}
\addplot[teal!75!black,line width=1.1pt] table[x=G,y=forks_direct_outer_mean_re] {additional-tree-euclidean.dat};
\addlegendentry{Adjacent: $a_{i\alpha}$}
\addplot[blue!75!black,line width=1.1pt] table[x=G,y=arms_direct_inner_mean_re] {additional-tree-euclidean.dat};
\addlegendentry{Extended: $a_i$}
\addplot[violet!85!black,line width=1.1pt] table[x=G,y=arms_direct_outer_mean_re] {additional-tree-euclidean.dat};
\addlegendentry{Extended: $a'_i$}
\nextgroupplot[title={Direct measure},ylabel={$\Delta a/\<a\>$}]
\addplot[orange!85!black,line width=1.1pt] table[x=G,y=forks_direct_inner_relative_sd_re] {additional-tree-euclidean.dat};
\addplot[teal!75!black,line width=1.1pt] table[x=G,y=forks_direct_outer_relative_sd_re] {additional-tree-euclidean.dat};
\addplot[blue!75!black,line width=1.1pt] table[x=G,y=arms_direct_inner_relative_sd_re] {additional-tree-euclidean.dat};
\addplot[violet!85!black,line width=1.1pt] table[x=G,y=arms_direct_outer_relative_sd_re] {additional-tree-euclidean.dat};
\nextgroupplot[title={Liouville measure},ylabel={$\<a\>/L$},scaled y ticks=false,
 ytick={.01,.02,.03,.04},
 yticklabel style={/pgf/number format/fixed,/pgf/number format/precision=2}]
\addplot[orange!85!black,line width=1.1pt] table[x=G,y=forks_liouville_inner_mean_re] {additional-tree-euclidean.dat};
\addplot[teal!75!black,line width=1.1pt] table[x=G,y=forks_liouville_outer_mean_re] {additional-tree-euclidean.dat};
\addplot[blue!75!black,line width=1.1pt] table[x=G,y=arms_liouville_inner_mean_re] {additional-tree-euclidean.dat};
\addplot[violet!85!black,line width=1.1pt] table[x=G,y=arms_liouville_outer_mean_re] {additional-tree-euclidean.dat};
\nextgroupplot[title={Liouville measure},ylabel={$\Delta a/\<a\>$}]
\addplot[orange!85!black,line width=1.1pt] table[x=G,y=forks_liouville_inner_relative_sd_re] {additional-tree-euclidean.dat};
\addplot[teal!75!black,line width=1.1pt] table[x=G,y=forks_liouville_outer_relative_sd_re] {additional-tree-euclidean.dat};
\addplot[blue!75!black,line width=1.1pt] table[x=G,y=arms_liouville_inner_relative_sd_re] {additional-tree-euclidean.dat};
\addplot[violet!85!black,line width=1.1pt] table[x=G,y=arms_liouville_outer_relative_sd_re] {additional-tree-euclidean.dat};
\end{groupplot}
\end{tikzpicture}
\caption{Metric means and relative fluctuations for \(\widetilde D_5\) and the extended-arm graph of Section~\ref{sec:exte}.}
\label{fig:extra-metrics}
\end{figure}
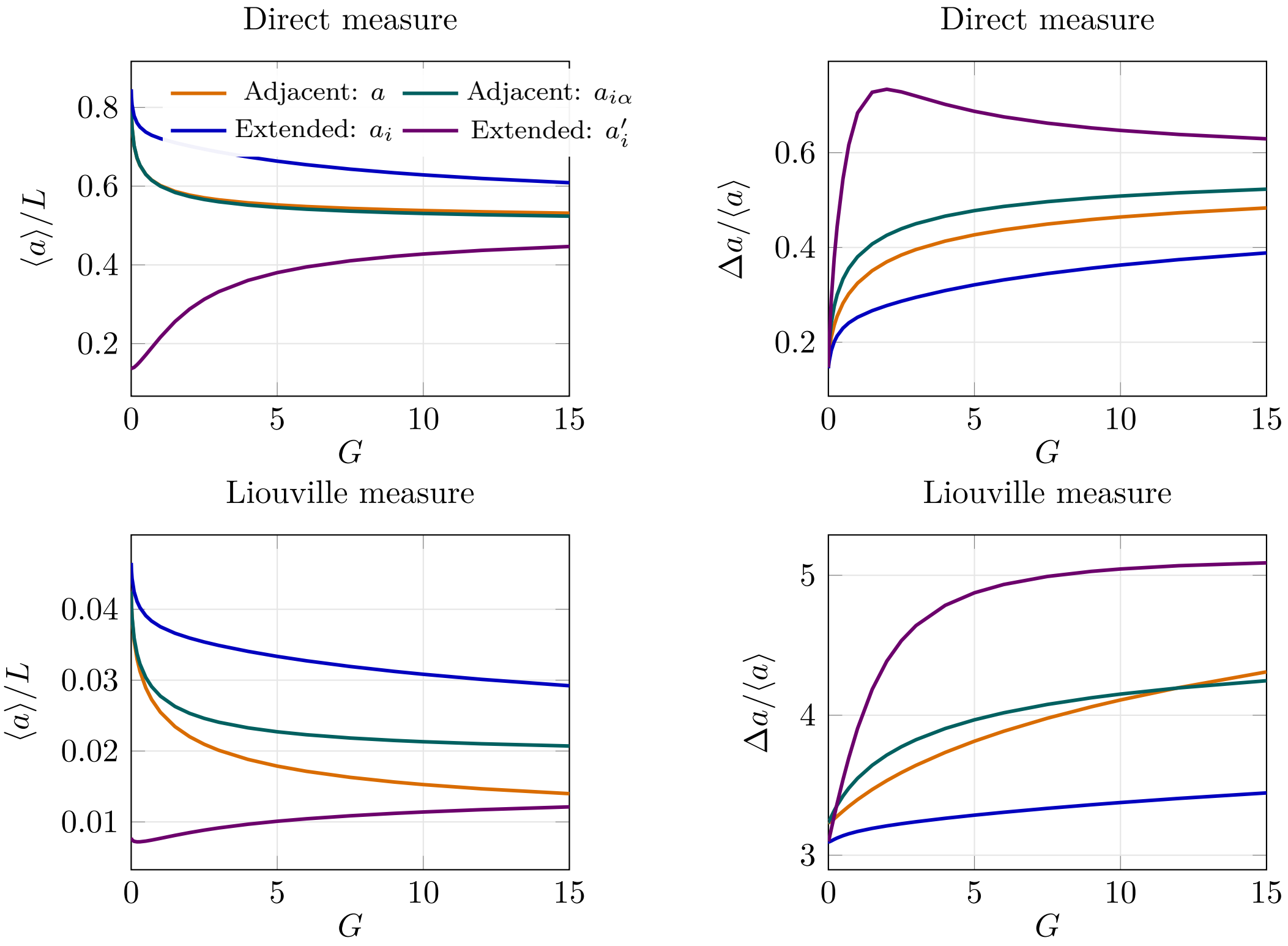

For quantum gravity, cut off each of the five free parameters to
\((\delta,L)\), and use \(e^{S[a]/G}\), \(G>0\).
The four leaf integrations are independent at fixed \(a\).
For direct measure (\(\nu=1\)) or Liouville measure (\(\nu=0\)),
put
\[
 J_\nu(a;c)=\int_\delta^L b^{\nu-1}
       \exp\left[-\frac{(a-cb)^2}{Gcab}\right]\,\mathrm db.
\]
The partition function reduces to
\[
 Z_\nu=\int_\delta^L a^{\nu-1}J_\nu(a;1)^4\,\mathrm da.
\]
Observable moments are obtained by inserting the observable before
performing the leaf integrations. Scale invariance makes
dimensionless expectations functions of \(G\) and \(\delta/L\).
Direct measure allows \(\delta=0\) for the partition function and
bounded metric observables, since \(0<e^{S/G}\le1\).
For Liouville measure the undamped common scale gives a logarithmic
divergence as \(\delta\to0\).

\begin{figure}[tbp]
\centering
\begin{tikzpicture}
\begin{groupplot}[
 group style={group name=extra18,group size=2 by 2,horizontal sep=2.6cm,vertical sep=2.4cm},
 width=.40\linewidth,height=5.2cm,xmin=0,xmax=15,xlabel={$G$},xlabel style={yshift=4pt},
 grid=major,grid style={gray!20},tick label style={font=\small},
 label style={font=\small},title style={font=\small},scaled y ticks=false,
 yticklabel style={/pgf/number format/fixed,/pgf/number format/precision=3},
 legend style={font=\scriptsize,draw=none},legend columns=-1,
]
\nextgroupplot[title={Direct measure},ylabel={$L\<\Scal\>$},legend to name=leg:extra-curv-forks]
\addplot[orange!85!black,line width=1.1pt] table[x=G,y=forks_direct_curv0_mean_re] {additional-tree-euclidean.dat};
\addlegendentry{$i$}
\addplot[teal!75!black,line width=1.1pt] table[x=G,y=forks_direct_curv1_mean_re] {additional-tree-euclidean.dat};
\addlegendentry{$i\alpha$}
\nextgroupplot[title={Liouville measure},ylabel={$\delta\<\Scal\>$}]
\addplot[orange!85!black,line width=1.1pt] table[x=G,y=forks_liouville_curv0_mean_re] {additional-tree-euclidean.dat};
\addplot[teal!75!black,line width=1.1pt] table[x=G,y=forks_liouville_curv1_mean_re] {additional-tree-euclidean.dat};
\nextgroupplot[title={Direct measure},ylabel={$L\Delta\Scal$}]
\addplot[orange!85!black,line width=1.1pt] table[x=G,y=forks_direct_curv0_sd_re] {additional-tree-euclidean.dat};
\addplot[teal!75!black,line width=1.1pt] table[x=G,y=forks_direct_curv1_sd_re] {additional-tree-euclidean.dat};
\nextgroupplot[title={Liouville measure},ylabel={$\delta\Delta\Scal$}]
\addplot[orange!85!black,line width=1.1pt] table[x=G,y=forks_liouville_curv0_sd_re] {additional-tree-euclidean.dat};
\addplot[teal!75!black,line width=1.1pt] table[x=G,y=forks_liouville_curv1_sd_re] {additional-tree-euclidean.dat};
\end{groupplot}
\node[anchor=center,font=\scriptsize,inner sep=0pt] at ($(extra18 c1r1.south west)!0.5!(extra18 c2r2.north east)$) {
\pgfplotslegendfromname{leg:extra-curv-forks}
};
\end{tikzpicture}
\caption{Adjacent 3-stars: local scalar-curvature means and fluctuations.}
\label{fig:extra-curvature}
\end{figure}
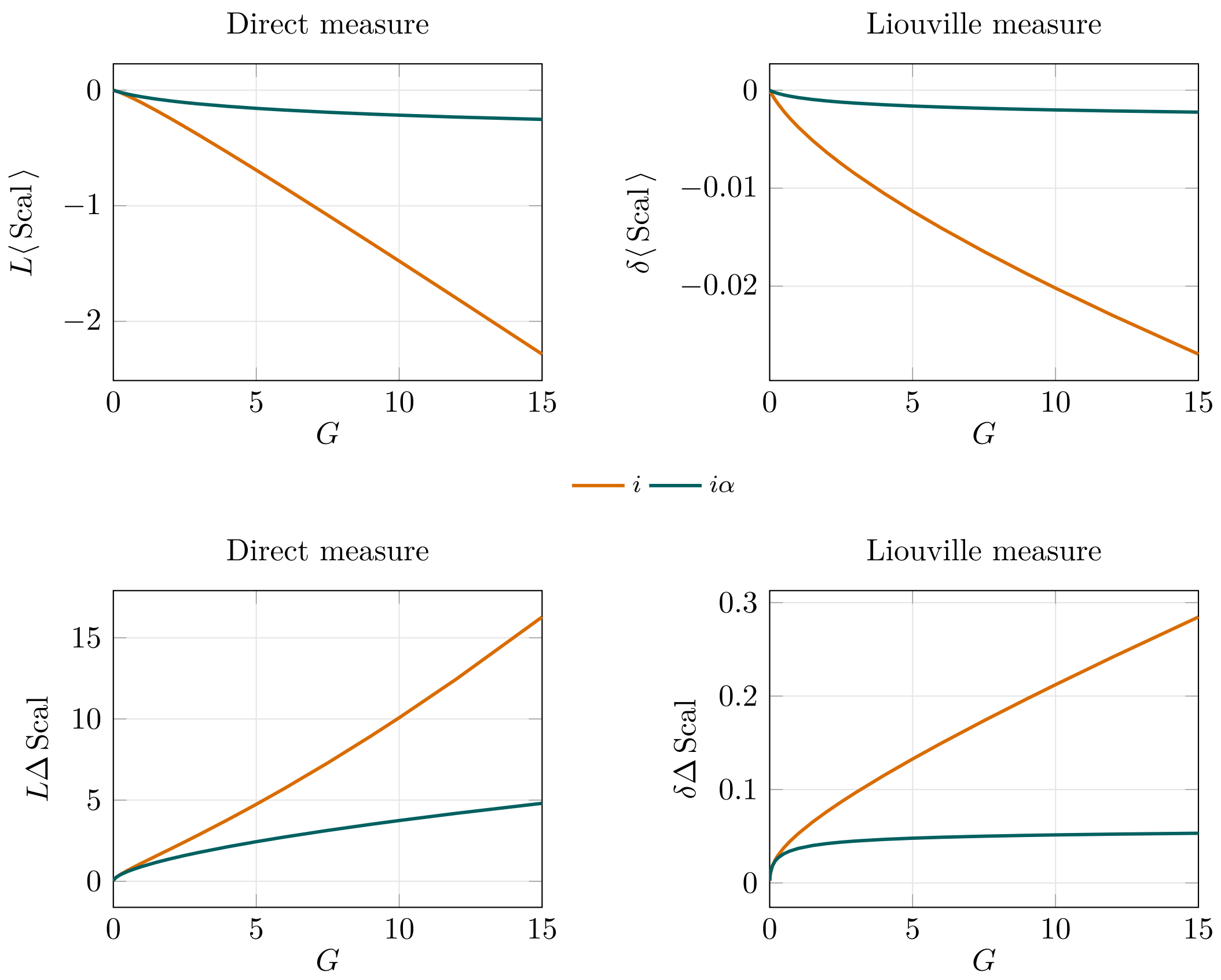

\begin{figure}[tbp]
\centering
\begin{tikzpicture}
\begin{groupplot}[
 group style={group name=extra20,group size=2 by 2,horizontal sep=2.6cm,vertical sep=2.4cm},
 width=.40\linewidth,height=5.2cm,xmin=0.1,xmax=5,xlabel={$G$},xlabel style={yshift=4pt},
 grid=major,grid style={gray!20},tick label style={font=\small},
 label style={font=\small},title style={font=\small},scaled y ticks=false,
 yticklabel style={/pgf/number format/fixed,/pgf/number format/precision=3},
 legend style={font=\scriptsize,draw=none},legend columns=-1,
]
\nextgroupplot[ylabel={$\Re\<a\>/L$},legend to name=leg:extra-lor-met-forks]
\addplot[orange!85!black,line width=1.1pt] table[x=G,y=O_forks_direct_inner_mean_re] {j6-lorentzian-decade.dat};
\addlegendentry{$a$}
\addplot[orange!85!black,dashed,line width=1.1pt,forget plot] table[x=G,y=E_forks_direct_inner_mean_re] {j6-lorentzian-decade.dat};
\addplot[teal!75!black,line width=1.1pt] table[x=G,y=O_forks_direct_outer_mean_re] {j6-lorentzian-decade.dat};
\addlegendentry{$a_{i\alpha}$}
\addplot[teal!75!black,dashed,line width=1.1pt,forget plot] table[x=G,y=E_forks_direct_outer_mean_re] {j6-lorentzian-decade.dat};
\nextgroupplot[ylabel={$\Im\<a\>/L$}]
\addplot[orange!85!black,line width=1.1pt] table[x=G,y=O_forks_direct_inner_mean_im] {j6-lorentzian-decade.dat};
\addplot[teal!75!black,line width=1.1pt] table[x=G,y=O_forks_direct_outer_mean_im] {j6-lorentzian-decade.dat};
\nextgroupplot[ylabel={$|\Delta a|/|\<a\>|$}]
\addplot[orange!85!black,line width=1.1pt] table[x=G,y=O_forks_direct_inner_relative_sd_re] {j6-lorentzian-decade.dat};
\addplot[orange!85!black,dashed,line width=1.1pt,forget plot] table[x=G,y=E_forks_direct_inner_relative_sd_re] {j6-lorentzian-decade.dat};
\addplot[teal!75!black,line width=1.1pt] table[x=G,y=O_forks_direct_outer_relative_sd_re] {j6-lorentzian-decade.dat};
\addplot[teal!75!black,dashed,line width=1.1pt,forget plot] table[x=G,y=E_forks_direct_outer_relative_sd_re] {j6-lorentzian-decade.dat};
\nextgroupplot[ylabel={$L|\Delta\Scal|$},legend to name=leg:extra-lor-curv-forks]
\addplot[orange!85!black,line width=1.1pt] table[x=G,y=O_forks_direct_curv0_sd_re] {j6-lorentzian-decade.dat};
\addlegendentry{$i$}
\addplot[orange!85!black,dashed,line width=1.1pt,forget plot] table[x=G,y=E_forks_direct_curv0_sd_re] {j6-lorentzian-decade.dat};
\addplot[teal!75!black,line width=1.1pt] table[x=G,y=O_forks_direct_curv1_sd_re] {j6-lorentzian-decade.dat};
\addlegendentry{$i\alpha$}
\addplot[teal!75!black,dashed,line width=1.1pt,forget plot] table[x=G,y=E_forks_direct_curv1_sd_re] {j6-lorentzian-decade.dat};
\end{groupplot}
\node[anchor=center,font=\scriptsize,inner sep=0pt] at ($(extra20 c1r1.south west)!0.5!(extra20 c2r2.north east)$) {
\begin{tabular}{@{}l@{\quad}l@{}}
Metrics: & \pgfplotslegendfromname{leg:extra-lor-met-forks}\\[2pt]
Vertices: & \pgfplotslegendfromname{leg:extra-lor-curv-forks}
\end{tabular}
};
\end{tikzpicture}
\caption{Adjacent 3-stars: Lorentzian (solid) and Euclidean (dashed) expectations.}
\label{fig:extra-lorentzian}
\end{figure}
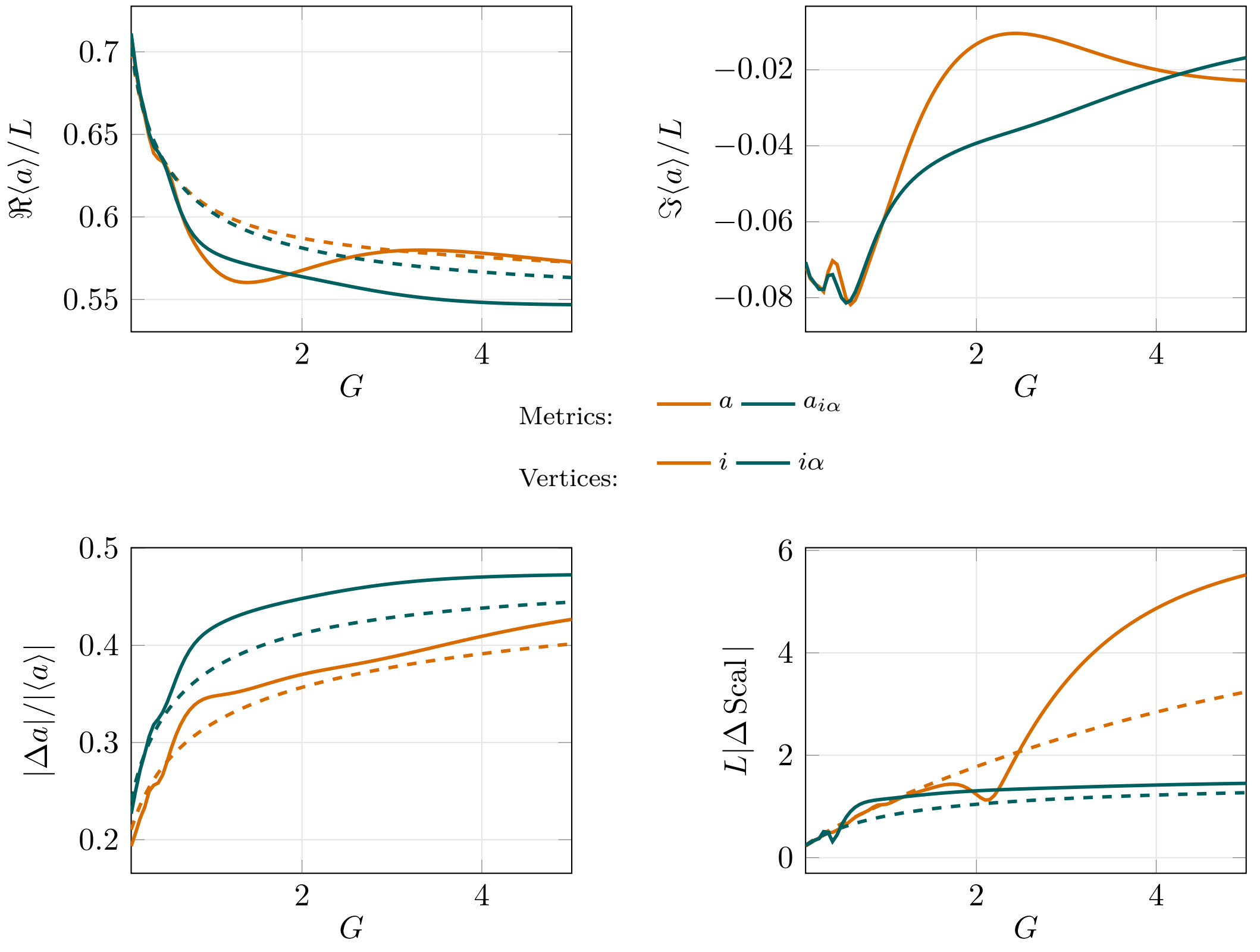

For direct measure with \(\delta=0\), the four transverse Gaussian
integrations at weak coupling leave the common scale with density
proportional to \(a^4\,\mathrm da\), \(0<a<L\). Hence
\[
 \frac{\< a\>}{L},\frac{\< a_{i\alpha}\>}{L}
 \longrightarrow\frac56,\qquad
 \frac{\Delta a}{\< a\>},
 \frac{\Delta a_{i\alpha}}{\< a_{i\alpha}\>}
 \longrightarrow\frac1{\sqrt{35}}\quad(G\downarrow0).
\]
At strong coupling every free metric parameter is uniform on
\((0,L)\), giving mean \(L/2\) and relative uncertainty
\(1/\sqrt3\). Results for the metric are given in Figure~\ref{fig:extra-metrics} and for the scalar curvature
in Figure~\ref{fig:extra-curvature}. The Lorentzian theory uses
\(e^{-\imath{}S[a]/(G-\imath{}\eta)}\), as in
Section~\ref{sec:lore} with results in Figure~\ref{fig:extra-lorentzian}. As for the $n$-star we were also able to do this for $\delta=1/100$ but have retained the larger value to better illustrate  interference effects.

\subsection{Quantum gravity on the extended 3-star  \(\widetilde E_6\)}\label{sec:exte}
For the seven-vertex tree with arms \((2,2,2)\), label the centre \(0\), the intermediate vertices \(1,2,3\),
and the leaves \(1',2',3'\). Write
\(u_i=\omega_{0i}\), \(u'_i=\omega_{i,i'}\), with bars denoting
reverse arrows and \(i\in\{1,2,3\}\). See Figure~\ref{fig:three-arm}(b).

\subsubsection{Calculus and QRG}

The return relations in \(\Om^2_{\min}\) are
\[
 \sum_{i=1}^3u_i\wedge\bar u_i=0,\qquad
 \bar u_i\wedge u_i+u'_i\wedge\overline{u'_i}=0,\qquad
 \overline{u'_i}\wedge u'_i=0.
\]
All non-returning two-paths vanish. Thus \(\Om^2_{\min}\) has
basis \(u_1\wedge\bar u_1,u_2\wedge\bar u_2,
 u'_i\wedge\overline{u'_i}\) for \(i=1,2,3\), and dimension five.
The exterior differential has the same arrow formula as above,
and \(\Om^m_{\min}=0\) for \(m\ge3\). Take the metric
\[
 g=\sum_{i=1}^3\left[
 a_i\left(u_i\ten\bar u_i+\tfrac32\bar u_i\ten u_i\right)
 +a'_i\left(u'_i\ten\overline{u'_i}
                         +2\overline{u'_i}\ten u'_i\right)\right],
 \qquad a_i,a'_i>0.
\]
Theorem~\ref{thm:one-branch} applies with \(\lambda=2\),
\(z=q=1\), and arm lengths \((2,2,2)\). Its unique QLC has
\(Q_{0i}=Q_{i0}=-a_i/2\) and \(Q_{ii'}=Q_{i'i}=a'_i\).
Here \(\alpha=0\), and the return parameters are
\(r_{0i}=2/3\), \(r_{i0}=1/2\), \(r_{i,i'}=3/2\),
\(r_{i',i}=2\). Equations
\eqref{eq:return-sigma}--\eqref{eq:straight-sigma} specify the
braiding. The connection on all twelve arrows is
\begin{align*}
 \nabla u_i={}&\sum_{j=1}^3\bar u_j\ten u_i+u_i\ten\bar u_i
 -\frac23\sum_{j=1}^3u_j\ten\bar u_j
 -\frac{2a'_i}{a_i}u_i\ten u'_i,\\
 \nabla\bar u_i={}&u_i\ten\bar u_i
 +\overline{u'_i}\ten\bar u_i
 +\frac12\left(\bar u_i\ten u_i-u'_i\ten\overline{u'_i}\right)
 -\sum_{j\ne i}\frac{a_j}{a_i}\bar u_i\ten u_j,\\
 \nabla u'_i={}&u_i\ten u'_i+\overline{u'_i}\ten u'_i
 -\frac12u'_i\ten\overline{u'_i}-\frac32\bar u_i\ten u_i,\\
 \nabla\overline{u'_i}={}&u'_i\ten\overline{u'_i}
 -\overline{u'_i}\ten u'_i
 -\frac{a_i}{2a'_i}\overline{u'_i}\ten\bar u_i.
\end{align*}

\subsubsection{Einstein--Hilbert action and quantum gravity}

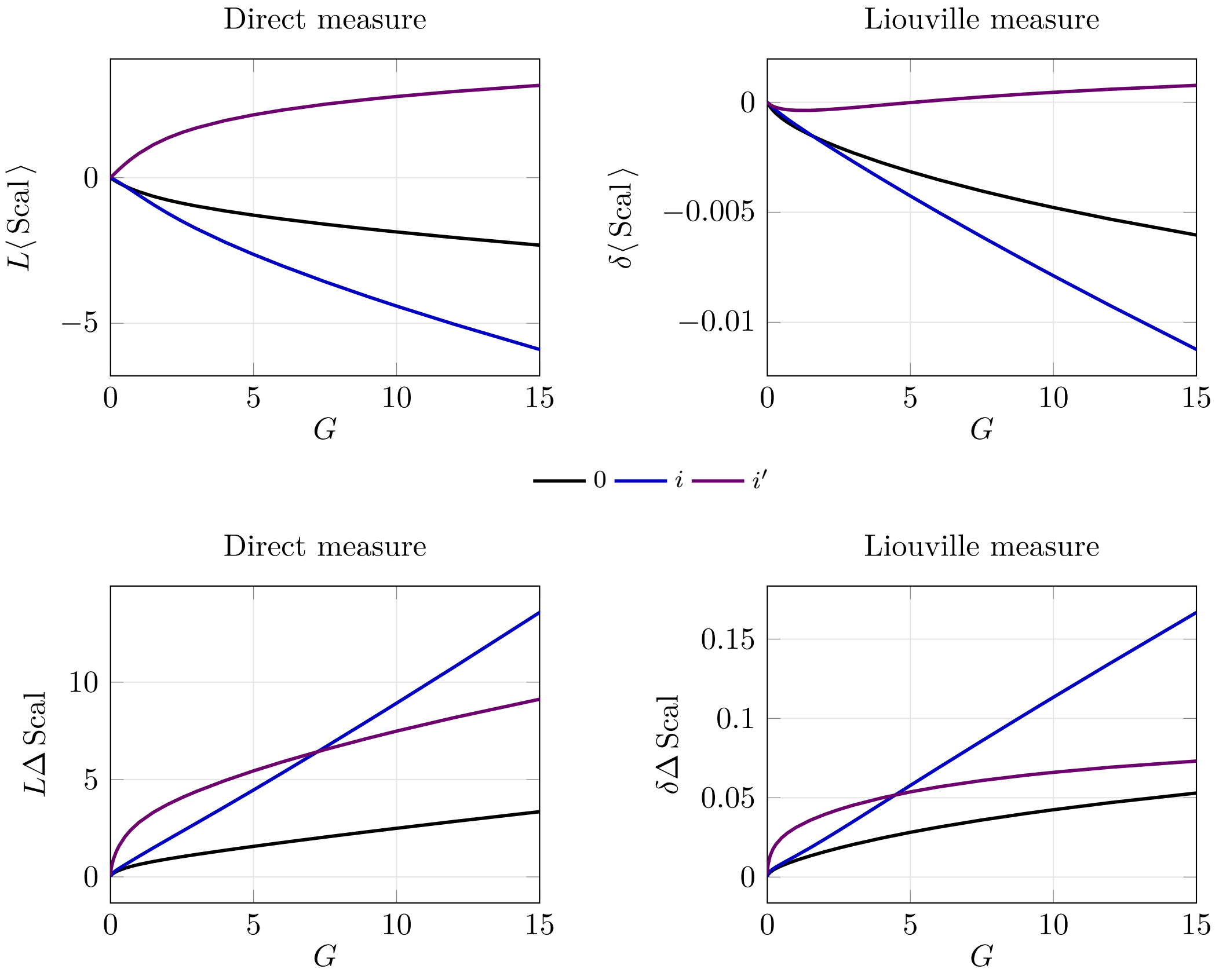
\begin{figure}[bp]
\centering
\begin{tikzpicture}
\begin{groupplot}[
 group style={group name=extra19,group size=2 by 2,horizontal sep=2.6cm,vertical sep=2.4cm},
 width=.40\linewidth,height=5.2cm,xmin=0,xmax=15,xlabel={$G$},xlabel style={yshift=4pt},
 grid=major,grid style={gray!20},tick label style={font=\small},
 label style={font=\small},title style={font=\small},scaled y ticks=false,
 yticklabel style={/pgf/number format/fixed,/pgf/number format/precision=3},
 legend style={font=\scriptsize,draw=none},legend columns=-1,
]
\nextgroupplot[title={Direct measure},ylabel={$L\<\Scal\>$},legend to name=leg:extra-curv-arms]
\addplot[black,line width=1.1pt] table[x=G,y=arms_direct_curv0_mean_re] {additional-tree-euclidean.dat};
\addlegendentry{$0$}
\addplot[blue!75!black,line width=1.1pt] table[x=G,y=arms_direct_curv1_mean_re] {additional-tree-euclidean.dat};
\addlegendentry{$i$}
\addplot[violet!85!black,line width=1.1pt] table[x=G,y=arms_direct_curv2_mean_re] {additional-tree-euclidean.dat};
\addlegendentry{$i'$}
\nextgroupplot[title={Liouville measure},ylabel={$\delta\<\Scal\>$}]
\addplot[black,line width=1.1pt] table[x=G,y=arms_liouville_curv0_mean_re] {additional-tree-euclidean.dat};
\addplot[blue!75!black,line width=1.1pt] table[x=G,y=arms_liouville_curv1_mean_re] {additional-tree-euclidean.dat};
\addplot[violet!85!black,line width=1.1pt] table[x=G,y=arms_liouville_curv2_mean_re] {additional-tree-euclidean.dat};
\nextgroupplot[title={Direct measure},ylabel={$L\Delta\Scal$}]
\addplot[black,line width=1.1pt] table[x=G,y=arms_direct_curv0_sd_re] {additional-tree-euclidean.dat};
\addplot[blue!75!black,line width=1.1pt] table[x=G,y=arms_direct_curv1_sd_re] {additional-tree-euclidean.dat};
\addplot[violet!85!black,line width=1.1pt] table[x=G,y=arms_direct_curv2_sd_re] {additional-tree-euclidean.dat};
\nextgroupplot[title={Liouville measure},ylabel={$\delta\Delta\Scal$}]
\addplot[black,line width=1.1pt] table[x=G,y=arms_liouville_curv0_sd_re] {additional-tree-euclidean.dat};
\addplot[blue!75!black,line width=1.1pt] table[x=G,y=arms_liouville_curv1_sd_re] {additional-tree-euclidean.dat};
\addplot[violet!85!black,line width=1.1pt] table[x=G,y=arms_liouville_curv2_sd_re] {additional-tree-euclidean.dat};
\end{groupplot}
\node[anchor=center,font=\scriptsize,inner sep=0pt] at ($(extra19 c1r1.south west)!0.5!(extra19 c2r2.north east)$) {
\pgfplotslegendfromname{leg:extra-curv-arms}
};
\end{tikzpicture}
\caption{Extended 3-star: local scalar-curvature means and fluctuations.}
\label{fig:extra-curvature-arms}
\end{figure}

The Ricci scalar values are
\begin{align*}
 \Scal(0)&=-\frac4{27}\sum_{i=1}^3\frac{\sum_{j\ne i}a_j}{a_i^2}
           -\frac1{27}\sum_{i=1}^3\frac1{a_i}
           +\frac1{18}\sum_{i=1}^3\frac1{a'_i},\\
 \Scal(i)&=-\frac{3a'_i}{a_i^2}-\frac1{6a_i}
                              +\frac13\sum_{j\ne i}\frac1{a_j},
 \qquad \Scal(i')=\frac3{a_i}-\frac1{2a'_i}.
\end{align*}
Choose
\[
 \mu(0)=\frac32\sum_i a_i,\qquad \mu(i)=2a_i,\qquad
 \mu(i')=\frac{3a_i+\sum_{j\ne i}a_j}{6}.
\]
Substitution gives
\begin{equation}\label{eq:three-arm-action}
 S[a]=-\frac29\sum_{i=1}^3
      \left(\frac{\sum_{j\ne i}(a_j-a_i)}{a_i}\right)^2
      -\sum_{i=1}^3\frac{(a_i-6a'_i)^2}{6a_i a'_i}.
\end{equation}
This is nonpositive and vanishes exactly when
\(a_1=a_2=a_3=a\) and \(a'_1=a'_2=a'_3=a/6\).
All seven scalar curvature values then vanish.
For \(\iota_\sigma\), the same action is obtained by replacing
\(\mu(0)\) by \(2\mu(0)/3\) and each \(\mu(i')\) by
\(2\mu(i')\), keeping \(\mu(i)\) unchanged.

For direct or Liouville measure on the six free parameters in
\((\delta,L)\), the three outer integrations again separate.
Using \(J_\nu\) from the preceding model, the Euclidean partition
function is
\[
 Z_\nu=\int_{(\delta,L)^3}
 \exp\left[-\frac2{9G}\sum_{i=1}^3
        \left(\frac{\sum_{j\ne i}(a_j-a_i)}{a_i}\right)^2\right]
 \prod_{i=1}^3a_i^{\nu-1}J_\nu(a_i;6)\,\mathrm da_i.
\]
Dimensionless expectations depend on \(G\) and \(\delta/L\).
The direct-measure partition function and bounded metric moments
converge at \(\delta=0\); Liouville measure needs a lower cutoff
because of the undamped common scale.

For direct measure at \(\delta=0\), the five transverse Gaussian
integrations give a weak-coupling scale density proportional to
\(a^5\,\mathrm da\). The preferred ratio \(a'_i=a_i/6\) gives
\[
 \frac{\< a_i\>}{L}\longrightarrow\frac67,\qquad
 \frac{\< a'_i\>}{L}\longrightarrow\frac17,\qquad
 \frac{\Delta a_i}{\< a_i\>},
 \frac{\Delta a'_i}{\< a'_i\>}
 \longrightarrow\frac1{\sqrt{48}}\quad(G\downarrow0).
\]
At strong coupling both inner and outer metric parameters have
mean \(L/2\) and relative uncertainty \(1/\sqrt3\).

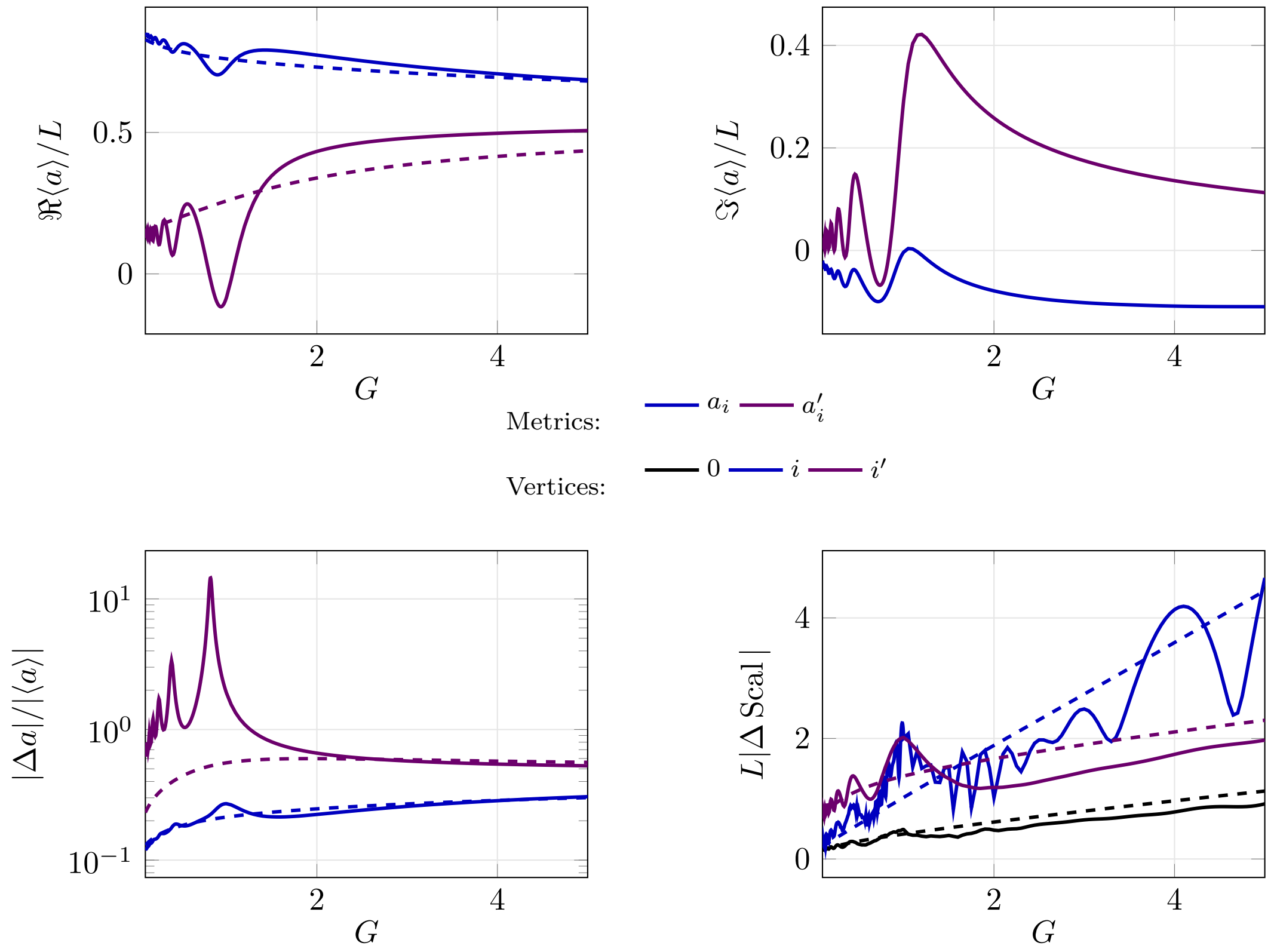
\begin{figure}[tbp]
\centering
\begin{tikzpicture}
\begin{groupplot}[
 group style={group name=extra21,group size=2 by 2,horizontal sep=2.6cm,vertical sep=2.4cm},
 width=.40\linewidth,height=5.2cm,xmin=0.1,xmax=5,xlabel={$G$},xlabel style={yshift=4pt},
 grid=major,grid style={gray!20},tick label style={font=\small},
 label style={font=\small},title style={font=\small},scaled y ticks=false,
 yticklabel style={/pgf/number format/fixed,/pgf/number format/precision=3},
 legend style={font=\scriptsize,draw=none},legend columns=-1,
]
\nextgroupplot[ylabel={$\Re\<a\>/L$},legend to name=leg:extra-lor-met-arms]
\addplot[blue!75!black,line width=1.1pt] table[x=G,y=O_arms_direct_inner_mean_re] {extended-arm-lorentzian-decade.dat};
\addlegendentry{$a_i$}
\addplot[blue!75!black,dashed,line width=1.1pt,forget plot] table[x=G,y=E_arms_direct_inner_mean_re] {extended-arm-lorentzian-decade.dat};
\addplot[violet!85!black,line width=1.1pt] table[x=G,y=O_arms_direct_outer_mean_re] {extended-arm-lorentzian-decade.dat};
\addlegendentry{$a'_i$}
\addplot[violet!85!black,dashed,line width=1.1pt,forget plot] table[x=G,y=E_arms_direct_outer_mean_re] {extended-arm-lorentzian-decade.dat};
\nextgroupplot[ylabel={$\Im\<a\>/L$}]
\addplot[blue!75!black,line width=1.1pt] table[x=G,y=O_arms_direct_inner_mean_im] {extended-arm-lorentzian-decade.dat};
\addplot[violet!85!black,line width=1.1pt] table[x=G,y=O_arms_direct_outer_mean_im] {extended-arm-lorentzian-decade.dat};
\nextgroupplot[ylabel={$|\Delta a|/|\<a\>|$},ymode=log]
\addplot[blue!75!black,line width=1.1pt] table[x=G,y=O_arms_direct_inner_relative_sd_re] {extended-arm-lorentzian-decade.dat};
\addplot[blue!75!black,dashed,line width=1.1pt,forget plot] table[x=G,y=E_arms_direct_inner_relative_sd_re] {extended-arm-lorentzian-decade.dat};
\addplot[violet!85!black,line width=1.1pt] table[x=G,y=O_arms_direct_outer_relative_sd_re] {extended-arm-lorentzian-decade.dat};
\addplot[violet!85!black,dashed,line width=1.1pt,forget plot] table[x=G,y=E_arms_direct_outer_relative_sd_re] {extended-arm-lorentzian-decade.dat};
\nextgroupplot[ylabel={$L|\Delta\Scal|$},legend to name=leg:extra-lor-curv-arms]
\addplot[black,line width=1.1pt] table[x=G,y=O_arms_direct_curv0_sd_re] {extended-arm-lorentzian-decade.dat};
\addlegendentry{$0$}
\addplot[black,dashed,line width=1.1pt,forget plot] table[x=G,y=E_arms_direct_curv0_sd_re] {extended-arm-lorentzian-decade.dat};
\addplot[blue!75!black,line width=1.1pt] table[x=G,y=O_arms_direct_curv1_sd_re] {extended-arm-lorentzian-decade.dat};
\addlegendentry{$i$}
\addplot[blue!75!black,dashed,line width=1.1pt,forget plot] table[x=G,y=E_arms_direct_curv1_sd_re] {extended-arm-lorentzian-decade.dat};
\addplot[violet!85!black,line width=1.1pt] table[x=G,y=O_arms_direct_curv2_sd_re] {extended-arm-lorentzian-decade.dat};
\addlegendentry{$i'$}
\addplot[violet!85!black,dashed,line width=1.1pt,forget plot] table[x=G,y=E_arms_direct_curv2_sd_re] {extended-arm-lorentzian-decade.dat};
\end{groupplot}
\node[anchor=center,font=\scriptsize,inner sep=0pt] at ($(extra21 c1r1.south west)!0.5!(extra21 c2r2.north east)$) {
\begin{tabular}{@{}l@{\quad}l@{}}
Metrics: & \pgfplotslegendfromname{leg:extra-lor-met-arms}\\[2pt]
Vertices: & \pgfplotslegendfromname{leg:extra-lor-curv-arms}
\end{tabular}
};
\end{tikzpicture}
\caption{Extended 3-star: Lorentzian (solid) and Euclidean (dashed) expectations.}
\label{fig:extra-lorentzian-arms}
\end{figure}

The extended-arm metric means were already included in Figure~\ref{fig:extra-metrics} for comparison. We see that it has widely separated inner and outer metric means at weak coupling, reflecting its preferred ratio
\(a'_i=a_i/6\).  Scalar curvatures are plotted in Figure~\ref{fig:extra-curvature-arms}. Unlike the adjacent-stars case, the extended-arm mean curvatures have mixed signs with direct measure.  The Lorentzian model is shown in Figure~\ref{fig:extra-lorentzian-arms} and was  regularised as in Section~\ref{sec:lore}.  We note the peak in relative uncertainty for the outer metric parameter $a'$ at $G=0.82$ for the $\delta=1/10$ case plotted, as well as interference effects.

\FloatBarrier



\begin{thebibliography}{99}
\bibitem{ArgMa:qua} J.~N.~Argota-Quiroz and S.~Majid,
Quantum gravity on polygons and \(\R\times\Z_n\) FLRW model,
\emph{Class. Quantum Grav.} \textbf{37} (2020) 245001.
\href{https://doi.org/10.1088/1361-6382/abbaa8}{doi:10.1088/1361-6382/abbaa8};
\bibitem{ArgMa:dis} J.~N.~Argota-Quiroz and S.~Majid,
Quantum Riemannian geometry of the discrete interval and
\(q\)-deformation, \emph{J. Math. Phys.} \textbf{64} (2023) 051701.
\bibitem{BegMa} E.~J.~Beggs and S.~Majid,
\emph{Quantum Riemannian Geometry}, Grundlehren der mathematischen
Wissenschaften \textbf{355}, Springer (2020).
\href{https://doi.org/10.1007/978-3-030-30294-8}{doi:10.1007/978-3-030-30294-8}.
\bibitem{BegMa:gra} E.~Beggs and S.~Majid,
Quantum geodesic flows on graphs,
\emph{Lett. Math. Phys.} \textbf{114} (2024) 112.
\href{https://doi.org/10.1007/s11005-024-01860-6}{doi:10.1007/s11005-024-01860-6}.
\bibitem{BliMa:cur} S.~Blitz and S.~Majid,
Quantum curvature fluctuations and the cosmological constant in a
single plaquette quantum gravity model,
\emph{Class. Quantum Grav.} \textbf{42} (2025) 04LT01.
\href{https://doi.org/10.1088/1361-6382/ada90c}{doi:10.1088/1361-6382/ada90c}.
\bibitem{BSZ} A.~Bochniak, A.~Sitarz and P.~Zalecki,
Riemannian geometry of a discretized circle and torus,
\emph{SIGMA} \textbf{16} (2020) 143 (28pp).
\href{https://doi.org/10.3842/SIGMA.2020.143}{doi:10.3842/SIGMA.2020.143}.
\bibitem{Car}S. Carlip, Hiding the cosmological constant, \emph{Phys. Rev. Lett.} 123 (2019) 131302
\bibitem{Con} A.~Connes,
\emph{Noncommutative Geometry}, Academic Press, San Diego (1994).
\bibitem{Dri} V.~G.~Drinfeld,
Quantum groups, in \emph{Proc. ICM}, Berkeley, 1986, Vol.~1, A.~M.~Gleason (ed.),
Amer. Math. Soc., Providence (1987) 798--820.
\bibitem{Jim} M.~Jimbo,
A \(q\)-difference analogue of \(U(\mathfrak g)\) and the Yang--Baxter equation,
\emph{Lett. Math. Phys.}  \textbf{10} (1985) 63--69.
\href{https://doi.org/10.1007/BF00704588}{doi:10.1007/BF00704588}.
\bibitem{Kac} V.~G.~Kac,
\emph{Infinite-Dimensional Lie Algebras}, Cambridge University Press (1990).
\bibitem{LirMa:fuz} E.~Lira-Torres and S.~Majid,
Quantum gravity and Riemannian geometry on the fuzzy sphere,
\emph{Lett. Math. Phys.} \textbf{111} (2021) 29.
\href{https://doi.org/10.1007/s11005-021-01359-4}{doi:10.1007/s11005-021-01359-4}.
\bibitem{LiuMa:kal} C.~Liu and S.~Majid,
Kaluza--Klein ansatz from Lorentzian quantum gravity on the fuzzy
sphere, \emph{Eur. Phys. J. C} \textbf{85} (2025) 1464.
\href{https://doi.org/10.1140/epjc/s10052-025-15160-3}{doi:10.1140/epjc/s10052-025-15160-3}.
\bibitem{Ma:book} S.~Majid,
\emph{Foundations of Quantum Group Theory}, Cambridge University Press (1995,2000)
\bibitem{Ma:gra} S.~Majid,
Noncommutative Riemannian geometry of graphs,
\emph{J. Geom. Phys.} \textbf{69} (2013) 74--93.
\bibitem{Ma:squ} S.~Majid,
Quantum gravity on a square graph,
\emph{Class. Quantum Grav.} \textbf{36} (2019) 245009.
\href{https://doi.org/10.1088/1361-6382/ab4975}{doi:10.1088/1361-6382/ab4975}.
\bibitem{Ma:par} S.~Majid,
Quantum Riemannian geometry and particle creation on the integer line,
\emph{Class. Quantum Grav.} \textbf{36} (2019) 135011.
\bibitem{Ma:alg4} S.~Majid,
Algebraic approach to quantum gravity IV: applications,
\emph{PoS} \textbf{CORFU2025} (2025)
\href{https://pos.sissa.it/509/335/}{335}.
\bibitem{Ma:yet} S.~Majid,
Quantum gravity: are we there yet?,
\emph{Phil. Trans. Roy.  Soc. A} \textbf{383} (2025) 20230377.
\href{https://doi.org/10.1098/rsta.2023.0377}{doi:10.1098/rsta.2023.0377}.
\bibitem{Smi} J.~H.~Smith,
Some properties of the spectrum of a graph,
in \emph{Combinatorial Structures and Their Applications},
R.~Guy, H.~Hanani, N.~Sauer and J.~Sch\"onheim (eds.),
Gordon and Breach, New York (1970) 403--406.
\bibitem{WanUnr}Q. Wang and W. Unruh, Vacuum fluctuation, microcyclic universes and the cosmological constant problem,  \emph{Phys. Rev. D} 99 (2019) 023537
\bibitem{Whe} J.~A.~Wheeler,
Geons, \emph{Phys. Rev.} \textbf{97} (1955) 511--536.
\href{https://doi.org/10.1103/PhysRev.97.511}{doi:10.1103/PhysRev.97.511}.
\end{thebibliography}
\end{document}